\documentclass[onecolumn,final]{elsart3p}

\usepackage{amsmath} % equations and subequations
\usepackage{amsfonts} % fonts of math equations  
\usepackage{bbm}
\usepackage{mathrsfs} % for an alternative font for categories and sheaves
\usepackage{xifthen} % for writing IfThenElse
\usepackage[usenames,dvipsnames,table]{xcolor}
\usepackage[numbers]{natbib}
\usepackage{bigints}
\usepackage{graphicx}
\usepackage{booktabs}
\usepackage{para list}
\usepackage[shortlabels]{enumitem}

\newcommand{\cref}[1]{\eqref{#1}}
\newcommand{\Cref}[1]{\eqref{#1}}

\usepackage{multicol}

\usepackage{cancel}
\usepackage{changes}
\newcommand{\soutthick}[1]{%
    \renewcommand{\ULthickness}{1.0pt}%
       \sout{#1}%
    \renewcommand{\ULthickness}{.4pt}% Resetting to ulem default
}
\newcommand{\stkout}[1]{
  \ifmmode\text{\soutthick{\ensuremath{#1}}}\else\soutthick{#1}\fi
}
\setdeletedmarkup{\stkout{#1}}

\definechangesauthor[name={Mario}, color={blue}]{MP}
\definechangesauthor[name={Gianmarco}, color={green}]{GM}
\definechangesauthor[name={Elena}, color={red}]{EB}

\usepackage{mathrsfs}
\usepackage{scalerel}

\DeclareMathAlphabet{\Mymathbb}{U}{bbold}{m}{n}
\DeclareMathAlphabet{\mathpzc}{OT1}{pzc}{m}{it}

\newcommand{\T}{\mbox{\tiny T}}

\newcommand{\Forall}{\forall \,}%{\mbox{ for all }}

\newcommand{\FirstSymbol}{\mathcal{G}}
\newcommand{\First}[1][]{
  \ifthenelse{\equal{#1}{}}
  {\FirstSymbol}
  {\FirstSymbol_{\scriptscriptstyle{#1}}}
}
\newcommand{\FirstForm}[1][]{
  \ifthenelse{\equal{#1}{}}       
  {\operatorname{I_{\point}}}             
  {\operatorname{I_{#1}}}             
}
\newcommand{\GradSymbol}{\operatorname{\mathbf{\nabla}}}
\newcommand{\Grad}{\GradSymbol}
\newcommand{\Div}{\Grad\cdot}
\newcommand{\Lap}{\Delta}
\newcommand{\GradSurf}{\GradSymbol_{\scriptscriptstyle\First}}
\newcommand{\DivSurf}{\GradSurf\cdot}

\newcommand{\LapSurf}{\Lap_{\scriptscriptstyle\First}}

\newcommand{\Der}[1][]
{
  \ifthenelse{\equal{#1}{}}
  {\partial}
  {\partial_{\scriptscriptstyle{#1}}}
}

\newcommand{\DerPar}[2]{\frac{\Der #1}{\Der #2}}

\newcommand{\DerParDS}[2]{\dfrac{\Der #1}{\Der #2}}

\newcommand{\DerTot}[2][t]
{
  \ifthenelse{\equal{#2}{}}
  {\frac{d #2}{dt}}
  {\frac{d #2}{d #1}}
}
\newcommand{\Diffsymbol}{\operatorname{d}\!}
\newcommand{\Diff}[2][]{
  \ifthenelse{\equal{#1}{}}
  {\Diffsymbol{#2}}
  {\Diffsymbol{#2}_{{#1}}}
}
\newcommand{\DirDerSymb}{D}
\newcommand{\DirDer}[2][]
{
  \ifthenelse{\equal{#1}{}}
  {\DirDerSymb^{#2}}
  {\DirDerSymb^{#2}_{#1}}
}

\newcommand{\REALsymbol}{\mathbb{R}}
\newcommand{\REAL}[1][]{
  \ifthenelse{\equal{#1}{}}
  {\REALsymbol}
  {{\REALsymbol}^{#1}}
}
\newcommand{\EUCLsymbol}{\mathbb E}
\newcommand{\EUCL}[1][]{
  \ifthenelse{\equal{#1}{}}
  {\EUCLsymbol}
  {{\EUCLsymbol}^{#1}}
}
\newcommand{\NATURALsymbol}{\mathbb N}
\newcommand{\NATURAL}[1][]{
  \ifthenelse{\equal{#1}{}}
  {\NATURALsymbol}
  {{\NATURALsymbol}^{#1}}
}

\newcommand{\ABS}[2][]
{
  \ifthenelse{\equal{#1}{}}
  {\left| #2 \right|}
  {\left| #2 \right|_{#1}}
}

\newcommand{\NORM}[2][]
{
  \ifthenelse{\equal{#1}{}}
  {\left\| #2 \right\|}
  {\left\| #2 \right\|_{#1}}
}

\newcommand{\SCAL}[3][]
{
  \ifthenelse{\equal{#1}{}}
  {\left\langle{#2},{#3}\right\rangle}
  {\left\langle{#2},{#3}\right\rangle_{#1}}
}
\newcommand{\SCALF}[3][]
{
  \ifthenelse{\equal{#1}{}}
  {\left({#2},{#3}\right)}
  {\left({#2},{#3}\right)_{#1}}
}

\newcommand{\scalprodSurf}[3][]
{
  \ifthenelse{\equal{#1}{}}
  {\left\langle {#2},{#3} \right\rangle_{\scriptscriptstyle\First}}
  {\left\langle {#2},{#3} \right\rangle_{{#1}}}
}

\newcommand{\DET}[1]{\ensuremath{\operatorname{det}#1}}

\newcommand{\Cond}{\kappa}

\newcommand{\point}[1][]
{
  \ifthenelse{\equal{#1}{}}
  {\mathbf{p}}
  {\mathbf{p}_{#1}}
}

\newcommand{\qpoint}{\mathbf{q}}

\newcommand{\midPoint}{\mathbf{m}}

\newcommand{\From}{:}
\newcommand{\To}{\rightarrow}
\newcommand{\Mapsto}{\mapsto}

\newcommand{\RegionSymb}{\Cellsymb}
\newcommand{\Region}[1][]
{
  \ifthenelse{\equal{#1}{}}
  {\RegionSymb}
  {\RegionSymb_{#1}}
}

\newcommand{\Interval}[1][]
{
  \ifthenelse{\equal{#1}{}}
  {I}
  {I_{#1}}
}

\newcommand{\SurfDomainsymb}{\Gamma}
\newcommand{\SurfDomain}[1][]{
  \ifthenelse{\equal{#1}{}}
  {\SurfDomainsymb}
  {\SurfDomainsymb_{#1}}
}
\newcommand{\SurfDomainBndsymb}{\partial\Gamma}
\newcommand{\SurfDomainBnd}[1][]{
  \ifthenelse{\equal{#1}{}}
  {\SurfDomainBndsymb}
  {\SurfDomainBndsymb_{#1}}
}
\newcommand{\ClosedSurfDomain}[1][]{
  \ifthenelse{\equal{#1}{}}
  {\Closedsymb{\SurfDomain}}
  {\Closedsymb{\SurfDomain[#1]}}
}

\newcommand{\heightsymb}{\mathcal{H}}
\newcommand{\height}[1][]{
  \ifthenelse{\equal{#1}{}}
  {\heightsymb}
  {\heightsymb_{#1}}
}

\newcommand{\BSMsymbol}{\mathcal{B}}

\newcommand{\BSM}[1][]
{
  \ifthenelse{\equal{#1}{}}
  {\BSMsymbol}
  {\BSMsymbol_{#1}}
}

\newcommand{\Surf}{\mathcal{S}}

\newcommand{\SurfBndsymb}{\partial\Surf}
\newcommand{\SurfBnd}[1][]{
  \ifthenelse{\equal{#1}{}}
  {\SurfBndsymb}
  {\SurfBndsymb_{#1}}
}
\newcommand{\Closedsymb}[1]{\bar{#1}}
\newcommand{\ClosedSurf}[1][]{
  \ifthenelse{\equal{#1}{}}
  {\Closedsymb{\Surf}}
  {\Closedsymb{\Surf[#1]}}
}

\newcommand{\Vector}[1]{\mathbf{#1}}

\newcommand{\MapUsymb}{\phi}
\newcommand{\MapU}[1][]
{
  \ifthenelse{\equal{#1}{}}
    {\MapUsymb}
    {\MapUsymb_{#1}}
}
\newcommand{\MapVsymb}{\psi}
\newcommand{\MapV}[1][]
{
  \ifthenelse{\equal{#1}{}}
    {\MapVsymb}
    {\MapVsymb_{#1}}
}
\newcommand{\Transsymb}{\Phi}
\newcommand{\Trans}[1][]
{
  \ifthenelse{\equal{#1}{}}
    {\Transsymb} 
    {\Transsymb_{\scriptscriptstyle{#1}}}
}
\newcommand{\InvMapsymb}{\Psi}
\newcommand{\InvMap}[1][]
{
  \ifthenelse{\equal{#1}{}}
    {\InvMapsymb}
    {\InvMapsymb_{\scriptscriptstyle{#1}}}
}

\newcommand{\fsymb}{f}
\newcommand{\scalFun}[1][]
{
  \ifthenelse{\equal{#1}{}}
  {\fsymb}
  {\fsymb_{#1}}
}
\newcommand{\tscalFun}[1][]
{
  \ifthenelse{\equal{#1}{}}
  {\tilde{\fsymb}}
  {\tilde{\fsymb}_{#1}}
}
\newcommand{\hscalFun}[1][]
{
  \ifthenelse{\equal{#1}{}}
  {\hat{\fsymb}}
  {\hat{\fsymb}_{#1}}
}
\newcommand{\gsymb}{g}
\newcommand{\scalFung}[1][]
{
  \ifthenelse{\equal{#1}{}}
  {\gsymb}
  {\gsymb_{#1}}
}
\newcommand{\Fsymb}{F}
\newcommand{\FvecFun}[1][]
{
  \ifthenelse{\equal{#1}{}}
  {\Fsymb}
  {\Fsymb_{#1}}
}
\newcommand{\tFvecFun}[1][]
{
  \ifthenelse{\equal{#1}{}}
  {\tilde{\Fsymb}}
  {\tilde{\Fsymb}_{#1}}
}
\newcommand{\Ffunc}[2][]
{
  \ifthenelse{\equal{#1}{}}
  {\ensuremath{\Fsymb_{#2}}}
  {\ensuremath{\Fsymb^{#1}_{#2}}}
}
\newcommand{\inverse}[1]{#1^{-1}}

\newcommand{\Jac}{\mathbf{J}}

\newcommand{\PrincipalK}[1][]
{
  \ifthenelse{\equal{#1}{}}
  {k}
  {k_{#1}}
}

\newcommand{\vecsymb}{u}
\newcommand{\vecFun}[1][]{
  \ifthenelse{\equal{#1}{}}
  {\Vector{\vecsymb}}
  {\vecsymb^{#1}}
}
\newcommand{\tvecFun}[1][]{
  \ifthenelse{\equal{#1}{}}
  {\tilde{\vecsymb}}
  {\tilde{\vecsymb}_{#1}}
}

\newcommand{\vvsymb}{v}
\newcommand{\vv}[1][]{
  \ifthenelse{\equal{#1}{}}
  {\mathbf{\vvsymb}}
  {\vvsymb^{#1}}
}
\newcommand{\wwsymb}{w}
\newcommand{\ww}[1][]
{
  \ifthenelse{\equal{#1}{}}
  {\mathbf{\wwsymb}}
  {\wwsymb^{#1}}
}
\newcommand{\uusymb}{u}
\newcommand{\uu}[1][]
{
  \ifthenelse{\equal{#1}{}}
  {\mathbf{\uusymb}}
  {\uusymb^{#1}}
}

\newcommand{\VecFieldSymbol}{X}
\newcommand{\VecField}[1][]
{
  \ifthenelse{\equal{#1}{}}
  {\VecFieldSymbol}
  {\VecFieldSymbol^{#1}}
}
\newcommand{\VecFieldYSymbol}{Y}
\newcommand{\VecFieldY}[1][]
{
  \ifthenelse{\equal{#1}{}}
  {\VecFieldYSymbol}
  {\VecFieldYSymbol^{#1}}
}

\newcommand{\xvsymb}{x}
\newcommand{\xv}[1][]
{
  \ifthenelse{\equal{#1}{}}
  {\mathbf{\xvsymb}}
  {\mathbf{\xvsymb}_{\scriptscriptstyle{#1}}}
}
\newcommand{\xvcomp}[1][]{
  \ifthenelse{\equal{#1}{}}
  {\xvsymb}
  {\xvsymb^{\scriptscriptstyle{#1}}}
}
\newcommand{\xcg}[1][]{
  \ifthenelse{\equal{#1}{}}
  {\xvcomp[1]}
  {\xvcomp[1]_{\scriptscriptstyle{#1}}}
}
\newcommand{\ycg}[1][]{
  \ifthenelse{\equal{#1}{}}
  {\xvcomp[2]}
  {\xvcomp[2]_{\scriptscriptstyle{#1}}}
}
\newcommand{\zcg}[1][]{
  \ifthenelse{\equal{#1}{}}
  {\xvcomp[3]}
  {\xvcomp[3]_{\scriptscriptstyle{#1}}}
}
\newcommand{\gcscomp}{\xcg,\ycg,\zcg}

\newcommand{\svsymb}{s}
\newcommand{\sv}[1][]{
  \ifthenelse{\equal{#1}{}}
  {\mathbf{\svsymb}}
  {\mathbf{\svsymb}_{\scriptscriptstyle{#1}}}
}
\newcommand{\svcomp}[1][]
{
  \ifthenelse{\equal{#1}{}}
  {\svsymb}
  {\svsymb^{\scriptscriptstyle{#1}}}
}
\newcommand{\xcl}[1][]{
  \ifthenelse{\equal{#1}{}}
   {\svcomp[1]}
   {\svcomp[1]_{\scriptscriptstyle{#1}}}
}
\newcommand{\ycl}[1][]{
  \ifthenelse{\equal{#1}{}}
   {\svcomp[2]}
   {\svcomp[2]_{\scriptscriptstyle{#1}}}
}
\newcommand{\zcl}[1][]{
  \ifthenelse{\equal{#1}{}}
   {\svcomp[3]}
   {\svcomp[3]_{\scriptscriptstyle{#1}}}
}

\newcommand{\ProjSymb}{\operatorname{\pi}}
\newcommand{\ProjFun}[2][]{
  \ifthenelse{\equal{#1}{}}
  {\ProjSymb\left(#2\right)}
  {\ProjSymb_{\scriptscriptstyle{#1}}\left(#2\right)}
}
\newcommand{\Prm}[1][]
{
  \ifthenelse{\equal{#1}{}}
  {\operatorname{pr}}
  {\operatorname{pr}_{\scriptscriptstyle{#1}}}
}
\newcommand{\TanPlane}[2][]
{
  \ifthenelse{\equal{#1}{}}
  {T_{\scriptscriptstyle{\point}}#2}
  {T_{\scriptscriptstyle{#1}}#2}
}

\newcommand{\SubsetSymbol}{\mathcal{U}}

\newcommand{\SubsetU}[1][]
{
  \ifthenelse{\equal{#1}{}}
  {{U}}
  {{U}_{#1}}
}
\newcommand{\SubsetV}[1][]
{
  \ifthenelse{\equal{#1}{}}
  {{V}}
  {{V}_{#1}}
}
\newcommand{\SubsetW}[1][]
{
  \ifthenelse{\equal{#1}{}}
  {{W}}
  {{W}_{#1}}
}
\newcommand{\NeighSymbol}{\mathcal{N}}
\newcommand{\Neigh}[1][]
{
  \ifthenelse{\equal{#1}{}}
  {\NeighSymbol_{\point}}
  {\NeighSymbol_{#1}}
}
\newcommand{\NeighSurf}[1][]
{
  \ifthenelse{\equal{#1}{}}
  {\SubsetSymbol_{\point}}
  {\SubsetSymbol_{#1}}
}
\newcommand{\NormSymb}{\mathbf{N}}

\newcommand{\normalSurf}[1][]
{
  \ifthenelse{\equal{#1}{}}
  {\NormSymb}
  {\NormSymb(#1)}
}
\newcommand{\normalInterp}[1][]
{
  \ifthenelse{\equal{#1}{}}
   {\tilde{\NormSymb}}
   {\tilde{\NormSymb}_{\scriptscriptstyle{#1}}}
}

\newcommand{\normalEdge}{\mathbf{\nu}} 
\newcommand{\basisCC}{t}
\newcommand{\basisGC}{e}
\newcommand{\vecBaseGC}[1][]
{
  \ifthenelse{\equal{#1}{}}
  {\mathbf{\basisGC}}
  {\mathbf{\basisGC}_{#1}}
}
\newcommand{\vecBasePhys}[1][]
{
  \ifthenelse{\equal{#1}{}}
  {\mathbf{\basisGC}}
  {\mathbf{\basisGC}_{#1}}
}

\newcommand{\vecBaseCCcv}[1][]
{
  \ifthenelse{\equal{#1}{}}
  {\mathbf{\basisCC}}
  {\mathbf{\basisCC}_{#1}}
}

\newcommand{\tvecBaseCCcv}[1][]
{
  \ifthenelse{\equal{#1}{}}
  {\tilde{\mathbf{\basisCC}}}
  {\tilde{\mathbf{\basisCC}}_{#1}}
}
\newcommand{\hvecBaseCCcv}[1][]
{
  \ifthenelse{\equal{#1}{}}
  {\hat{\mathbf{\basisCC}}}
  {\hat{\mathbf{\basisCC}}_{#1}}
}
\newcommand{\vecBaseCCctrv}[1][]
{
  \ifthenelse{\equal{#1}{}}
  {\mathbf{\basisCC}}
  {\mathbf{\basisCC}^{#1}}
}

\newcommand{\metricsymbol}{g}
\newcommand{\metrTensCv}[1]{\metricsymbol_{\scriptscriptstyle{#1}}}

\newcommand{\first}[1]{
  \IfEqCase{#1}{
    {1}{\operatorname{E}}
    {2}{\operatorname{F}}
    {3}{\operatorname{G}}
  }
  [\PackageError{first}{Undefined option to first: #1}{}]%
}
\newcommand{\SecondFormSymbol}{\ensuremath{\operatorname{II}}}
\newcommand{\SecondForm}[1][]
{
  \ifthenelse{\equal{#1}{}}
  {\SecondFormSymbol_{\point}}
  {\SecondFormSymbol_{#1}}
}
\newcommand{\second}[1]{
  \IfEqCase{#1}{
    {1}{\operatorname{e}}
    {2}{\operatorname{f}}
    {3}{\operatorname{g}}
  }
  [\PackageError{first}{Undefined option to first: #1}{}]
}
\newcommand{\WeigSymbol}{\mathcal{W}}
\newcommand{\Weig}[1][]
{
  \ifthenelse{\equal{#1}{}}
  {\WeigSymbol}
  {\WeigSymbol_{#1}}
}

\newcommand{\velSymbol}{u}
\newcommand{\vectvel}[1][]
{
   \ifthenelse{\equal{#1}{}}
   {\vec{\velSymbol}}
   {\vec{\velSymbol}(#1)}
}

\newcommand{\velcompContr}[2][i]
{
   \ifthenelse{\equal{#2}{}}
   {\velSymbol^{#1}}
   {\velSymbol^{#1}(#2)}}
\newcommand{\velcompPhys}[2][i]
{
   \ifthenelse{\equal{#2}{}}
   {\velSymbol_{(#1)}}
   {\velSymbol_{(#1)}(#2)}}
\newcommand{\velSymbolRP}{v}
\newcommand{\velcompContrRP}[2][i]
{
   \ifthenelse{\equal{#2}{}}
   {\velSymbolRP^{#1}}
   {\velSymbolRP^{#1}(#2)}}
\newcommand{\velRP}[1][]
{
  \ifthenelse{\equal{#1}{}}
  {\velSymbolRP}
  {\velSymbolRP_{#1}}
}

\newcommand{\velcompApprox}[2][i]
{
   \ifthenelse{\equal{#2}{}}
   {\velSymbol^{#1)}}
   {\velSymbol^{#1}_{(#2)}}
}

\newcommand{\VelSymbol}{U}
\newcommand{\vectVel}[1][]
{
   \ifthenelse{\equal{#1}{}}
   {\vec{\VelSymbol}}
   {\vec{\VelSymbol}(#1)}
}
\newcommand{\Velcomp}[2][i]
{
   \ifthenelse{\equal{#2}{}}
   {\VelSymbol^{#1}}
   {\VelSymbol^{#1}(#2)}
}
\newcommand{\VprimoSymbol}{\tilde{u}}
\newcommand{\Vprimo}[1][]
{
   \ifthenelse{\equal{#1}{}}
   {\VprimoSymbol}
   {\VprimoSymbol(#1)}
}
\newcommand{\VprimoComp}[2][i]
{
   \ifthenelse{\equal{#2}{}}
   {\VprimoSymbol^{#1}}
   {\VprimoSymbol^{#1}(#2)}
}
\newcommand{\ttvelSymbol}{\tilde{\Mymathbb{u}}}
\newcommand{\ttvel}[1][]
{
   \ifthenelse{\equal{#1}{}}
   {\mathbf{\ttvelSymbol}}
   {\mathbf{\ttvelSymbol}(#1)}
}
\newcommand{\ttvelComp}[2][i]
{
   \ifthenelse{\equal{#2}{}}
   {\ttvelSymbol^{#1}}
   {\ttvelSymbol^{#1}(#2)}
}

\newcommand{\MatAlphaSymbol}{\mathbb{A}}
\newcommand{\MatAlpha}[1][]{%
  \ifthenelse{\equal{#1}{}}
  {\MatAlphaSymbol}
  {\MatAlphaSymbol_{#1}}
}

\newcommand{\QSymbol}{q}
\newcommand{\Qdisch}[1][]
{
   \ifthenelse{\equal{#1}{}}
   {\mathbf{\QSymbol}}
   {\mathbf{\QSymbol}(#1)}
}
\newcommand{\Qcomp}[2][i]
{
   \ifthenelse{\equal{#2}{}}
   {\QSymbol^{#1}}
   {\QSymbol^{#1}(#2)}
}
\newcommand{\Qvect}[1][]
{
   \ifthenelse{\equal{#1}{}}
   {\mathbf{\QSymbol}}
   {\mathbf{\QSymbol}(#1)}
}
\newcommand{\FricSymbol}{f}
\newcommand{\vectFric}[1][]
{
   \ifthenelse{\equal{#1}{}}
   {\mathbf{\FricSymbol}}
   {\mathbf{\FricSymbol}_{\scriptscriptstyle{#1}}}
}
\newcommand{\Friccomp}[2][i]
{
   \ifthenelse{\equal{#2}{}}
   {\FricSymbol_{#1}}
   {\FricSymbol_{#1}(#2)}}
\newcommand{\BFsymbol}{\tau}
\newcommand{\BottomFriction}[1][]{
  \ifthenelse{\equal{#1}{}}
  {\BFsymbol_{b}}
  {\BFsymbol_{b}^{#1}}
}

\newcommand{\ProjMatSymb}{\mathbb{P}}
\newcommand{\ProjMat}[1][]{ 
  \ifthenelse{\equal{#1}{}}
  {\ProjMatSymb}
  {\ProjMatSymb_{#1}}
}
\newcommand{\IDSymbol}{\mathbb{I}}
\newcommand{\IDtens}[1][]{
  \ifthenelse{\equal{#1}{}}
  {\IDSymbol}
  {\IDSymbol(#1)}
}

\newcommand{\tensSymbol}{\mathbb{T}}
\newcommand{\tenscompSymbol}{\tau}
\newcommand{\tenscomp}[2][ij]
{
  \ifthenelse{\equal{#2}{}}
  {\tenscompSymbol^{#1}}
  {\tenscompSymbol^{#1}(#2)}
}
\newcommand{\tensrow}[2][i]
{
  \ifthenelse{\equal{#2}{}}
  {\tensSymbol^{(#1)}}
  {\tensSymbol^{(#1)}(#2)}
}
\newcommand{\TensSymbol}{\mathbf{T}}
\newcommand{\Tens}[1][]{
  \ifthenelse{\equal{#1}{}}
  {\TensSymbol}
  {\TensSymbol_{#1}}
}

\newcommand{\TensCompSymbol}{\TensSymbol}
\newcommand{\TensComp}[2][ij]
{
  \ifthenelse{\equal{#2}{}}
  {\TensCompSymbol^{#1}}
  {\TensCompSymbol^{#1}(#2)}
}

\newcommand{\tensPrimoSymbol}{\tilde{\mathbf{\tau}}}
\newcommand{\tensPrimo}[1][]{
  \ifthenelse{\equal{#1}{}}
  {\tensPrimoSymbol}
  {\tensPrimoSymbol(#1)}
}
\newcommand{\tensPrimoCompSymbol}{\tensPrimoSymbol}
\newcommand{\tensPrimoComp}[2][ij]
{
  \ifthenelse{\equal{#2}{}}
  {\tensPrimoCompSymbol^{#1}}
  {\tensPrimoCompSymbol^{#1}(#2)}
}
\newcommand{\MCxl}[1][]{\ifthenelse{\equal{#1}{}}{h_{(1)}}{h_{(1),#1}}} 
\newcommand{\MCyl}[1][]{\ifthenelse{\equal{#1}{}}{h_{(2)}}{h_{(2),#1}}} 
\newcommand{\MCzl}[1][]{\ifthenelse{\equal{#1}{}}{h_{(3)}}{h_{(3),#1}}}
\newcommand{\MPsymb}{h}
\newcommand{\Lev}{\ell}
\newcommand{\smallestH}[1][]
{
  \ifthenelse{\equal{#1}{}}
    {{l}}
    {{l}_{#1}}
}
\newcommand{\meshparam}[1][]
{
  \ifthenelse{\equal{#1}{}}
    {\MPsymb}
    {\MPsymb_{#1}}
}
\newcommand{\InradiusSymbol}{r}
\newcommand{\Inradius}[1][]
{
  \ifthenelse{\equal{#1}{}}
  {\InradiusSymbol}
  {\InradiusSymbol_{\scriptscriptstyle{#1}}}
}
\newcommand{\Tsymb}{\mathcal{T}}
\newcommand{\Triang}[1][]
{
  \ifthenelse{\equal{#1}{}}
    {\Tsymb}
    {\Tsymb_{#1}}
}
\newcommand{\TriangH}[1][]
{
  \ifthenelse{\equal{#1}{}}
    {\Tsymb_{\meshparam}}
    {\Tsymb_{#1}}
}
\newcommand{\Edgesymb}{\sigma}
\newcommand{\Edge}[1][]{
  \ifthenelse{\equal{#1}{}}
    {\Edgesymb}
    {\Edgesymb_{#1}}
}
\newcommand{\EdgeH}[1][]{
  \ifthenelse{\equal{#1}{}}
    {\Edgesymb_{\meshparam}}
    {\Edgesymb_{\meshparam,#1}}
}
\newcommand{\NEdge}[1][]{
  \ifthenelse{\equal{#1}{}}
    {N_{\Edgesymb}}
    {N_{\Edgesymb({#1})}}
}
\newcommand{\Cellsymb}{T}
\newcommand{\Cell}[1][]{
  \ifthenelse{\equal{#1}{}}
    {\Cellsymb}
    {\Cellsymb_{#1}}
}
\newcommand{\CellH}[1][]{
  \ifthenelse{\equal{#1}{}}
    {\Cellsymb_{\meshparam}}
    {\Cellsymb_{\meshparam,#1}}
}
\newcommand{\areaSymb}{\mathcal{A}}
\newcommand{\CellArea}[1][]
{
  \ifthenelse{\equal{#1}{}}
    {\areaSymb_{\Cell}}
    {\areaSymb_{#1}}
}
\newcommand{\CellHArea}[1][]
{
  \ifthenelse{\equal{#1}{}}
    {\areaSymb_{\CellH}}
    {\areaSymb_{\meshparam,#1}}
}

\newcommand{\NCell}[1][]{
  \ifthenelse{\equal{#1}{}}
    {N_{\Cellsymb}}
    {N_{\Cellsymb({#1})}}
}

\newcommand{\lengthSymb}{{\ell}}
\newcommand{\edgeLength}[1][]
{
  \ifthenelse{\equal{#1}{}}
  {\lengthSymb_{\Edge}}
  {\lengthSymb_{#1}}
}
 
\newcommand{\edgeHLength}[1][]
{
  \ifthenelse{\equal{#1}{}}
  {\lengthSymb_{\EdgeH}}
  {\lengthSymb_{\meshparam,#1}}
}

\newcommand{\Sourcesymb}{\mathbf{S}}
\newcommand{\Source}[1][]
{
  \ifthenelse{\equal{#1}{}}
    {\Sourcesymb}
    {\Sourcesymb_{#1}}
}

\newcommand{\SourceEdge}[1][]{
  \ifthenelse{\equal{#1}{}}
    {\Sourcesymb_{ij}}
    {\Sourcesymb_{#1}}
}

\newcommand{\FluxEdgesymb}{\mathbf{F}}
\newcommand{\FluxEdge}[1][]{
  \ifthenelse{\equal{#1}{}}
    {\FluxEdgesymb_{ij}}
    {\FluxEdgesymb_{#1}}
}

\newcommand{\numFlux}[1][]
{
  \ifthenelse{\equal{#1}{}}
  {\tilde{\fsymb}}
  {\tilde{\fsymb}_{#1}}
}

\newcommand{\FluxFuncNormSymbol}{\mathbf{F}}
\newcommand{\FluxFuncNorm}[1][]{
  \ifthenelse{\equal{#1}{}}
    {\FluxFuncNormSymbol^{\normalEdge}}
    {\FluxFuncNormSymbol^{\normalEdge}_{#1}}    
}

\newcommand{\JacobianSymbol}{\mathbf{A}}
\newcommand{\Jacobian}[1][]{
  \ifthenelse{\equal{#1}{}}
    {\JacobianSymbol}
    {\JacobianSymbol_{#1}}    
}
\newcommand{\EValSymbol}{\lambda}
\newcommand{\EVal}[1][]{
  \ifthenelse{\equal{#1}{}}
  {\EValSymbol}
  {\EValSymbol_{#1}}
}
\newcommand{\EVecSymbol}{\mathbf{r}}
\newcommand{\EVec}[2][]{
  \ifthenelse{\equal{#1}{}}
  {\EVecSymbol^{(#2)}}
  {\EVecSymbol^{(#2)}_{#1}}
}

\newcommand{\midPointEdge}[1][]
{
  \ifthenelse{\equal{#1}{}}
  {\midPoint_{\scriptscriptstyle\Edge}}
  {\midPoint_{\scriptscriptstyle\Edge[#1]}}
}
\newcommand{\gpPointEdgeDG}[1][]
{
  \ifthenelse{\equal{#1}{}}
  {\point_{\scriptscriptstyle\Edge}}
  {\point_{\scriptscriptstyle\Edge,#1}}
}
\newcommand{\midPointCell}[1][]
{
  \ifthenelse{\equal{#1}{}}
  {\midPoint_{\scriptscriptstyle\Cell}}
  {\midPoint_{\scriptscriptstyle\Cell[#1]}}
}

\newcommand{\Eoc}{\mbox{eoc}_{\Lev}}

\newcommand{\speedRS}[2][]
{
  \ifthenelse{\equal{#2}{}}
  {S_{#2}}
  {S_{#2}^{#1}}
}
\newcommand{\ContSymbol}{C}
\newcommand{\Cont}[1][]{
  \ifthenelse{\equal{#1}{}}
  {\ContSymbol^{0}}
  {\ContSymbol^{#1}}
}
\newcommand{\Cinf}[1][]{
  \ifthenelse{\equal{#1}{}}
  {\ContSymbol^{\infty}}
  {\ContSymbol^{\infty}(#1)}
}
\newcommand{\SobSymbol}{W}
\newcommand{\HilbSymbol}{H}
\newcommand{\Sob}[2][]{
  \ifthenelse{\equal{#1}{}}
  {\HilbSymbol^{#2}}
  {\SobSymbol^{#2,#1}}  
}
\newcommand{\Hilb}[2][]{
  \ifthenelse{\equal{#1}{}}
  {\HilbSymbol^{#2}}
  {\HilbSymbol^{#2}_{#1}}  
}

\newcommand{\LspaceSymb}{L}
\newcommand{\Lspace}[1][]{
  \ifthenelse{\equal{#1}{}}
  {\LspaceSymb^{2}}
  {\LspaceSymb^{#1}}  
}

\newcommand{\TestSpSymbol}{\mathcal{V}}
\newcommand{\TestSpace}[1][]{
  \ifthenelse{\equal{#1}{}}
  {\TestSpSymbol({\SurfDomain})}
  {\TestSpSymbol_{#1}({\SurfDomain})}
}
\newcommand{\TestSpaceEmbedded}[1][]{
  \ifthenelse{\equal{#1}{}}
  {\TestSpSymbol(\TriangH(\SurfDomain))}
  {\TestSpSymbol_{#1}(\TriangH(\SurfDomain))}
}
\newcommand{\TestSpaceIntrinsic}[1][]{
  \ifthenelse{\equal{#1}{}}
  {\TestSpSymbol(\Triang(\SurfDomain))}
  {\TestSpSymbol_{#1}(\Triang(\SurfDomain))}
}

\newcommand{\PolySymb}{\mathcal{P}}
\newcommand{\PC}[1]{\PolySymb_{#1}}

\newcommand{\AngleSymbol}{\theta}
\newcommand{\DevAngle}[1][]
{
  \ifthenelse{\equal{#1}{}}
  {\AngleSymbol}
  {\AngleSymbol_{\scriptscriptstyle{#1}}}
}
\newcommand{\relheightSymb}{\pi}
\newcommand{\relheight}[1][]
{
  \ifthenelse{\equal{#1}{}}
  {\relheightSymb_{\scriptscriptstyle{\SurfDomain}}}
  {\relheightSymb_{\scriptscriptstyle{#1}}}
}

\newcommand{\conormal}{\mathbf{\mu}}

\newcommand{\SolSymbol}{u}
\newcommand{\Sol}{\SolSymbol}

\newcommand{\force}{f}
\newcommand{\AdvVel}{\Vector{w}}

\newcommand{\MassCoef}{\gamma}

\newcommand{\BilinearStiffSymbol}{a}
\newcommand{\BilinearStiff}[3][]
{
  \ifthenelse{\equal{#1}{}}
  {\BilinearStiffSymbol(#2,#3)}    
  {\BilinearStiffSymbol_{#1}(#2,#3)}    
}
\newcommand{\BilinearAdvSymbol}{b}
\newcommand{\BilinearAdv}[3][]
{
  \ifthenelse{\equal{#1}{}}
  {\BilinearAdvSymbol(#2,#3)}    
  {\BilinearAdvSymbol_{#1}(#2,#3)}    
}
\newcommand{\BilinearMassSymbol}{m}
\newcommand{\BilinearMass}[3][]
{
  \ifthenelse{\equal{#1}{}}
  {\BilinearMassSymbol(#2,#3)}    
  {\BilinearMassSymbol_{#1}(#2,#3)}    
}
\newcommand{\BilinearReactSymbol}{c}
\newcommand{\BilinearReact}[3][]
{
  \ifthenelse{\equal{#1}{}}
  {\BilinearReactSymbol(#2,#3)}    
  {\BilinearReactSymbol_{#1}(#2,#3)}    
}

\newcommand{\RhsSymbol}{F}
\newcommand{\Rhs}[1]{\RhsSymbol(#1)}

\newcommand{\BilinearStiffSymbolGamma}{a_{\Gamma}}
\newcommand{\BilinearStiffGamma}[3][]
{
  \ifthenelse{\equal{#1}{}}
  {\BilinearStiffSymbolGamma(#2,#3)}    
  {\BilinearStiffSymbolGamma_{#1}(#2,#3)}    
}
\newcommand{\BilinearAdvSymbolGamma}{b_{\Gamma}}
\newcommand{\BilinearAdvGamma}[3][]
{
  \ifthenelse{\equal{#1}{}}
  {\BilinearAdvSymbolGamma(#2,#3)}    
  {\BilinearAdvSymbolGamma_{#1}(#2,#3)}    
}
\newcommand{\BilinearMassSymbolGamma}{m_{\Gamma}}
\newcommand{\BilinearMassGamma}[3][]
{
  \ifthenelse{\equal{#1}{}}
  {\BilinearMassSymbolGamma(#2,#3)}    
  {\BilinearMassSymbolGamma_{#1}(#2,#3)}    
}
\newcommand{\BilinearReactSymbolGamma}{c_{\Gamma}}
\newcommand{\BilinearReactGamma}[3][]
{
  \ifthenelse{\equal{#1}{}}
  {\BilinearReactSymbolGamma(#2,#3)}    
  {\BilinearReactSymbolGamma_{#1}(#2,#3)}    
}

\newcommand{\TestSymbol}{v}
\newcommand{\Test}[1][]
{
  \ifthenelse{\equal{#1}{}}
  {\TestSymbol}  % test function
  {\TestSymbol_{\scriptscriptstyle{#1}}}
}

\newcommand{\nNodes}[1][]
{
  \ifthenelse{\equal{#1}{}}
  {N^{\scriptscriptstyle{dof}}}
  {N^{\scriptscriptstyle{dof}}_{\scriptscriptstyle{#1}}}
}

\newcommand{\TestApprox}[1][]
{
  \ifthenelse{\equal{#1}{}}
  {\TestSymbol_{\scriptscriptstyle{\meshparam}}}  % test function
  {\TestSymbol_{\scriptscriptstyle{\meshparam,#1}}}
}

\newcommand{\RhsApprox}[1]{\RhsSymbol_{\meshparam}(#1)}

\newcommand{\viscositySymb}{\nu}
\newcommand{\viscosity}[1][]
{
  \ifthenelse{\equal{#1}{}}
  {\viscositySymb}
  {\viscositySymb_{\scriptscriptstyle{#1}}}
}

\newcommand{\ResidualSymbol}{R}
\newcommand{\Residual}[1][]
{
  \ifthenelse{\equal{#1}{}}
  {\ResidualSymbol}
  {\ResidualSymbol_{#1}}
}

\newcommand{\NDOFS}{N^{\tiny\textrm{DOF}}}

\usepackage{comment}
\usepackage{chapterbib}    %   Bibtex
\usepackage{color}
\usepackage{comment}
\usepackage{bm}
\usepackage{overpic}
\usepackage{times}
\usepackage{epsfig}
\usepackage{graphicx,graphics,rotating}
\usepackage{multicol}
\usepackage{multirow}

\newtheorem{assumption} {Assumption}
\newtheorem{problem}    {Problem}
\newtheorem{remark}     {Remark}
\newtheorem{definition} {Definition}
\newtheorem{proposition}{Proposition}
\newtheorem{lemma}      {Lemma}
\newtheorem{theorem}    {Theorem}
\newtheorem{acknowledgements}{Acknowledgements}

\def\trait #1 #2 #3 {\vrule width #1pt height #2pt depth #3pt}
\def\fin{\hfill
        \trait .3 5 0
        \trait 5 .3 0
        \kern-5pt
        \trait 5 5 -4.7
        \trait 0.3 5 0
\medskip}

\newcommand{\PGRAPH}[1]{\medskip\noindent\textbf{#1}.}

\newcommand{\qv}{\mathbf{q}}

\newcommand{\as}{a}
\newcommand{\bs}{b}
\newcommand{\cs}{c}

\newcommand{\fs}{f}

\newcommand{\ms}{m}

\newcommand{\qs}{q}

\newcommand{\us}{u}
\newcommand{\vs}{v}
\newcommand{\ws}{w}

\newcommand{\Cs}{C}

\newcommand{\Fs}{F}

\newcommand{\matG}{\First}%\mathsf{G}} %EB

\newcommand{\matK}{\mathsf{K}}

\newcommand{\PS}[1]{\mathbbm{P}_{#1}}

\newcommand{\HONE}  {H^1}
\newcommand{\HONEzr}{H^1_0}

\newcommand{\LTWO}  {L^2}

\newcommand{\HS}[1] {H^{#1}}

\newcommand{\CS}[1] {C^{#1}}

\newcommand{\VS}[1] {V^{#1}}

\renewcommand{\P} {E}% {\textsf{E}}            % polyhedral element
\newcommand  {\E} {e}% {\textsf{e}}            % edge
\renewcommand{\S} {\sigma} %%{\textsf{s}} % sides

\newcommand{\hh}{h}
\newcommand{\Th}{\Omega_{\hh}} %EB

\newcommand{\xvP}{\sv_{\P}}        % element
\newcommand{\hP}{\hh_{\P}}

\newcommand{\hE}{\hh_{\E}}

\newcommand{\mP}{\ABS{\P}}

\newcommand{\dV}{\,d\sv} %V} %EB
\newcommand{\dS}{\,d\S}

\newcommand{\norPE}{\mathbf{n}_{\P,\E}}

\newcommand{\ush}{\us_{\hh}}

\newcommand{\vsh}{\vs_{\hh}}

\newcommand{\wsh}{\ws_{\hh}}

\newcommand{\Fsh}{\Fs_{\hh}}

\newcommand{\asP}{\as^{\P}}

\newcommand{\ash}{\as_{\hh}}
\newcommand{\bsh}{\bs_{\hh}}
\newcommand{\csh}{\cs_{\hh}}

\newcommand{\ashP}{\as^{\P}_{\hh}}
\newcommand{\bshP}{\bs^{\P}_{\hh}}
\newcommand{\cshP}{\cs^{\P}_{\hh}}

\newcommand{\SPh} {S^{\P}_{\hh}}

\newcommand{\nlen}{\hspace{-0.2mm}}

\newcommand{\snorm}  [2]{|#1|_{#2}}

\newcommand{\norm}   [2]{|\nlen|#1|\nlen|_{#2}}

\newcommand{\abs}    [1]{|#1|}

\newcommand{\Vhks}{\VS{\hh}_{k}}
\newcommand{\Vhk}{\VS{\hh}_{k}}

\newcommand{\calM} [1]{\mathcal{M}_{#1}}

\newcommand{\Piz}[1]{\Pi^{0}_{#1}}
\newcommand{\PinP}[1]{\Pi^{\nabla,\P}_{#1}}
\newcommand{\PizP}[1]{\Pi^{0,\P}_{#1}}

\newcommand{\cbot}{c_*}
\newcommand{\ctop}{c^*}

\newcommand{\restrict}[2]{{#1}_{|{#2}}}

\newcommand{\EOD}{\end{document}}

\newcommand{\RefDomain}{\Omega} %EB

\begin{document}

\begin{frontmatter}
  
  \title{Arbitrary--order intrinsic virtual element method for
    elliptic equations on surfaces}
  
  \author[PDMATH,PDGEOS]{E.~Bachini}
  \author[IMATI]{, G.~Manzini}
  \author[PDMATH]{, and M.~Putti}
  
  \address[PDMATH]{
    Department of Mathematics ``Tullio Levi-Civita'',
    University of Padua, Italy.\\
    Now at: Institute of Scientific Computing,
    Department of Mathematics, TU Dresden, Germany\\
    \emph{elena.bachini@tu-dresden.de}
    %%\emph{e-mail: elena.bachini@unipd.it,mario.putti@unipd.it}
  }
  
  \address[PDGEOS]{
    Department of Geosciences
    University of Padua, Italy
  }
  
  \address[IMATI]{
    Istituto di Matematica Applicata e Tecnologie Informatiche - CNR
    Pavia, Italy
    %%\emph{e-mail: marco.manzini@imati.cnr.it}
  }
  
  %------------------- Abstract  -----------------------------
  \begin{abstract}
    We develop a geometrically intrinsic formulation of the
    arbitrary-order Virtual Element Method (VEM) on polygonal cells
    for the numerical solution of elliptic surface partial
    differential equations (PDEs).
    The PDE is first written in covariant form using an appropriate
    local reference system.
    The knowledge of the local parametrization allows us to consider
    the two-dimensional VEM scheme, without any explicit approximation
    of the surface geometry.
    The theoretical properties of the classical VEM are extended to
    our framework by taking into consideration the highly anisotropic
    character of the final discretization.
    These properties are extensively tested on triangular and
    polygonal meshes using a manufactured solution.
    The limitations of the scheme are verified as functions of the
    regularity of the surface and its approximation.
  \end{abstract}

  
\end{frontmatter}

%% % ================== D O C U M E N T ====================

\raggedbottom
%\input{sec1_introduction.tex}
%\input{sec2_surface_pde.tex}
%\input{sec3_cfvem.tex}
%\input{sec4_numerical.tex}
%\input{sec5_conclusions.tex}

%% SECTION 1
\section{Introduction}
\label{sec:intro}

Surface partial differential equations of elliptic and parabolic types
are often used for the simulation of diverse phenomena in many fields
of applications, for example in biology, atmospheric dynamics, and
image processing~\cite{art:Osher1988, art:Neilson2011}.
One of the main motivation that prompted this work is related to the
modeling of gravity-driven flows in earth-sciences, such as, e.g.,
flood forecasting, landslide and debris flow dynamics, avalanche
simulations~\cite{art:Bouchut2004, art:Flyer2009, art:Fent2017,
  art:Bachini2019}.
The numerical solution of surface PDEs has seen wide-spread interest
in the last few years with different approaches being proposed,
including continuous Finite Element Methods (FEM), discontinuous
Galerkin (DG), finite volumes, trace-FEMs, etc.~\cite{art:Osher1988,
  art:Dziuk1988, art:Olshanskii2009, art:Flyer2009,
  art:Antonietti2015, art:Dede2015, art:Bachini2019, art:Bachini2020}.
A recent survey of surface FEM was published in~\cite{art:Dziuk2013},
where both steady and moving surfaces are considered, the latter
finding a unifying theory in~\cite{art:Elliott2017}.

% intro vem
Although FE-based approaches are very successful in the numerical
treatment of surface PDEs,
they share the limitation that an explicit form of the basis functions
is required in the formulation of the method, and thus are restricted
mostly to triangular/quadrilateral elements.
This restriction is overcome by the Virtual Element Method (VEM) that
was designed from the very beginning to work on generally shaped
elements with high order of accuracy.
In fact, in the VEM approach $(i)$ it is possible to decompose the
computational domain into very general polygonal elements; $(ii)$ an
explicit form of the basis functions is not required; $(iii)$
approximation of arbitrary order and arbitrary regularity are
straightforward in two and three dimensions.
%%
% The VEM formulation and its practical implementation is based on
% suitable polynomial projections that are always computable from a
% careful choice of the degrees of freedom.
%%
%
VEM was originally developed as a variational reformulation of the
\emph{nodal} mimetic finite difference (MFD) method~\cite{%
Brezzi-Buffa-Lipnikov:2009,%
BeiraodaVeiga-Lipnikov-Manzini:2011,% 
BeiraodaVeiga-Manzini-Putti:2015,%
Manzini-Lipnikov-Moulton-Shashkov:2017% 
} for solving diffusion
problems on unstructured polygonal meshes.
A survey on the MFD method can be found in the review
paper~\cite{Lipnikov-Manzini-Shashkov:2014} and the research
monograph~\cite{BeiraodaVeiga-Lipnikov-Manzini:2014}.
The scheme inherits the flexibility of the MFD method with respect to
the admissible meshes and this feature is well reflected in the many
significant applications that have been developed so far, see, for
example,~\cite{%
  BeiraodaVeiga-Manzini:2014,%
  Benedetto-Berrone-Pieraccini-Scialo:2014,%
  BeiraodaVeiga-Manzini:2015, Berrone-Pieraccini-Scialo-Vicini:2015,%
  Mora-Rivera-Rodriguez:2015,%
  Paulino-Gain:2015,%
  Antonietti-BeiraodaVeiga-Scacchi-Verani:2016,%
  BeiraodaVeiga-Chernov-Mascotto-Russo:2016,%
  BeiraodaVeiga-Brezzi-Marini-Russo:2016a,%
  BeiraodaVeiga-Brezzi-Marini-Russo:2016b,%
  Cangiani-Georgoulis-Pryer-Sutton:2016,%
  Perugia-Pietra-Russo:2016,%
  Wriggers-Rust-Reddy:2016,
  Certik-Gardini-Manzini-Vacca:2018:ApplMath:journal,%%
  Dassi-Mascotto:2018,%
  Benvenuti-Chiozzi-Manzini-Sukumar:2019:CMAME:journal,%
  Antonietti-Manzini-Verani:2019:CAMWA:journal,%%
  Certik-Gardini-Manzini-Mascotto-Vacca:2020%
}.
Because of its origins, VEM is intimately connected with other
FE-based approaches.
The connection between the VEM and finite elements on
polygonal/polyhedral meshes is thoroughly investigated
in~\cite{Manzini-Russo-Sukumar:2014,
  Cangiani-Manzini-Russo-Sukumar:2015, DiPietro-Droniou-Manzini:2018},
between VEM and discontinuous skeletal gradient discretizations
in~\cite{DiPietro-Droniou-Manzini:2018}, and between the
VEM and the BEM-based FEM method
in~\cite{Cangiani-Gyrya-Manzini-Sutton:2017:GBC:chbook}.
VEM was originally formulated
in~\cite{BeiraodaVeiga-Brezzi-Cangiani-Manzini-Marini-Russo:2013} as a
conforming FEM for the Poisson problem, and was later extended to
convection-reaction-diffusion problems with variable coefficients
in~\cite{BeiraodaVeiga-Brezzi-Marini-Russo:2016b}.
However, the VEM technology has seen so far very few applications to
surface PDEs, and only with first-order polynomial
accuracy~\cite{art:Frittelli2018}.

One of the major difficulties in the high-order numerical solution of
surface PDEs is the achievement of a consistent approximation of both
the geometry and the PDE.
The work of~\cite{art:Demlow2009} develops a general technique for
high-order polygonal approximation of a smooth manifold, but this
method requires the explicit knowledge of the distance function within
the tubular neighborhood of the surface.
The task of approximating this distance function to high order is
still an open problem~\cite{art:Memoli2005}.
A recent approach based on this idea was presented
in~\cite{art:Antonietti2015}, where the authors study a high-order (up
to four) DG scheme based on a piecewise polynomial approximation of
the surface triangulation.
However, extensions to polygonal grids with high polynomial orders
have not yet been addressed.
While VEM provides an ideal framework to work at high-order on
generally shaped cells, according to \cite{art:Frittelli2018} the main
difficulty is the high-order approximation of the surface, limiting
their current developments to polynomials of order one.
The same authors suggest the use of the approach
in~\cite{art:Demlow2009} to extend their VEM scheme to higher order
polynomials, without however eliminating the difficulty of the
approximation to a consistent order of the distance function.

In this paper, we develop a novel VEM-based approach for the solution
of elliptic surface PDEs that works at all polynomial orders. We avoid
the difficulties related to high-order surface approximation by
employing intrinsic geometry and following the approach described
in~\cite{art:Bachini2020} to adapt the virtual element technology to
the surface PDE. Using this approach, we first rewrite the partial
differential equation in covariant form in such a way that the
geometric information, essentially the metric tensor, is completely
encoded in the equation itself.
As a consequence, the numerical scheme can be constructed directly on the
two-dimensional local chart where the surface parametrization is
defined, thus enabling the full exploitation of the VEM machinery.
Here, we restrict our attention to the case when the surface is
defined by a single chart, a case of great interest for example in
gravity-driven flows on terrain surfaces, such as water
flow and sediment transport in mountain areas~\cite{%
  art:Bouchut2003,%
  art:Bouchut2004,%
  art:FernandezNieto2008,%
  art:Moretti2015,%
  art:Fent2017}.
In principle, our proposed approach can be applied to the more general
situation of a surface defined by an atlas if the transition between
charts is done with care by enforcing proper smoothness as described
in~\cite{art:Lindblom2016}.  This is shown by testing our approach on
the sphere by using the well-known charts arising from the
stereographical projection.

Our method starts from a partition of the surface into polygons with
curvilinear edges, assuming that the parametrization of the surface is
known at relevant quadrature points.
Proceeding from the covariant PDE, we construct a high-order scheme
exploiting the ability of the VEM approach to discretize problems that
are anisotropic and with spatially variable
coefficients~\cite{BeiraodaVeiga-Brezzi-Marini-Russo:2016b}.
In practice, we re-define the PDE on a local coordinate system using
intrinsic geometric quantities and operators, which contain explicitly
the metric information deriving from the surface.
Then, all the VEM projection operators are calculated using this local
coordinate system and the knowledge of the parametrization is used to
define the needed quantities, thus incurring in no explicit geometric
error.
Hence the final scheme is defined on a planar two-dimensional domain
(the surface chart) and all the available machinery to achieve
high-order on polygonal cells can be exploited.
The price we pay is that now the PDE contains the anisotropic metric
tensor and all the coefficients vary in space as a function of the
regularity of the surface.
The virtual element method has proved its efficiency in handling these
situations~\cite{Mazzia:2020} and can be implemented directly in this
two-dimensional setting.
In addition, with our approach the convergence theory extends
straight-forwardly to surface problems without additional efforts.

Our development of the intrinsic VEM proceeds as follows.
In Section~\ref{sec:geom}, we describe the local reference system of
choice and the geometric setting.
Next, we define the differential operators and the corresponding PDE
in covariant form.
In Section~\ref{sec:VEM}, we summarize the VEM adopted in this work
and discuss the necessary adaptations to the problem at hand.
The final Section~\ref{sec:numerical-results} reports the results of
extensive numerical experiments that assess the effectiveness,
accuracy, and robustness of the proposed approach.

%% SECTION 2
\section{The surface partial differential equation and its Galerkin discretization}
\label{sec:geom}

\PGRAPH{Notation} Throughout the paper, we use the standard definition
and notation of Sobolev spaces, norms and seminorms
(see~\cite{Adams-Fournier:2003}), which can be directly extended to a
compact manifold $\SurfDomain$ (see~\cite{book:Taylor2010I}). Given
$\omega$ an open and bounded subset of $\REAL^d$, $d=2,3$, we denote
with $\Lspace[p](\omega)$ and $\Sob[p]{k}(\omega)$ the Lebesgue and
Sobolev spaces, with $\Sob[2]{k}(\omega)=\Hilb{k}(\omega)$ the
classical Hilbert space.
Norms and seminorms in $\HS{k}(\omega)$ are denoted by
$\norm{\cdot}{\HS{k}(\omega)}$ and $\snorm{\cdot}{\HS{k}(\omega)}$, respectively,
and $(\cdot,\cdot)_{\omega}$ denotes the inner product in
$\LTWO(\omega)$.
We omit the subscript in the inner product notation when
$\omega$ is the whole computational domain.
In a few situations, for the sake of clarity, we may prefer to use the
integral notation of the inner product.

\bigskip
\noindent
Consider a compact surface $\SurfDomain\subset\REAL[3]$ with boundary
$\partial\SurfDomain$ over which the following elliptic partial
differential equation is defined:
\begin{align} 
  -\LapSurf\Sol +\scalprodSurf{\AdvVel}{\GradSurf\Sol}+\MassCoef\;
  \Sol &= \force\phantom{0} \qquad\mbox{on~}\SurfDomain,\nonumber\\[-0.7em]
  \label{eq:pde_surf}\\[-0.7em]
    \Sol &= 0\phantom{\force} \qquad\mbox{on~}\SurfDomainBnd,\nonumber
\end{align}
where the solution $\Sol\From\SurfDomain\rightarrow\REAL$ is a scalar
function defined on the surface,
$\AdvVel:\SurfDomain\rightarrow\REAL[2]$ is a
given divergence-free velocity field tangent to the surface, the
function $\MassCoef:\SurfDomain\to\REAL$ is a non-negative reaction
coefficient.
We denote by $\LapSurf$ and $\GradSurf$ the Laplace-Beltrami and the
tangential gradient operators, respectively, and by
$\scalprodSurf{\cdot}{\cdot}$ the intrinsic scalar product.
These operators will be given precise definition depending on the
chosen coordinate system.
Classically, we assume $\force\in \HS{-1}(\SurfDomain)$,
$\AdvVel\in[\Sob[\infty]{1}(\SurfDomain)]^2$, and
$\MassCoef-\tfrac{1}{2}\DivSurf\AdvVel>0$.
Here we consider homogeneous Dirichlet problems with the more
general boundary conditions described in \cite{art:Burman2018}.

\smallskip
\noindent
The variational formulation of equation~\eqref{eq:pde_surf} reads:

\smallskip
\begin{problem}[Intrinsic variational formulation]
  \label{pb:varformSurf}

  Find $\Sol\in\Hilb[0]{1}(\SurfDomain)$ such that
  \begin{align}
    \BilinearStiff{\Sol}{\Test} + \BilinearAdv{\Sol}{\Test} + \BilinearReact{\Sol}{\Test} = \Rhs{\Test} \quad \Forall\Test\in\Hilb[0]{1}(\SurfDomain),
    \label{eq:varform:surface}
  \end{align}
  where the bilinear forms
  $\BilinearStiff{\cdot}{\cdot},\,\BilinearAdv{\cdot}{\cdot},\,\BilinearReact{\cdot}{\cdot}:\,\Hilb[]{1}(\SurfDomain)\times\Hilb[]{1}(\SurfDomain)\to\REAL$
  are given by
  \begin{align*}
    \BilinearStiff{\Sol}{\Test}:=
    \;
    \int_{\SurfDomain}
    \scalprodSurf{\GradSurf\Sol}{\GradSurf\Test},
    \qquad
    \BilinearAdv{\Sol}{\Test} :=
    \int_{\SurfDomain}\scalprodSurf{\AdvVel}{\GradSurf\Sol}\Test,
    \qquad
    \BilinearReact{\Sol}{\Test} :=
    \int_{\SurfDomain}\MassCoef\,\Sol\,\Test,
  \end{align*}
  and the right-hand side linear functional
  $\Rhs{\cdot}:\,\HONEzr(\SurfDomain)\to\REAL$ is given by
  \begin{equation*}
    \Rhs{\Test} := \int_{\SurfDomain}\force\,\Test.
  \end{equation*}

  %% %%
  %% \begin{align*}
  %%   \BilinearStiff{\Sol}{\Test} &:= \int_{\SurfDomain}\scalprodSurf{\GradSurf\Sol}{\GradSurf\Test},\\[0.5em]
  %%   \BilinearAdv{\Sol}{\Test}   &:= \int_{\SurfDomain}\scalprodSurf{\AdvVel}{\GradSurf\Sol}\Test,\\[0.5em]
  %%   \BilinearReact{\Sol}{\Test} &:= \int_{\SurfDomain}\MassCoef\,\Sol\,\Test,
  %% \end{align*}
  %% %%
  %% and the linear functional
  %% $Rhs:\,\Hilb[]{1}(\SurfDomain)\to\REAL$ as
  %% %%
  %% \begin{equation*}
  %%   \Rhs{\Test}:=\int_{\SurfDomain}\force\,\Test.
  %% \end{equation*}
\end{problem}

\smallskip
\begin{remark}\label{rem:well-pos}
  The wellposedness of Problem~\ref{pb:varformSurf} follows from the
  application of the Lax-Milgram theorem since the classical theory of
  elliptic equations can be extended to surface PDEs in a
  straight-forward manner.
  In particular, due to the coercivity of the bilinear form
  $\as(\cdot,\cdot)$, the continuity of the bilinear forms
  $\as(\cdot,\cdot)$, $\bs(\cdot,\cdot)$, and $\cs(\cdot,\cdot)$ and
  the linear functional $\Fs(\cdot)$, and under the assumption that
  $\MassCoef-\tfrac{1}{2}\DivSurf\AdvVel>0$, the solution $\Sol$
  exists and is unique and belongs to $\HONE_0(\Gamma)$ if
  $\fs\in\HS{-1}(\Gamma)$.
  %% with the following stability property:
  %%  
  %% \begin{align*}
  %%    \norm{\Sol}{\HONE(\Gamma)}\leq \Cs\norm{\fs}{\LTWO(\Gamma)}.
  %% \end{align*}
  %%
  %%\begin{align*}
  %%    \norm{\Sol}{\RED{\HTWO(\Gamma)???}}\leq\Cs\norm{\fs}{\LTWO(\Gamma)}.
  %%\end{align*}
  %%
  For a detailed description of the properties, well-posedness, and
  regularity of the variational problem on manifolds we refer
  to~\cite{book:Gilbarg2001} and
  to~\cite{book:Taylor2010I,book:Taylor2010II}.
\end{remark}

\smallskip
\noindent
The discrete approximation of this problem reads as follows:

\smallskip
\begin{problem}[Intrinsic discrete Galerkin approximation]
  \label{problem:VEM}
  
  Find $\ush\in\Vhk$ such that
  \begin{align}
    \ash(\ush,\vsh) + \bsh(\ush,\vsh) + \csh(\ush,\vsh) = \Fsh(\vsh)
    \qquad\forall\vsh\in\Vhk,
    \label{eq:VEM}
  \end{align}
  where $\Vhk$ is the functional space that provides a conforming
  approximation of $\Hilb[0]{1}(\SurfDomain)$ in the virtual element
  setting, and $\ush$, $\ash(\cdot,\cdot)$, $\bsh(\cdot,\cdot)$,
  $\csh(\cdot,\cdot)$, and $\Fsh(\cdot)$ are the virtual element
  approximations to $\Sol$, $\as(\cdot,\cdot)$, $\bs(\cdot,\cdot)$,
  $\cs(\cdot,\cdot)$ and $\Fs(\cdot)$.
\end{problem}

\smallskip
These mathematical objects are defined and discussed in the next
sections.

\subsection{Geometrical setting.}

We assume that the surface $\SurfDomain$ is $\Cont[m]$ regular, i.e.:

\smallskip
\begin{definition}[Regular Surface]
  \label{def:reg-surface}
  A connected set $\SurfDomain\subset\REAL[3]$ is a $\Cont[m]$
  \emph{regular or embedded} surface if for all $\point\in\SurfDomain$
  there exists an open subset $\SubsetU\subseteq\REAL[2]$ and a map
  $\MapU[\point]:\SubsetU\To\REAL[3]$ of class $\Cont[m]$,
  $m\in\NATURAL\cup\{\infty\}$, such that:
  \begin{inparaenum}[i)]
  \item $\MapU[\point](\SubsetU)\subseteq\SurfDomain$ is an open
    neighborhood of $\point\in\SurfDomain$;
  \item $\MapU[\point]$ is a homeomorphism with its image (i.e., there
    exists an open neighborhood of $\point$,
    $\SubsetV\subseteq\REAL[3]$ such that
    $\MapU[\point](\SubsetU)=\SubsetV\cap\SurfDomain$);
  \item the differential $\Diff[\point]{\MapU}:\REAL[2]\To\REAL[3]$ is
    injective in $\SubsetU$ (i.e., it has maximum rank, in our case
    2).
  \end{inparaenum}
\end{definition}

\smallskip
The map $\MapU[\point]$ is the \emph{local parametrization} of
$\SurfDomain$ centered in $\point$ and we denote with
$\MapU[\point]^{-1}:\SubsetV\cap\SurfDomain\To\SubsetU$
its inverse map, called the \emph{local chart}, in $\point$.
The set $\MapU[\point](\SubsetU)\subset\SurfDomain$ is called a
\emph{coordinate neighborhood}, while
$(\xcl[\qpoint],\ycl[\qpoint])$ are the
\emph{local coordinates} of any point
$\qpoint\in\MapU[\point](\SubsetU)$.

\smallskip
\begin{remark}\label{rem:onechart}
  Throughout the paper we assume that $\SurfDomain$ is contained in
  only one chart.
  This is not a limitation under the assumption of $\Cont[\infty]$
  regularity of $\SurfDomain$ (or $\Cont[m]$ with $m$ sufficiently
  large) since we can always find compatible local parametrizations
  covering $\SurfDomain$.
  Indeed, given two points $\point$ and $\qpoint\in\SurfDomain$ with
  local parametrizations $\MapU[\point]$ and $\MapU[\qpoint]$ such
  that $\SubsetU_{\point}\cap\SubsetU_{\qpoint}\ne\emptyset$, the
  transition map $\MapU[\point]\circ\MapU[\qpoint]^{-1}$ is a
  $\Cont[\infty]$ (or $\Cont[m]$) diffeomorphism.
  Thus, it is always possible to find an atlas for $\SurfDomain$
  formed by appropriate charts that maintains all the required
  continuity properties.
  A proper selection of these charts is fundamental to obtain a
  numerically well-conditioned reference system in our approach.
  For an example of a constructive methodology for the definition of
  smooth multi-charts see~\cite{art:Lindblom2016}.
\end{remark}

\smallskip
For simplicity, from now on we will drop the subscripts $\point$ and
$\qpoint$ in both the local coordinates $\sv=(\xcl,\ycl)$ and the
global Cartesian coordinates $\xv=(\gcscomp)$.
In summary, we have the following explicit definitions of these
transformations:
\begin{align*}
  \MapU\From\SubsetU
  &  \To\SubsetV\cap\SurfDomain
  &  \MapU^{-1}\From\SubsetV\cap\SurfDomain&\To\SubsetU\\
  \sv & \Mapsto \xv
  &  
  \xv& \Mapsto \sv
\end{align*}
We want to choose a coordinate system to give a workable meaning to
the partial differential equation and related differential
operators.
We define the local reference system following the approach
in~\cite{art:Bachini2019,art:Bachini2020}.
To this aim, we compute the pair of tangent vectors
$\{\hvecBaseCCcv[1](\point),\hvecBaseCCcv[2](\point)\}$ on the tangent
plane $\TanPlane{\SurfDomain}$:
\begin{equation*}
  \hvecBaseCCcv[i] (\point) = 
%  \Diff[\point]{\MapU}(\vecBaseGC[j] (\point)) = 
  \left(
    \DerPar{\xcg}{\svcomp[i]}, \DerPar{\ycg}{\svcomp[i]},
    \DerPar{\zcg}{\svcomp[i]} 
  \right), \qquad i = 1,2 . %\mbox{ and } j=1,2,3,
\end{equation*} 
This pair is orthogonalized via Gram-Schmidt, yielding the orthogonal
frame $\{\vecBaseCCcv[1],\vecBaseCCcv[2]\}$.
The ensuing metric tensor is given by:
\begin{equation}
  \label{eq:first-FF}
  \First :=
  \left(
    \begin{array}{ccc}
      \NORM{\vecBaseCCcv[1](\point)}^2 & 0 \\
      0 &\NORM{\vecBaseCCcv[2](\point)}^2\\
    \end{array}
  \right)
  % =
  % \left(
  %   \begin{array}{ccc}
  %     \metrcoef{1}^2&     0       \\
  %     0             &\metrcoef{2}^2\\
  %   \end{array}
  % \right)
  .
\end{equation}
The associated scalar product between two vectors $\uu$ and $\vv$ is given by
$\scalprodSurf{\uu}{\vv}={\uu}\cdot{\First\,\vv}={\First\,\uu}\cdot{\vv}$, where
``$\cdot$'' is the canonical $\REAL^2$ scalar product.
Tensor $\First$ represents the realization of the first
fundamental form with respect to the chosen reference system
(chart).
For a $\Cont[m]$-regular surface (see
  definition~\ref{def:reg-surface}), the determinant $\DET{\First}$ of
  the metric tensor is a well-defined and bounded function, and the
  metric tensor itself is coercive, i.e., it is symmetric and
  positive-definite and has a symmetric and positive-definite inverse.
In other words, we can find constants $g_*$ and $g^*$ such that
\cite{book:Ciarlet2013}:
\begin{equation}
  \label{eq:G-coerc}
  g_*\NORM{\uu}^2\le\SCAL{\uu}{\First\,\uu}\le g^*\NORM{\uu}^2,
\end{equation}
where $\NORM{\uu}^2=\uu\cdot\uu$.

We can now write the intrinsic differential operators with respect to
the local coordinate system, and we collect the appropriate
definitions in the following proposition, which we state without
proof:

\smallskip
\begin{proposition}[Intrinsic Differential Operators]
  \label{def:LCS-ops}
  Given $\scalFun\From\SurfDomain\To\REAL$ a scalar differentiable
  function on $\SurfDomain$ and denoting with $\Grad$ and $\Div$ the
  gradient and divergence operators in $\REAL^2$, the
  \emph{intrinsic} differential operators expressed in the local
  coordinate system are given by the following expressions:
  \begin{itemize}
  \item The \emph{intrinsic} gradient of $\scalFun$ is:
    \begin{equation}
      \label{eq:gradf}
      \GradSurf\scalFun = 
      \First^{-1}\Grad\scalFun \, .
    \end{equation}    
  \item The \emph{intrinsic} Laplace-Beltrami operator of $\scalFun$
    is:
    \begin{equation}
      \label{eq:lapf}
      \LapSurf\scalFun =
      \DivSurf\GradSurf\scalFun=
      \frac{1}{\sqrt{\DET{\First}}} 
      \Div
      \left(
      \sqrt{\DET{\First}}\,\First^{-1}
      \Grad\scalFun
      \right).
    \end{equation}
  \end{itemize}
\end{proposition}

\smallskip
We would like to recall that our reference frame is covariant and thus
scalar products must act on vectors written in contravariant
components.
This applies both to velocity vector and to the divergence operator as
well.

\smallskip
\begin{remark}
  \label{rem:cartesian}
  We can use the orthonormal reference frame
  $\vecBaseGC[1],\vecBaseGC[2]$ as a base for $\REAL^2$ with which we
  can express differential operators and vector quantities.
  In this case we need to take into consideration the Jacobian matrix
  $\Jac=[\vecBaseCCcv[1],\vecBaseCCcv[2]]$ of the parametrization, and
  recall that $\First = \Jac^{\T}\Jac$.
  This applies in particular to the velocity vector field $\AdvVel$,
  which can be written as:
  \begin{equation*}
    \AdvVel = [w^1,w^2]^{\T}
    =\First^{-1/2}[w_{(1)},w_{(2)}]^{\T}=\First^{-1/2}\,\hat{\AdvVel}\,,
  \end{equation*}
  where $\hat{\AdvVel}=[w_{(1)}, w_{(2)}]^T$ is the velocity vector written with respect to
  $\vecBaseGC[1],\vecBaseGC[2]$.
\end{remark}

\smallskip
In this setting, we can give the definition of the integral of a
function over a surface as follows:

\smallskip
\begin{definition}
  Let $\scalFun:\SurfDomain\to\REAL$ be a continuous function defined
  on a regular surface $\SurfDomain$, which we assume contained in the
  image of a local parametrization $\MapU:\SubsetU\to\SurfDomain$.
  The \emph{integral of $\scalFun$ on $\SurfDomain$} is given by
  \begin{equation*}
    \int_{\SurfDomain}\scalFun =
    \int_{\inverse{\MapU}(\SurfDomain)}(\scalFun\circ\MapU)\,\sqrt{\DET{\First}}\Diff\sv
    \; .
  \end{equation*}
\end{definition}

\smallskip
%\begin{remark}
We can relate any function $\scalFun:\SurfDomain\to\REAL$ to a
specific coordinate system using the above coordinate transformations,
i.e.:
\begin{equation*}
  \scalFun(\xv)=\scalFun\circ\MapU(\sv) = \hscalFun(\sv)\;.
\end{equation*}
In the following we will make use only of the local coordinate system
and will write $\scalFun(\sv)$ omitting the hat symbol.
%\end{remark}

The classical tools deriving from Stokes theorems hold with the
intrinsic operators without any modification.
In particular, the intrinsic Green formula can be stated as in the
following lemma.

\smallskip
\begin{lemma}[Intrinsic Green formula]
  Let $\SurfDomain\subset\REAL[3]$ be a surface with smooth boundary
  $\SurfDomainBnd$ and given two functions
  $\Sol\in\Cont[2](\SurfDomain)$ and $\Test\in\Sob{1}(\SurfDomain)$,
  then:
  \begin{equation}
    \label{eq:intrinsic-green-lemma}
    \int_{\SurfDomain}\scalprodSurf{\GradSurf\Sol}{\GradSurf\Test}=
    -\int_{\SurfDomain}\LapSurf\Sol \;\Test
    + \int_{\SurfDomainBnd}
    \scalprodSurf{\GradSurf\Sol}{\conormal}\Test
    \;,
  \end{equation}
  where $\conormal\From\SurfDomain\To\REAL[2]$ denotes the vector
  tangent to $\SurfDomain$ and normal to $\SurfDomainBnd$ with
  components written with respect to the local reference frame
  (i.e. $\conormal=\mu^1\vecBaseCCcv[1]+\mu^2\vecBaseCCcv[2]$).
\end{lemma}

\smallskip
In view of remark~\ref{rem:cartesian}, we reformulate the bilinear
forms $\BilinearStiff{\cdot}{\cdot}$, $\BilinearAdv{\cdot}{\cdot}$,
and $\BilinearReact{\cdot}{\cdot}$ the linear functional $\Rhs{\cdot}$
of the intrinsic variational formulation \eqref{eq:varform:surface} on
the chart $\MapU^{-1}(\SurfDomain)$ through:
\begin{align*}
  \BilinearStiff{\Sol}{\Test} = \int_{\MapU^{-1}(\SurfDomain)}\left(\sqrt{\DET{\First}}\,\First^{-1}\right)\Grad\Sol\cdot\Grad\Test\Diff\sv\,,\qquad%\\[0.5em]
  \BilinearAdv  {\Sol}{\Test} = \int_{\MapU^{-1}(\SurfDomain)}\left(\sqrt{\DET{\First}}\,\First^{-1/2}\right)\hat{\AdvVel}\cdot\Grad\Sol\;\Test\Diff\sv\,,%\\[0.5em]
%\intertext{and}
\end{align*}
and
\begin{align*}
\BilinearReact{\Sol}{\Test} =
  \int_{\MapU^{-1}(\SurfDomain)}\sqrt{\DET{\First}}\,\MassCoef\,\Sol\,\Test\Diff\sv\,,\qquad\quad
  \Rhs{\Test}= \int_{\MapU^{-1}(\SurfDomain)}\sqrt{\DET{\First}}\force\,\Test\Diff\sv\,.
\end{align*}
Therefore, the intrinsic variational formulation
\eqref{eq:varform:surface} is equivalent to solving the
advection-diffusion-reaction equation in variational form:
\begin{equation}
  \label{eq:final-varform}
  \int_{\MapU^{-1}(\SurfDomain)} \Big(\matK\Grad\Sol\cdot\Grad\Test+\widetilde{\AdvVel}\cdot\Grad\Sol+\widetilde{\gamma}\Sol\Test\Big)\Diff\sv
  = \int_{\MapU^{-1}(\SurfDomain)}\widetilde{\scalFun}\Test\Diff\sv,
\end{equation}
where the equation coefficients are defined by
\begin{align*}
  \matK                = \sqrt{\DET{\First}}\,\First^{-1}\,, \qquad
  \widetilde{\AdvVel}  = \left(\sqrt{\DET{\First}}\,\First^{-1/2}\right)\hat{\AdvVel}\,,\qquad
  \widetilde{\gamma}   = \sqrt{\DET{\First}}\,\gamma\,,\qquad
  \widetilde{\scalFun} = \sqrt{\DET{\First}}\,\scalFun.
\end{align*}
The problem as above formulated is still well-posed and maintains all
the properties listed in remark~\ref{rem:well-pos}. Indeed, since
%\sout{$\sqrt{\DET{\First}}\le C_{\First}$}, 
our surface is assumed to be regular, 
there exist two positive constants $c_{_{\First}}$ and $C_{_{\First}}$ 
such that ${c_{_{\First}}\le\sqrt{\DET{\First}}\le C_{_{\First}}}$ and
\begin{equation}
  \kappa_*\ABS[\HONE(\SurfDomain)]{\Sol}^2\le\BilinearStiff{\Sol}{\Sol}\le\kappa^*\ABS[\HONE(\SurfDomain)]{\Sol}^2,
  \label{eq:coercivity-continuity}
\end{equation}
where $\kappa_*=c_{_{\First}}g^*$ and $\kappa^*=C_{_{\First}}g_*$,
and $g_*$ and $g^*$ are the constants introduced
in~\eqref{eq:G-coerc}.
The coercivity of $\as(\cdot,\cdot)$ with respect to the
$\HONE$-norm follows immediately by noting that
$\ABS[\HONE(\SurfDomain)]{\Sol}\leq\NORM[\HONE(\SurfDomain)]{\Sol}$  and applying the Poincar\'e inequality in
$\Hilb[0]{1}(\SurfDomain)$.
%     where \sout{$\kappa_*=C_{_{\First}}g_*$} \GREEN{$\kappa_*=c_{_{\First}}g^*$}
%   and \sout{$\kappa^*=C_{\First}g^*$} \GREEN{$\kappa^*=C_{_{\First}}g_*$}.
% 
Moreover, there exist two positive constants $\ws_{\max}$ and
$\gamma_{\max}$ such that
$\NORM[\infty]{\widetilde{\AdvVel}}\le
C_{_{\First}}\NORM[\infty]{\AdvVel}\le \ws_{\max}$ and
$\widetilde{\gamma}\le
C_{_{\First}}\,\gamma\le \gamma_{\max}$.

%% SECTION 3
\section{The virtual element method}
\label{sec:VEM}

\begin{figure}
  \centerline{
    \includegraphics[width=0.5\textwidth]{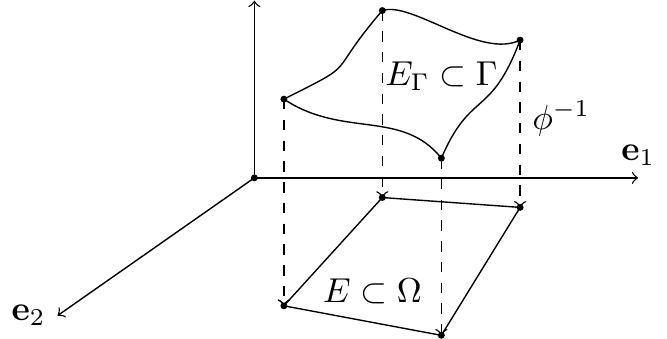} %% {reference_triang_vem}
  }
  \caption{Example of surface polygonal element $\P_{\SurfDomain}$ in
    $\SurfDomain$ and corresponding (planar) element
    $\P$ in $\RefDomain$.}
  \label{fig:triang}
\end{figure}

In this section, we discuss the virtual element
approximation of problem~\ref{eq:varform:surface}.
The numerical method that we use in this work is based on
refs.~\cite{Ahmad-Alsaedi-Brezzi-Marini-Russo:2013,BeiraodaVeiga-Brezzi-Cangiani-Manzini-Marini-Russo:2013,BeiraodaVeiga-Brezzi-Marini-Russo:2016b},
which define optimal approximations of the finite dimensional spaces
on polygonal meshes when the equation coefficients are variable in
space.
As already observed in remark~\ref{rem:onechart}, we work on a single
coordinate neighborhood $\RefDomain=\MapU^{-1}(\SurfDomain)$ and
consider the (global) parametrization
$\MapU:\RefDomain\to\SurfDomain$.
We start from a partition of $\SurfDomain$ formed by surface polygonal
elements $\P_{\SurfDomain}$ with edges denoted by
$\E_{\SurfDomain}$.
Through the parametrization $\MapU$, we can associate the
  partition of $\SurfDomain$ with a partition of $\RefDomain$ formed by
elements $\P$ and possibly curvilinear edges $\E$.
Because of the regularity assumption on the surface, every element
$\P_{\SurfDomain}$ is in a one-to-one relation with one and only one
polygonal element $\P$ in $\RefDomain$.
To avoid curvilinear edges, we use the surface vertices of
$\P_{\SurfDomain}$ to define the vertices of the polygons $\P$ in
$\RefDomain$ through the inverse parametrization and connect
  them with straight segments to define the partition
of $\RefDomain$.
This procedure maintains the above one-to-one relationship between
elements $\P_{\SurfDomain}$ in $\SurfDomain$ and $\P$ in $\RefDomain$
(see fig.~\ref{fig:triang}).

In addition, any function in $\P_{\SurfDomain}$ can be expressed in
$\P$ by composition with the inverse parametrization.
Thus, all the local functional spaces of interest can be defined
indifferently on $\P_{\SurfDomain}$ or $\P$.
The definitions of the building blocks of the virtual element method
is done in $\RefDomain$ using standard two-dimensional Cartesian
coordinates.
These constructions are needed to evaluate the surface bilinear forms
and the right-hand side linear functional of the weak
formulation~\eqref{problem:VEM} by a careful use of the metric tensor.

\PGRAPH{The conforming virtual element space}
%
%% The mesh
Let $\mathcal{T}=\{\Th\}_{\hh}$ be a set of decompositions $\Th$ of
the computational domain $\RefDomain$ into a finite set of
nonoverlapping polygonal elements $\P$.
The subindex label $\hh$ is the maximum of the diameters of the mesh
elements, i.e., $\hP=\sup_{\sv',\sv''\in\P}\abs{\sv'-\sv''}$.
Each element $\P$ has a nonintersecting boundary denoted by
$\partial\P$ formed by straight edges $\E$, center of gravity $\xvP$
and area $\mP$.
A few regularity assumptions are needed on the mesh family $\{\Th\}$
to prove the convergence of the VEM and derive the error estimates in
the $\LTWO$ and $\HONE$ norms.
We present these assumptions at the end of this section where we briefly
discuss the convergence of the proposed VEM.

%% [1]
Let $k\geq1$ be an integer number and $\P\in\Th$ a generic mesh
element.
The \emph{conforming virtual element space $\Vhk$ of order $k\geq1$
built on mesh $\Th$} is obtained by gluing together the local
approximation spaces denoted by $\Vhk(\P)$:
\begin{align}
  \Vhk:=\Big\{\,\vsh\in\HONEzr(\RefDomain)\,:\,\restrict{\vsh}{\P}\in\Vhk(\P)
  \,\,\,\forall\P\in\Th\,\Big\}.
  \label{eq:Vhk:def}
\end{align}
%%
%% [2]
The local virtual element space $\Vhk(\P)$ is defined in accordance
with the \emph{enhancement strategy} introduced
in~\cite{Ahmad-Alsaedi-Brezzi-Marini-Russo:2013}:
\begin{multline}
  \label{eq:VhkP:def}
  \Vhk(\P) = \bigg\{\,
  \vsh\in\HONE(\P)\cap\CS{0}(\overline{\P})\,:\,
  \restrict{\vsh}{\partial\P}\in\CS{0}(\partial\P),\,
  \restrict{\vsh}{\E}\in\PS{k}(\E)\,\forall\E\subset\partial\P,\,
  \Delta\vsh\in\PS{k}(\P)\,\\[0.25em]
  \int_{\P}(\vsh-\PinP{k}\vsh)\,\ms\dV=0
  \,\,\forall\ms\in\PS{k}(\P)\backslash\PS{k-2}(\P)
  \,\bigg\},
\end{multline}
where $\PS{k}(\P)$ and $\PS{k}(\E)$ denote the polynomial spaces of
degree at most $k$ defined over an element $\P$ or an edge $\E$,
respectively.
By definition, each space $\Vhk(\P)$
contains $\PS{k}(\P)$ and the global space $\Vhk$ is a conforming
subspace of $\HONEzr(\RefDomain)$.
The definition of the virtual element bilinear forms
$\ash(\cdot,\cdot)$, $\bsh(\cdot,\cdot)$, and $\csh(\cdot,\cdot)$, and
the forcing term $\Fsh(\cdot)$ requires the definition of the elliptic
and orthogonal projections operators.

\PGRAPH{Elliptic projection}
%
%% [3]
%%
The \emph{elliptic projection operator}
$\PinP{k}:\HONE(\P)\to\PS{k}(\P)$ can be defined for any
$\vsh\in\Vhk(\P)$ as:
\begin{align}  
  \int_{\P}\nabla\PinP{k}\vsh\cdot\nabla\qs\dV  &= \int_{\P}\nabla\vsh\cdot\nabla\qs\dV\quad\forall\qs\in\PS{k}(\P),\label{eq:PiFn}\\[0.5em]
  \int_{\partial\P}\big(\PinP{k}\vsh-\vsh\big)\dS &= 0.                                                              \label{eq:def:Pib_k}
\end{align}
Equation~\eqref{eq:def:Pib_k} allows the removal of the kernel of the
gradient operator.
The elliptic projection operator $\PinP{k}$ is a polynomial-preserving operator, i.e.,
$\PinP{k}\qs=\qs$ for every $\qs\in\PS{k}(\P)$.
One of its major property is that the
projection $\PinP{k}\vsh$ of any virtual element function
$\vsh\in\Vhk(\P)$ is computable from the degrees of freedom of
$\vsh$~\cite{BeiraodaVeiga-Brezzi-Cangiani-Manzini-Marini-Russo:2013},
which are defined as follows.

%% [4] DOFS
%%
The degrees of freedom of the virtual element function
$\vsh\in\Vhk(\P)$ are given by the set of values:
\begin{description}
\item[]\textbf{(D1)} for $k\geq1$, the values of $\vsh$ at the
  vertices of $\P$;

%%   \medskip
%% \item[]\textbf{(D2)} for $k\geq2$, the edge moments of $\vsh$ of order
%%   up to $k-2$ on every $\E\in\partial\P$:
%%   \begin{align}
%%     \frac{1}{\mE}\int_{\E}\vsh\,\ms\dS,
%%     \,\,\forall\ms\in\calM{k-2}(\E),\,
%%     \forall\E\in\partial\P;
%%     \label{eq:dofs:D2}
%%   \end{align}
  
  \medskip
\item[]\textbf{(D2)} for $k\geq2$, the values of $\vsh$ at the
  $k-1$ internal nodes of the $k$-th Gauss-Lobatto formula on
  every $\E\in\partial\P$;
  
  \medskip
\item[]\textbf{(D3)} for $k\geq2$, the cell moments of $\vsh$ of order
  up to $k-2$ on element $\P$:
  \begin{align}
    \frac{1}{\mP}\int_{\P}\vsh\,\ms\dV,
    \,\,\forall\ms\in\calM{k-2}(\P),
    \label{eq:dofs:D3}
  \end{align}
\end{description}
where $\calM{k-2}(\P)$ is the set of scaled monomials that span the
linear space of polynomials of degree up to $k-2$.
These set of values are unisolvent in $\Vhk(\P)$,
cf.~\cite{BeiraodaVeiga-Brezzi-Cangiani-Manzini-Marini-Russo:2013},
and thus, every virtual element function is uniquely identified by
them.
The degrees of freedom of a virtual element function in the global
space $\Vhk$ are given by collecting the elemental degrees of freedom
\textbf{(D1)}-\textbf{(D3)}.
Their unisolvence in $\Vhk$ is an immediate consequence of their
unisolvence in every elemental space $\Vhk(\P)$.

\medskip
\PGRAPH{Orthogonal projections}
%%
%% [5] orthogonal projections of the vem functions
From the degrees of freedom of a virtual element function
$\vsh\in\Vhk(\P)$ we can also compute the orthogonal projections
$\PizP{k}\vsh$ and $\PizP{k-1}\nabla\vsh$,
cf.~\cite{Ahmad-Alsaedi-Brezzi-Marini-Russo:2013}.
In fact, the definition of the orthogonal projection $\PizP{k}\vsh$
reads as
\begin{align}
  \int_{\P}\PizP{k}\vsh\,\qs\dV =
  \int_{\P}\vsh\,\qs\dV\qquad\forall\qs\in\PS{k}(\P).
\end{align}
The right-hand side is the integral of $\vsh$ against the polynomial
$\qs$, and is computable from the degrees of freedom \textbf{(D3)} of
$\vsh$ when $\qs$ is a polynomial of degree up to $k-2$, and from the
moments of $\PinP{k}\vsh$ when $\qs$ is a polynomial of degree $k-1$
and $k$, cf.~\eqref{eq:VhkP:def}.
Clearly, the orthogonal projection $\PizP{k-1}\vsh$ is also
computable.

%% [6] orthogonal projections of the gradients
In turn, using the definition of the orthogonal projection
$\PizP{k-1}\nabla\vsh$ and integrating by parts, we find that
\begin{align}
  \int_{\P}\PizP{k-1}\nabla\vsh\cdot\qv\dV
  = \int_{\P}\nabla\vsh\cdot\qv\dV
  = -\int_{\P}\vsh\nabla\cdot\qv\dV + \sum_{\E\in\partial\P}\int_{\E}\vsh\norPE\cdot\qv\dS
\end{align}
for every $\qv\in\left[\PS{k-1}(\P)\right]^{2}$, where $\norPE$
denotes the unit outward vector orthogonal to the boundary edge
$\E\in\partial\P$.
The first integral on the (last) right-hand side is computable from
the degrees of freedom \textbf{(D3)}, i.e., from the moments of $\vsh$
against the polynomials of degree $k-2$ over $\P$.
The edge integrals are computable from the degrees of freedom
\textbf{(D1)}-\textbf{(D2)} because we can compute the trace of $\vsh$
on each edge by interpolating these nodal values.
%% , i.e., from the moments of $\vsh$ against
%% the polynomials of degree $(k-1)$ over the edge.
%% %%

\PGRAPH{The virtual element bilinear forms}
Following the VEM gospel, we write the discrete bilinear forms
$\ash(\cdot,\cdot)$, $\bsh(\cdot,\cdot)$ and $\csh(\cdot,\cdot)$ as the
sum of elemental contributions
\begin{align}
  \ash(\ush,\vsh) = \sum_{\P\in\Th}\ashP(\ush,\vsh),\quad
  \bsh(\ush,\vsh) = \sum_{\P\in\Th}\bshP(\ush,\vsh),\quad
  \csh(\ush,\vsh) = \sum_{\P\in\Th}\cshP(\ush,\vsh).
\end{align}
The bilinear forms $\ashP(\ush,\vsh)$, $\bshP(\ush,\vsh)$ and
$\cshP(\ush,\vsh)$ on each element $\P$ are given by
\begin{align}
  \label{eq:discreteforms}
  \ashP(\ush,\vsh) 
  &= \int_{\P}\sqrt{\det{\matG}}\matG^{-1}\,\PizP{k-1}\nabla\ush\cdot\PizP{k-1}\nabla\vsh\dV
  + \SPh\Big( \big(I-\PinP{k}\big)\ush, \big(I-\PinP{k}\big)\vsh \Big),\\[0.5em]
  %% ------------------
  \bshP(\ush,\vsh)
  &= \int_{\P}\sqrt{\det{\matG}}\matG^{-1/2}\hat{\AdvVel}\cdot\PizP{k-1}\Grad\ush\;\PizP{k-1}\vsh\dV,\\[0.5em]
    %% ------------------
  \cshP(\ush,\vsh) 
  &= \int_{\P}\sqrt{\det{\matG}}\MassCoef\PizP{k-1}\ush\,\PizP{k-1}\vsh\dV.
\end{align}
The bilinear form $\SPh(\cdot,\cdot)$ in the definition of
$\ashP(\cdot,\cdot)$ provides the stability term and can be any
symmetric positive definite bilinear form defined on $\P$ for which
there exist two positive constants $\cbot$ and $\ctop$ such that
\begin{align}
  \cbot\asP(\vsh,\vsh)
  \leq\SPh(\vsh,\vsh)
  \leq\ctop\asP(\vsh,\vsh)
  \quad\forall\vsh\in\Vhks(\P)\textrm{~with~}\PinP{k}\vsh=0.
  \label{eq:SP:stability}
\end{align}
Note that $\SPh(\cdot,\cdot)$ must scale like the restriction of
$\as(\cdot,\cdot)$ on the mesh element $\P$.
Also, the stabilization term in the definition of $\ashP(\cdot,\cdot)$
gives a zero contribution if one of its two entries is a polynomial of
degree (at most) $k$ since $\PinP{k}$ is a projection on the
polynomial space.
In this work, we consider two possible implementations of the
stability term:
\begin{itemize}
\item the choice originally provided
  in~\cite{BeiraodaVeiga-Brezzi-Cangiani-Manzini-Marini-Russo:2013},
  which is sometimes called the ``\emph{dofi-dofi stabilization}'' in
  the virtual element literature, and reads as
  \begin{align}
    \SPh(\vsh,\wsh) =
    \sum_{i=1}^{\NDOFS}\textrm{DOF}_i(\vsh)\textrm{DOF}_i(\wsh),
  \end{align}
  where $\textrm{DOF}_i(\cdot)$ is the map between a virtual function
  and its degrees of freedom;

\item the formula proposed in~\cite{Mascotto:2018}, which is sometimes
  called the ``\emph{D-recipe stabilization}'' in
  the virtual element literature, and reads as
  \begin{align}
    \SPh(\vsh,\wsh) =
    \sum_{i=1}^{\NDOFS}\mathcal{A}_{ii}\textrm{DOF}_i(\vsh)\textrm{DOF}_i(\wsh),
  \end{align}
  where $\mathcal{A}$ is the matrix resulting from the implementation
  of the first term in the bilinear form $\ashP(\cdot,\cdot)$:
  \begin{align}
    \mathcal{A}_{ij}:=\int_{\P}\sqrt{\det{\matG}}\matG^{-1}\,\PizP{k-1}\nabla\varphi_i\cdot\PizP{k-1}\nabla\varphi_j\dV,
  \end{align}
  where $\varphi_i$ (and $\varphi_j$) are the ``canonical'' basis
    functions generating $\Vhk(\P)$, i.e., the functions whose $i-th$
    (or $j-th$) degree of freedom is equal to $1$ and all other
    degrees of freedom are $0$.
    We note that these basis function are unknown in the virtual
    element framework, but their projections $\PizP{k-1}\nabla\varphi_i$
    (and $\PizP{k-1}\nabla\varphi_j$) are computable.
\end{itemize}

\medskip
The stabilization term, and, in particular,
condition~\eqref{eq:SP:stability}, is designed in order that
$\ashP(\cdot,\cdot)$ satisfies the two fundamental properties:
\begin{description}
\item[-] {\emph{$k$-consistency}}: for all $\vsh\in\Vhks$ and for all
  $\qs\in\PS{k}(\P)$ it holds
  \begin{align}
    \label{eq:k-consistency}
    \ashP(\vsh,\qs) = \asP(\vsh,\qs);
  \end{align}
  
  \medskip
\item[-] {\emph{stability}}: there exist two positive constants
  $\alpha_*,\,\alpha^*$, independent of $\hh$ and $\P$, such that
  \begin{align}
    \label{eq:stability}
    \alpha_*\asP(\vsh,\vsh)
    \leq\ashP(\vsh,\vsh)
    \leq\alpha^*\asP(\vsh,\vsh)\quad\forall\vsh\in\Vhks.
  \end{align}
\end{description} 

%% The right-hand side
\PGRAPH{The virtual element forcing term}
To approximate the right-hand side of~\eqref{eq:VEM}, we split the
term into the sum of elemental contributions and approximate every
local linear functional by means of the orthogonal projection
$\PizP{k}\vsh$:
\begin{align}
  &\RhsApprox{\vsh} = \sum_{\P\in\Th}\big(\fs,\PizP{k}\vsh\big)_{\P}
  \textrm{~~where~~}\big(\fs,\PizP{k}\vsh\big)_{\P} = \int_{\P}\sqrt{\det{\matG}}\,\fs\,\PizP{k}\vsh\dV.
  \label{eq:extensionc}
\end{align}
With these definitions the VEM scheme in problem \ref{problem:VEM} is
completely determined.

\PGRAPH{Convergence properties}
The numerical analysis of the scheme requires the following hypotheses
on the mesh, typical of VEM methods. 
\begin{assumption}[Mesh regularity assumptions]
  \label{assum:mesh:regularity}
  There exists a positive constant $\varrho$ independent of $\hh$
  (and, hence, of $\Th$) such that
  \begin{description}
  \item[$(i)$] every element $\P$ of every mesh $\Th$ is star-shaped
    with respect to a disk with radius $\ge\varrho\hP$;
    
    \smallskip
  \item[$(ii)$] every edge $\E\in\partial\P$ has length
    $\hE\geq\varrho\hP$.
  \end{description}
\end{assumption}
The star-shapedness property $(i)$ implies that the polygonal elements
are \emph{simply connected} subsets of $\REAL^{2}$.
In turn, the scaling assumption $(ii)$ implies that the number of
edges in each elemental boundary is uniformly bounded over the whole
mesh family $\{\Th\}$.

\medskip
The following theorem summarizes the results for the virtual element
approximation in problem~\ref{eq:varform:surface}.
The proof of these results found
in~\cite{BeiraodaVeiga-Brezzi-Marini-Russo:2016b} is easily extended
to our setting.
Indeed, we choose to write the theorem in terms of the chart
$\RefDomain$ and its discretization $\RefDomain_{\hh}$, but it can be
written equivalently in terms of the surface $\SurfDomain$ and its
discretization $\SurfDomain_{\hh}$, since the norms of the
parametrization and its inverse are uniformly bounded by hypothesis.

\smallskip
\begin{theorem}\label{theorem:GDM:apriori:estimate}
  Let $\us\in\HS{k+1}(\RefDomain)\cap\HONEzr(\RefDomain)$, be the solution to the variational
  problem~\ref{eq:varform:surface} on a convex domain $\RefDomain$
  with ${\fs\in\HS{k}(\RefDomain)}$.
  Let $\ush\in\Vhks$ be the solution of the virtual element method
  \eqref{eq:VEM} on every mesh of a mesh family $\{\Th\}$ satisfying
  the mesh regularity assumption~\ref{assum:mesh:regularity}.
  Then, a strictly positive constant $\Cs$ independent of $\hh$ exists
  such that
  \begin{itemize}
  \item the $\HONE$-error estimate holds:
    \begin{align}
      \label{eq:source:problem:H1:error:bound}
      \norm{ \us-\ush }{\HONE(\RefDomain)}\leq
      \Cs\hh^{k}\left(
      \snorm{\us}{\HS{k+1}(\RefDomain)} 
      + \snorm{\fs}{\HS{k}(\Omega)}
      \right);
    \end{align}

    \medskip
  \item the $\LTWO$-error estimate holds:
    \begin{align}
      \label{eq:source:problem:L2:error:bound}
      \norm{ \us-\ush}{\LTWO(\RefDomain)}\leq 
      \Cs\hh^{k+1}\left(
      \snorm{\us}{\HS{k+1}(\RefDomain)}
      + \snorm{\fs}{\HS{k}(\Omega)}
      \right).
    \end{align}
  \end{itemize}
  The constant $\Cs$ may depend on the coefficient bounds $\kappa_{*}$,
  $\kappa^*$, $\ws_{\max}$ and $\gamma_{\max}$, the stability
  constants $\alpha_*$ and $\alpha^*$, the mesh regularity constant
  $\varrho$, the size of the computational domain $\ABS{\Omega}$, and
  the approximation degree $k$.
\end{theorem}

\smallskip
The approximate solution $\ush$ is not explicitly known inside
the elements.
Consequently, in the numerical experiments of
Section~\ref{sec:numerical-results}, we approximate the error norms as
follows:
\begin{align*}
  \norm{ \us-\ush }{\HONE(\RefDomain)}\approx\norm{ \us-\Piz{k}\ush }{\HONE(\Th)}
  \quad\textrm{and}\quad
  \norm{ \us-\ush }{\LTWO(\RefDomain)}\approx\norm{ \us-\Piz{k}\ush }{\LTWO(\RefDomain)}.
\end{align*}
Here, $\Piz{k}\ush$ is the global projector on the space of
discontinuous polynomials of degree at most $k$ built on mesh $\Th$,
and $\norm{ \us-\Piz{k}\ush }{\HONE(\Th)}$ is the norm in the broken
Sobolev space $\HONE(\Th)$ that is defined by summing the
$\HONE(\P)$-norms of each element $\P$.
Operator $\Piz{k}\ush$ is obtained by taking the elemental
$\LTWO$-orthogonal projections $\PizP{k}\ush$ in every mesh element
$\P$, which are computable from the degrees of freedom of $\ush$, so
that
$\restrict{\big(\Piz{k}\ush\big)}{\P}=\PizP{k}\big(\restrict{\ush}{\P}\big)$.

%% SECTION 4
\section{Numerical Results}
\label{sec:numerical-results}

\begin{figure}
  \centerline{
    \includegraphics[width=0.27\textwidth]{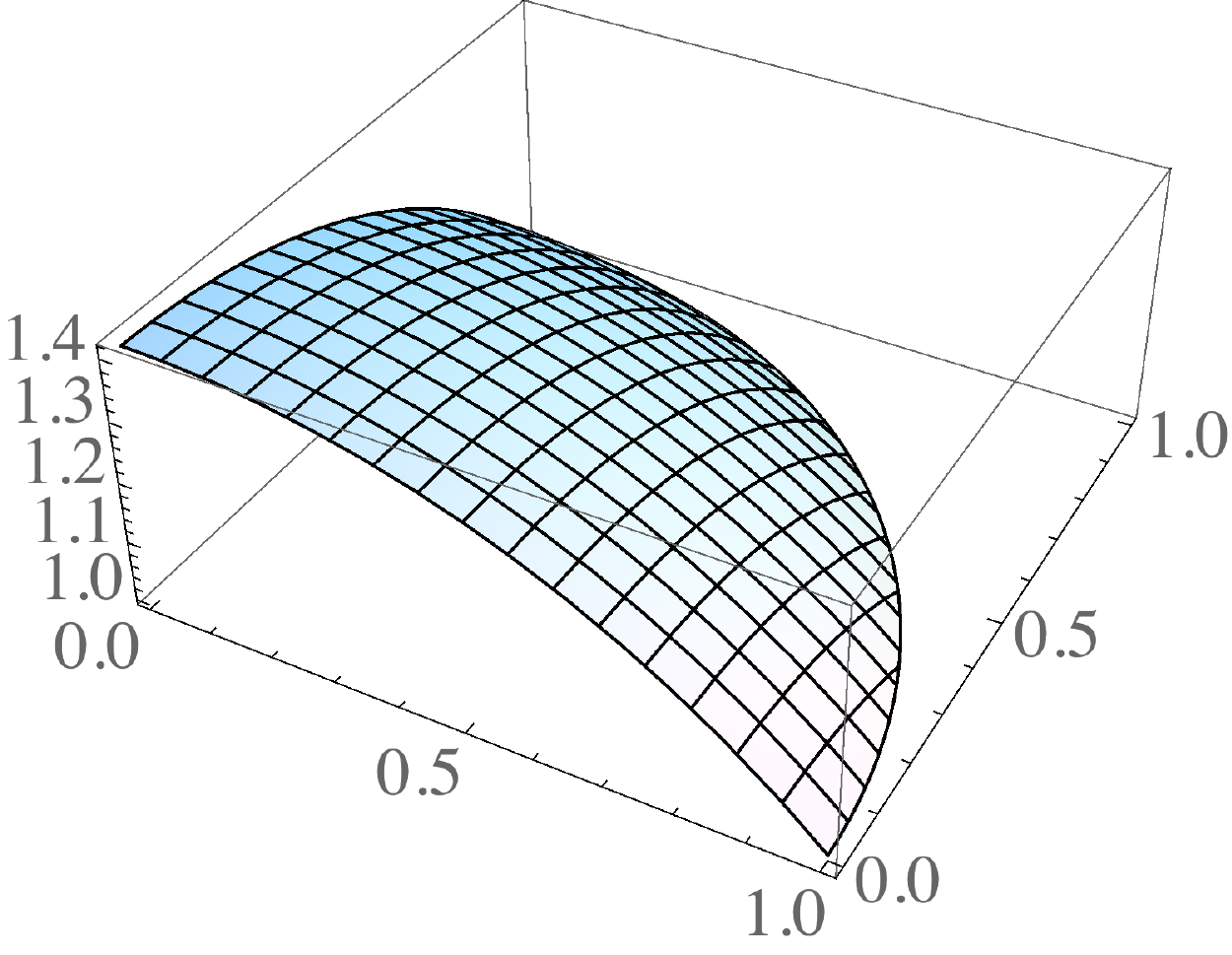}%%{surface_r=2-a=0}
    \includegraphics[width=0.22\textwidth]{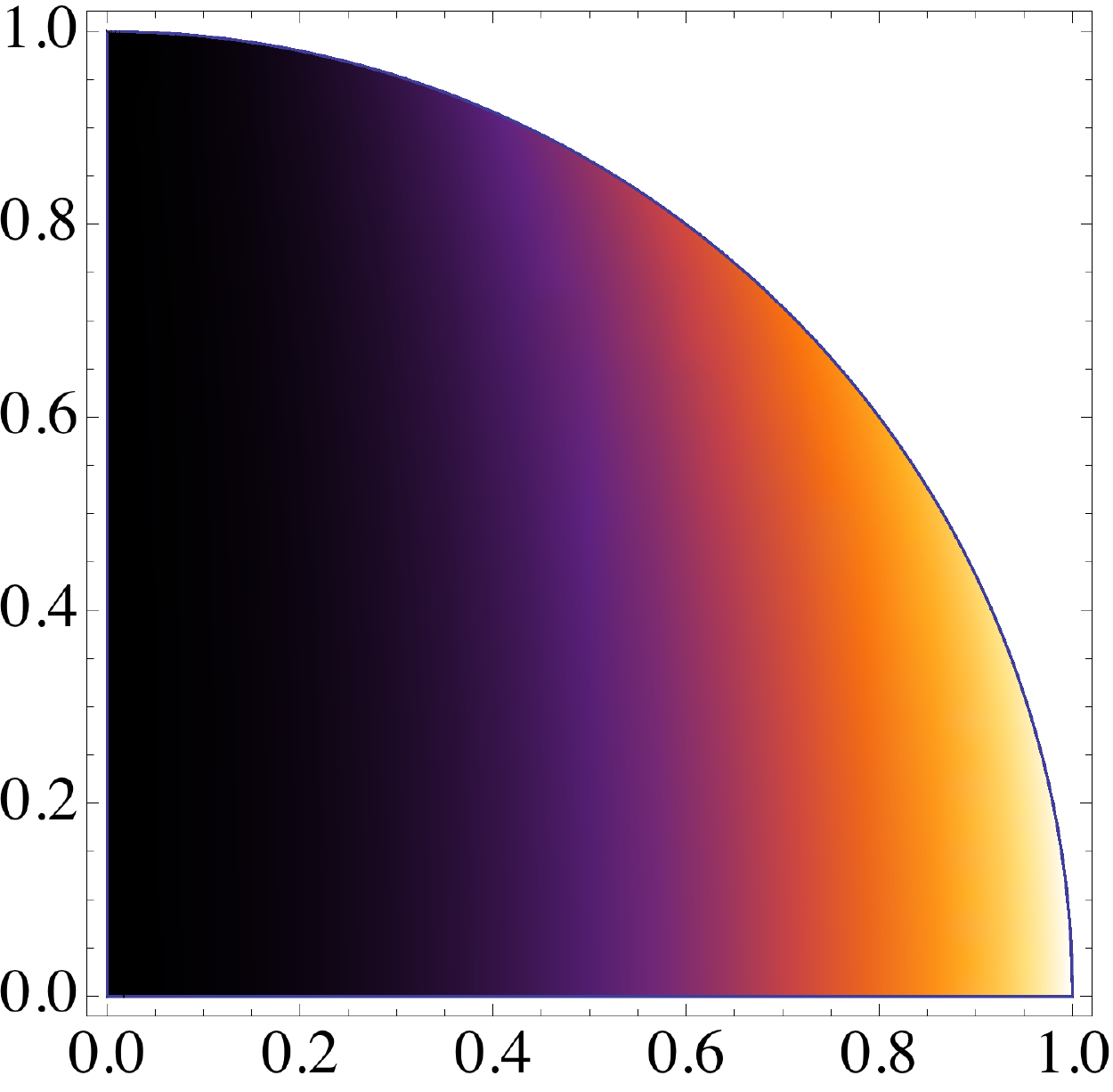}%%{g11_r=2-a=0}
    \includegraphics[width=0.22\textwidth]{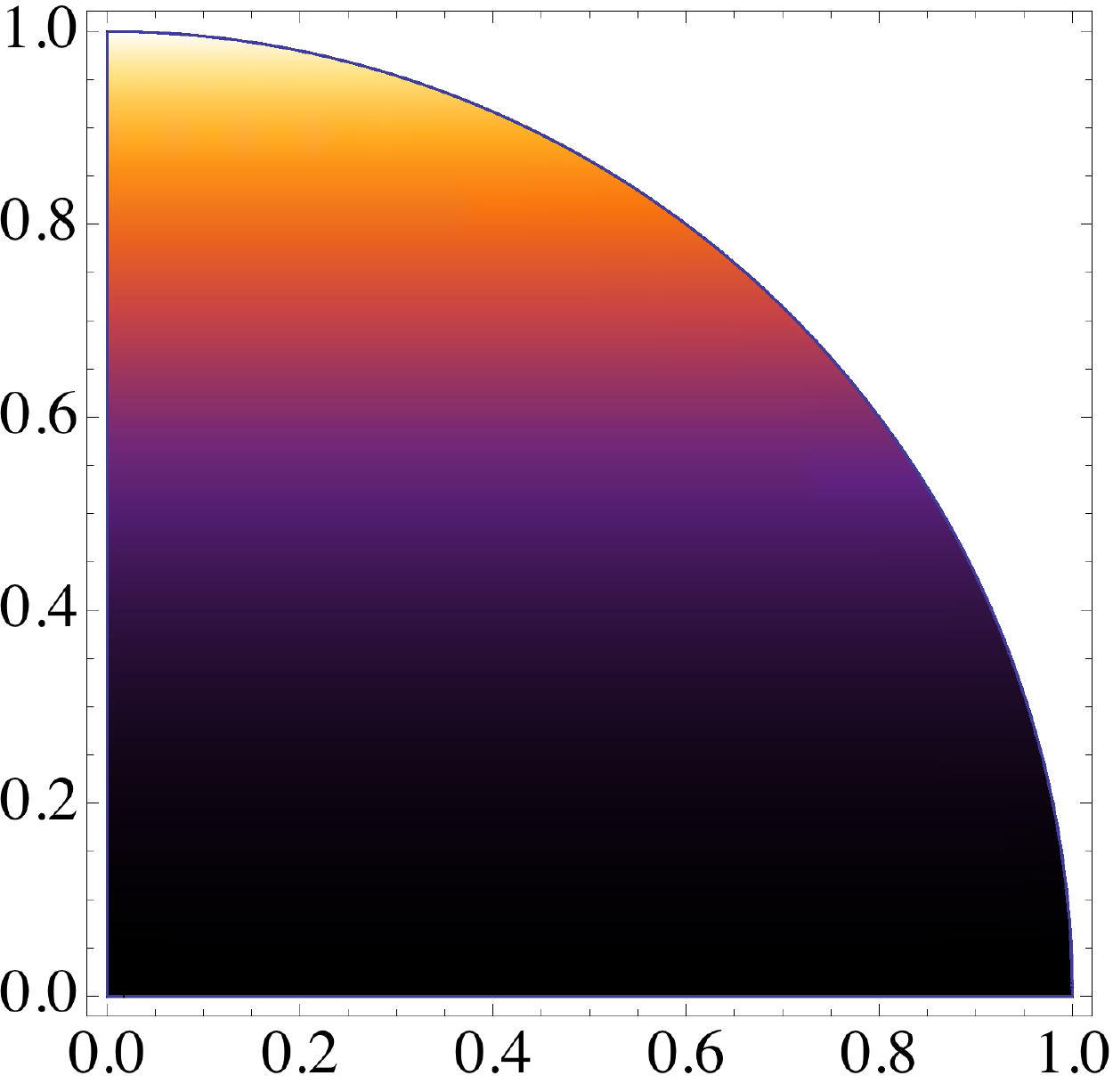}%%{g22_r=2-a=0}
    \includegraphics[width=0.22\textwidth]{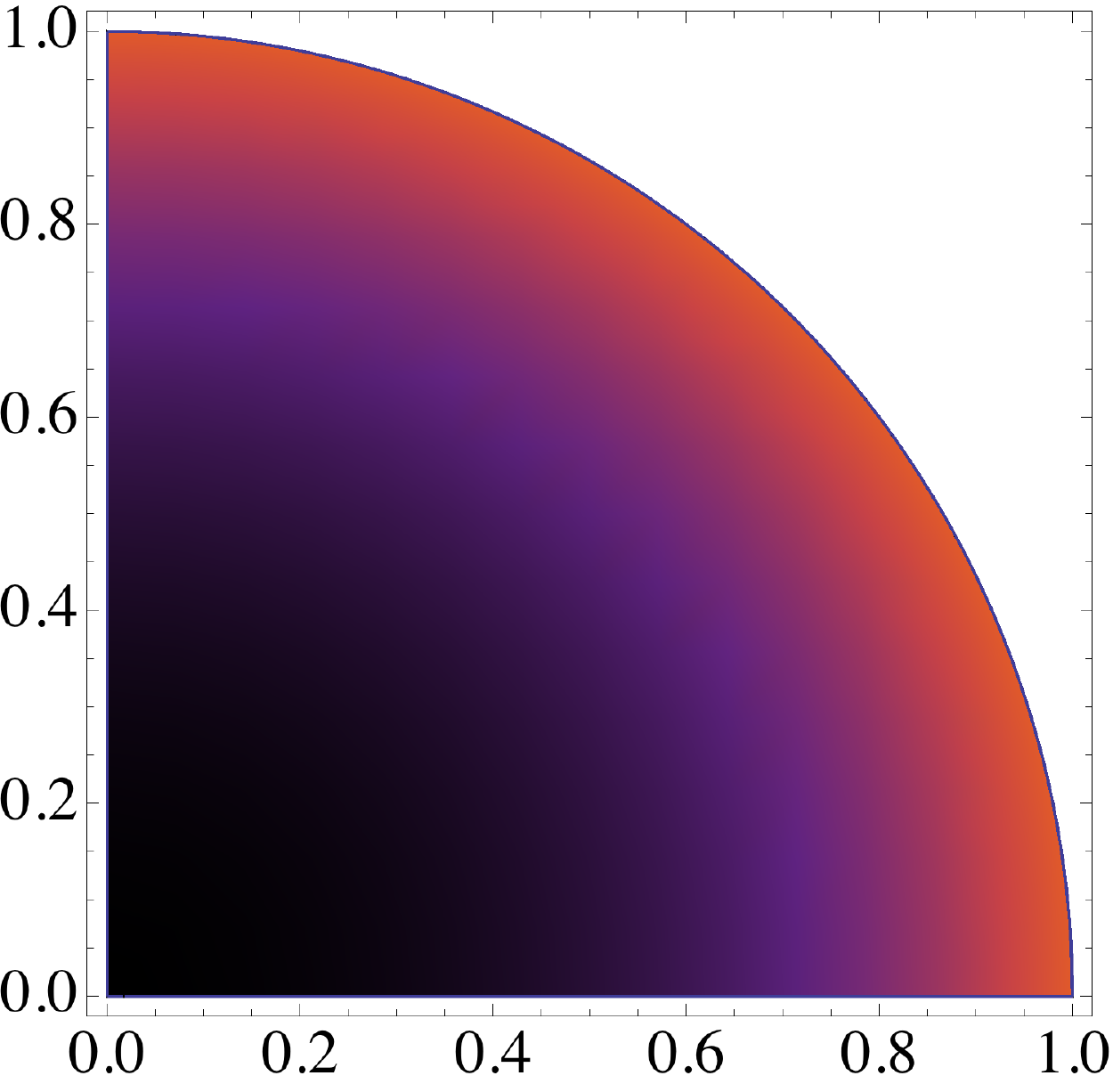}%%{sqrt-det-G_r=2-a=0}
    \includegraphics[width=0.03\textwidth]{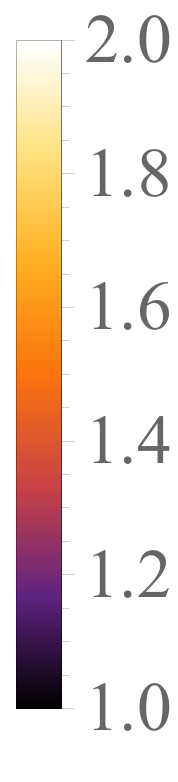}%%{legend-r=2-a=0}
  }
  \centerline{
    \includegraphics[width=0.27\textwidth]{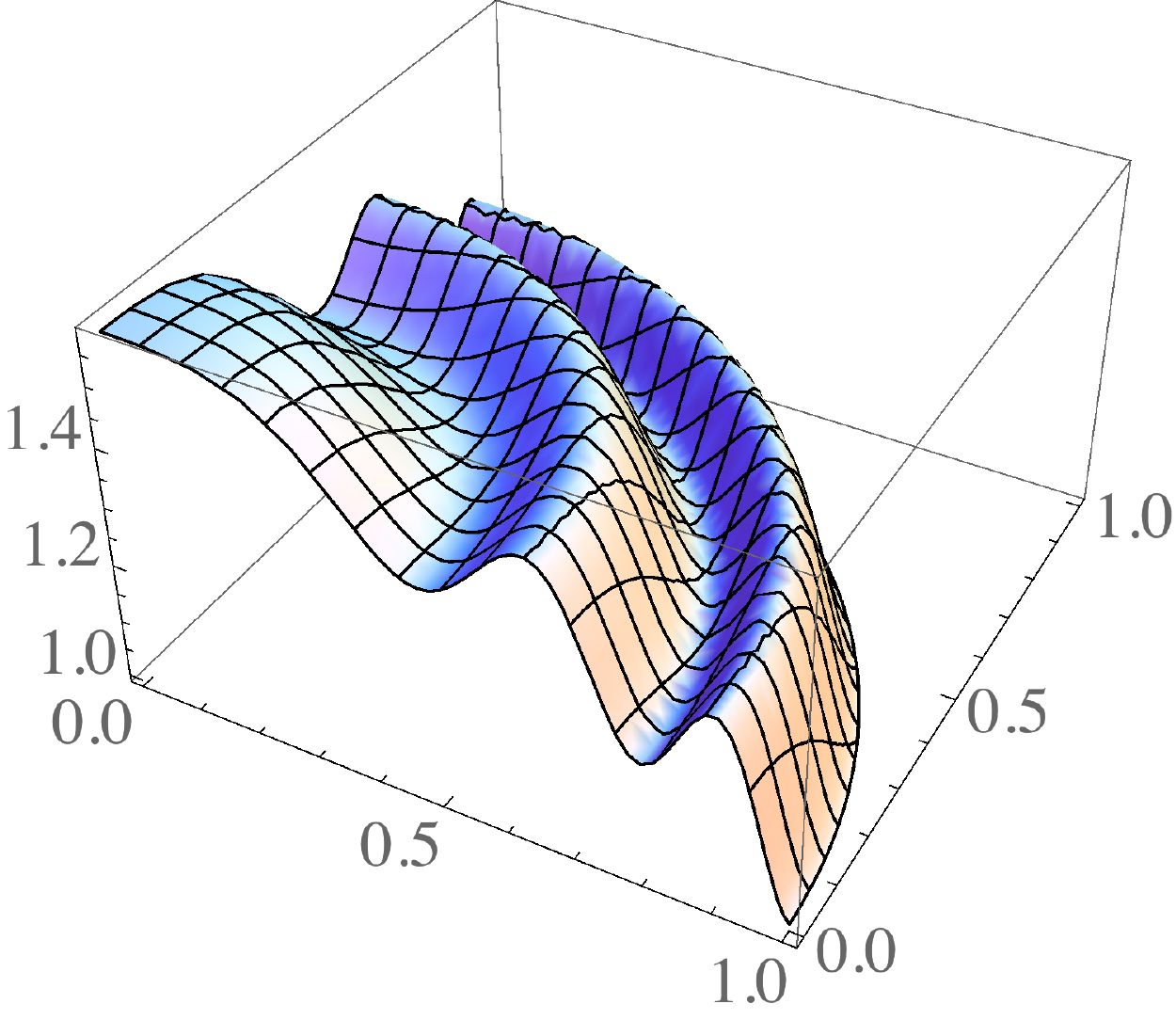}%%{surface_r=2-a=05-k=5}
    \includegraphics[width=0.22\textwidth]{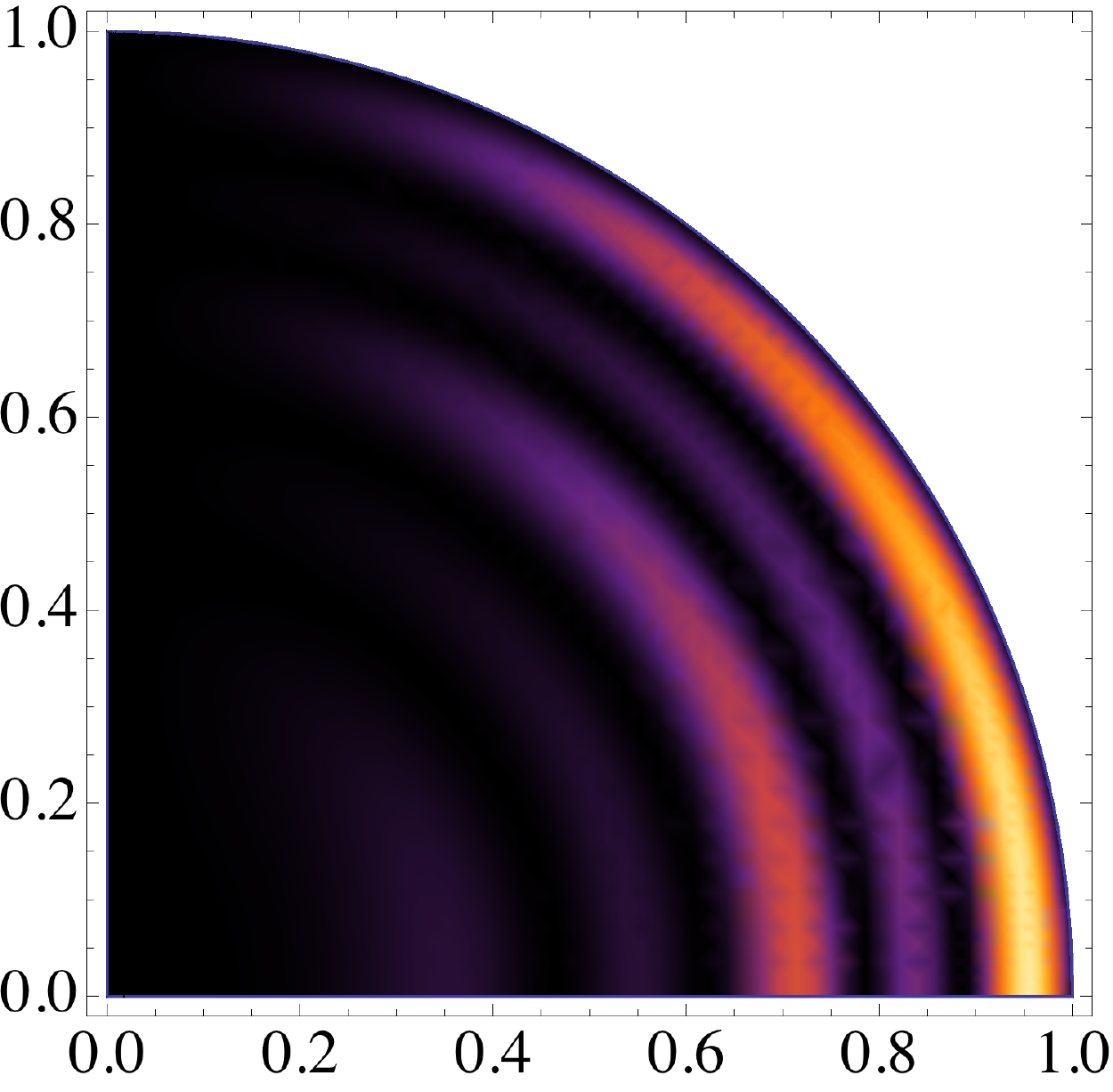}%%{g11_r=2-a=05-k=5}
    \includegraphics[width=0.22\textwidth]{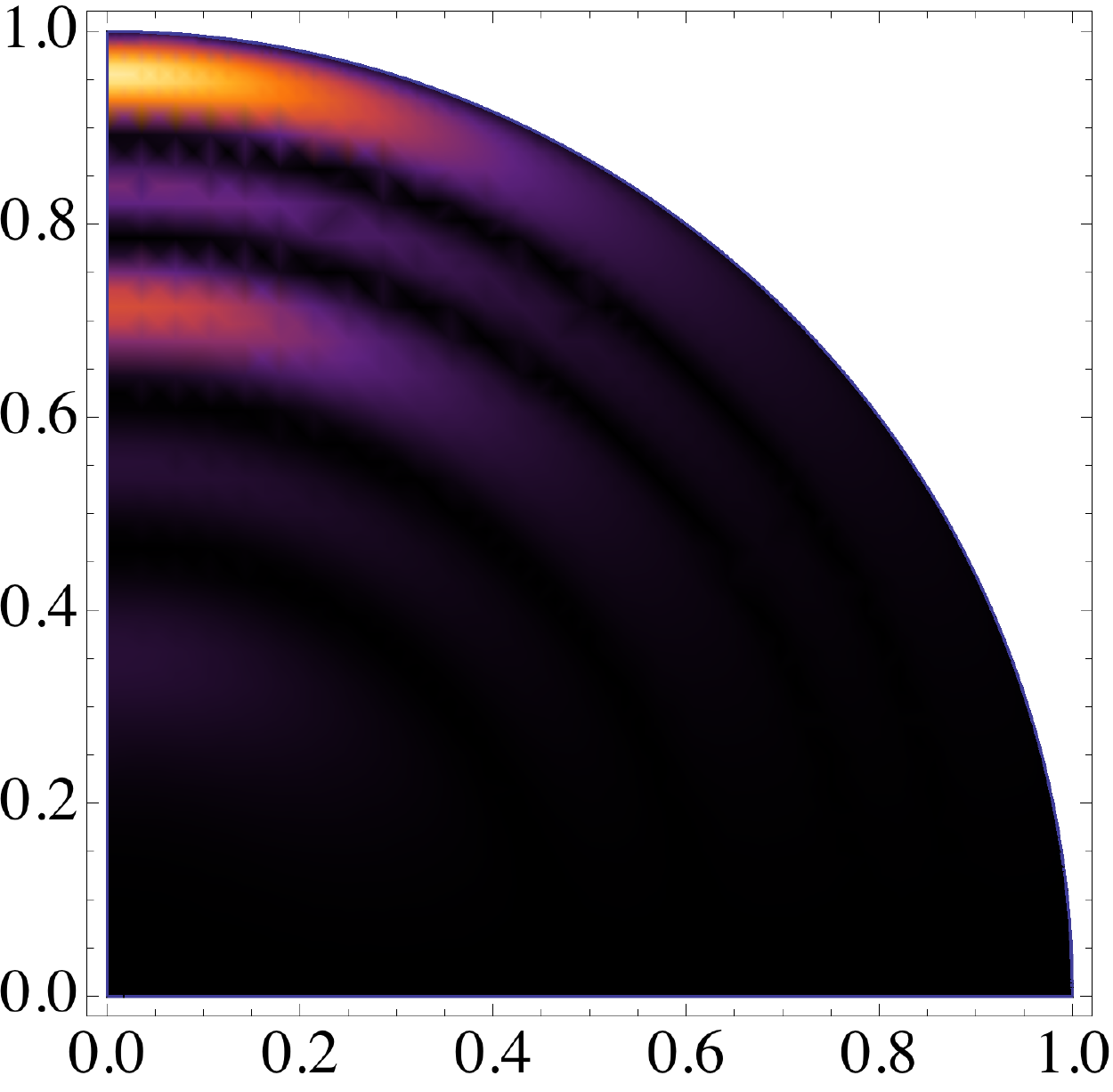}%%{g22_r=2-a=05-k=5}
    \includegraphics[width=0.22\textwidth]{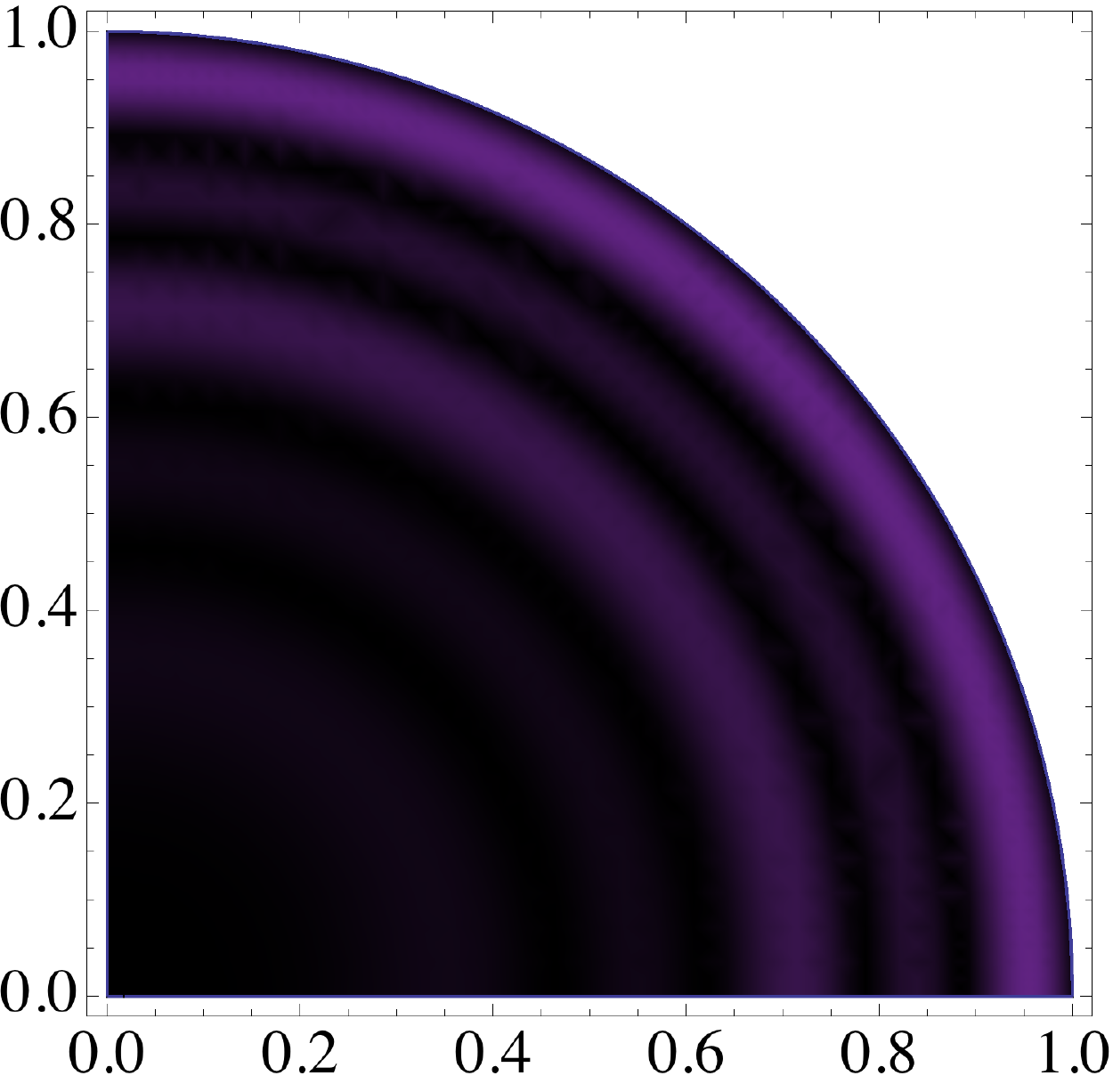}%%{sqrt-det-G_r=2-a=05-k=5}
    \includegraphics[width=0.03\textwidth]{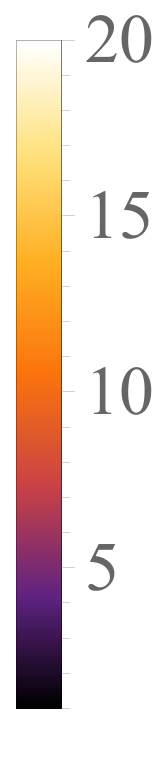}%%{legend-r=2-a=05-k=5}
  }
  \centerline{
    \includegraphics[width=0.27\textwidth]{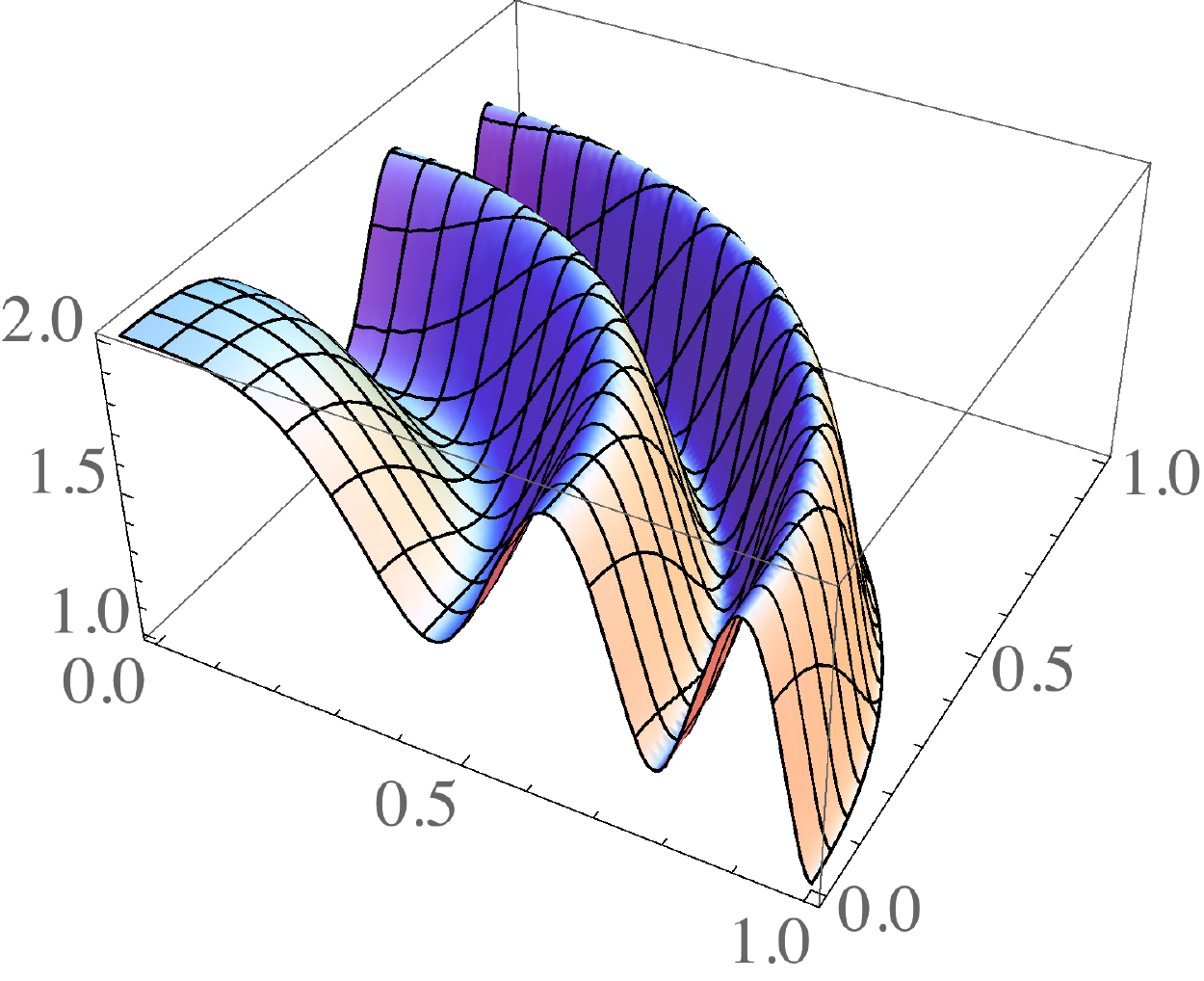}%%{surface_r=2-a=20-k=5}
    \includegraphics[width=0.22\textwidth]{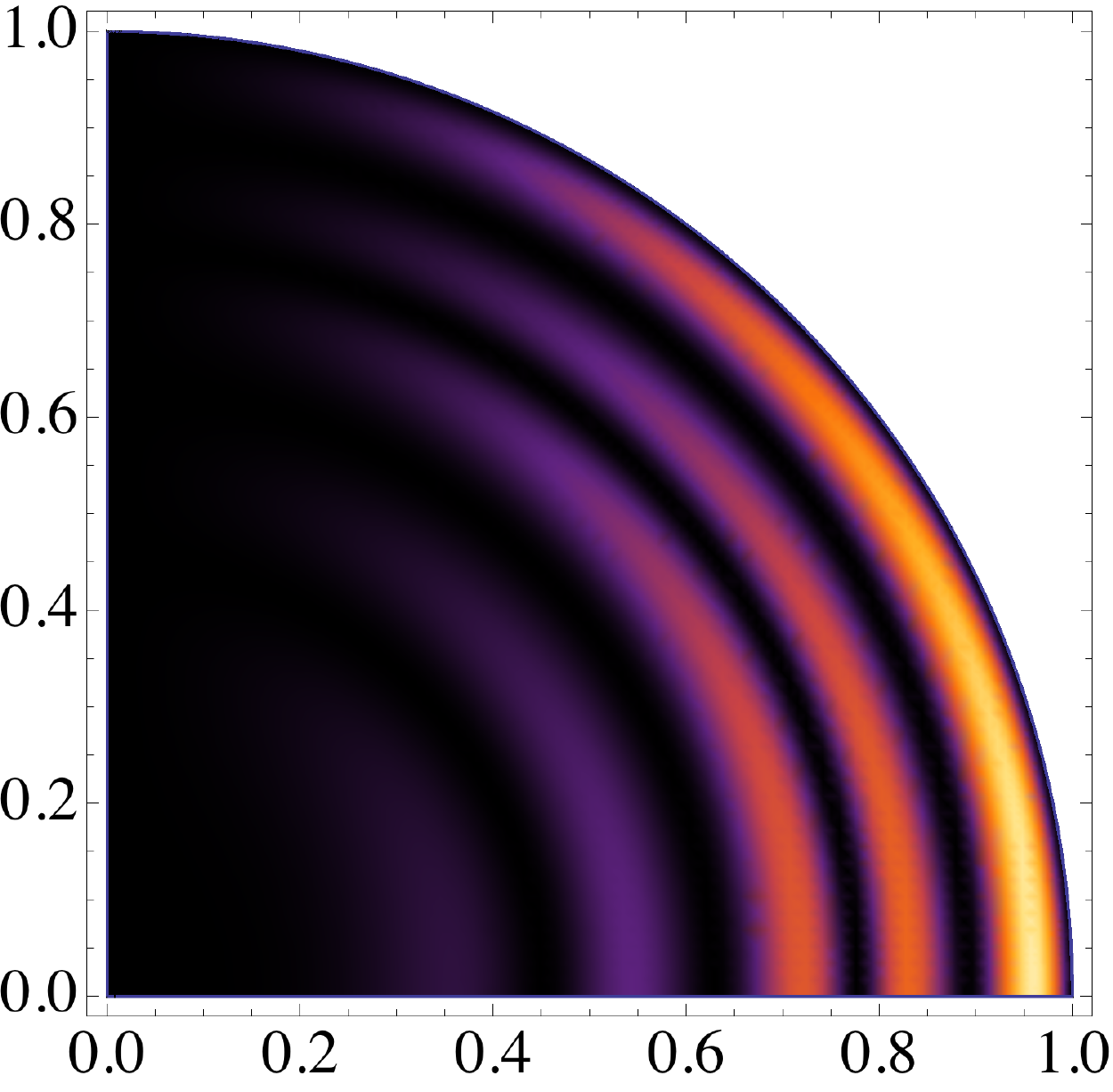}%%{g11_r=2-a=20-k=5}
    \includegraphics[width=0.22\textwidth]{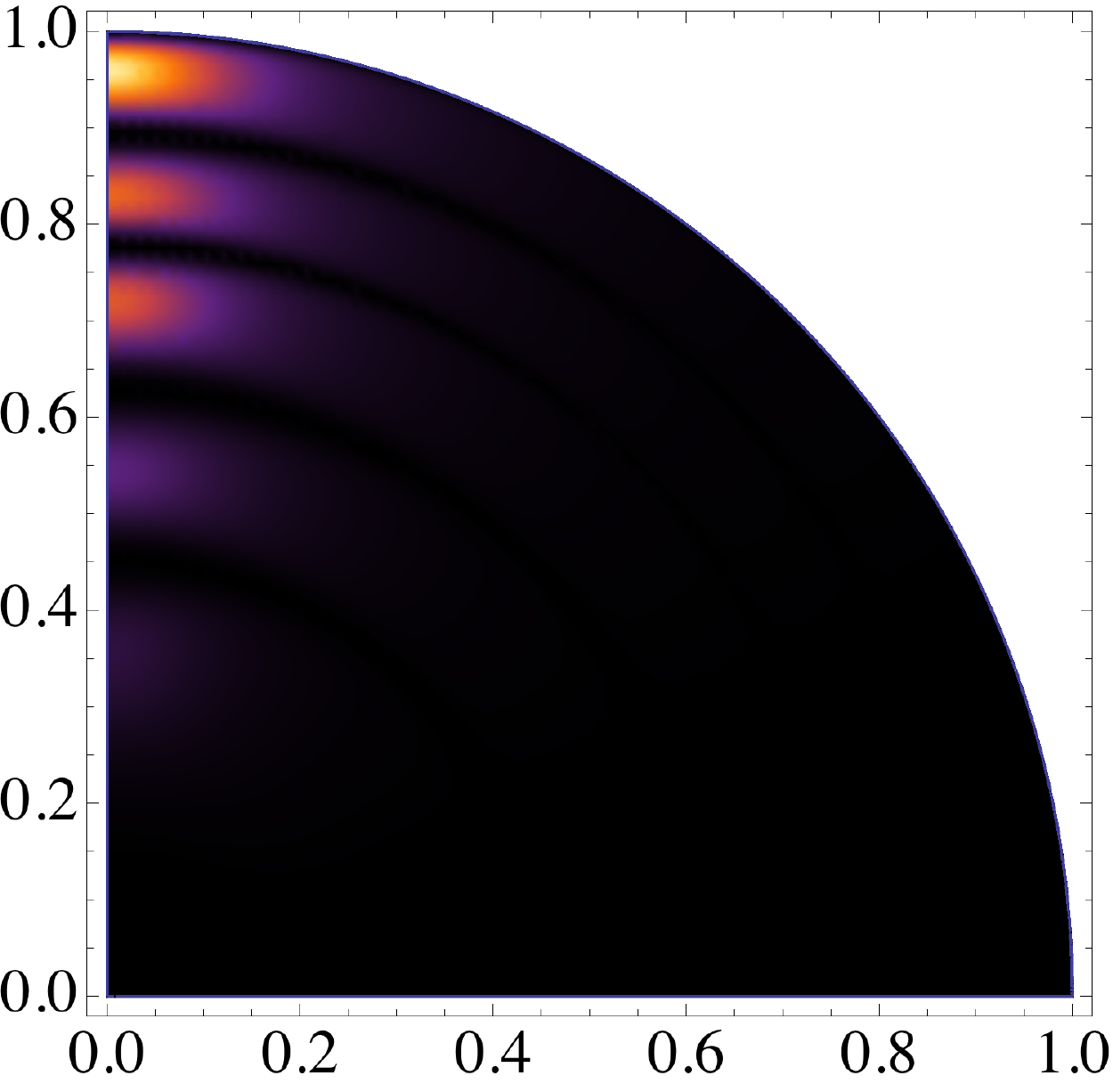}%%{g22_r=2-a=20-k=5}
    \includegraphics[width=0.22\textwidth]{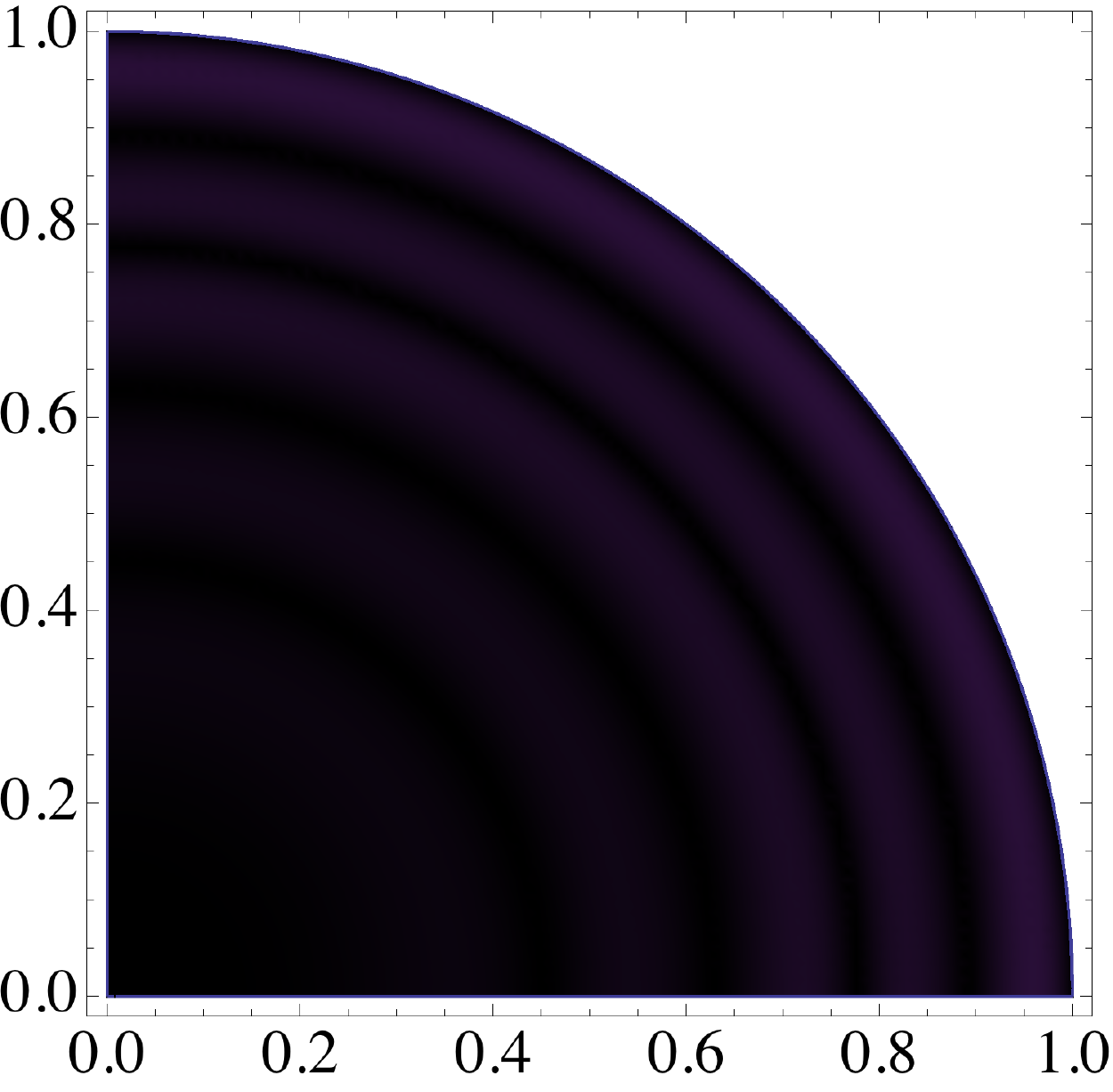}%%{sqrt-det-G_r=2-a=20-k=5}
    \includegraphics[width=0.03\textwidth,height=0.11\textheight]{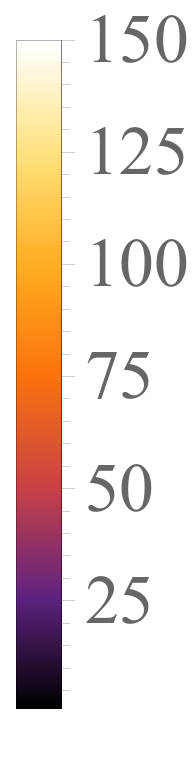}%%{legend-r=2-a=20-k=5}
  }
  \caption{Surfaces and metric components used in the numerical
    experiments (see eq.~(\ref{eq:surface-trig})).
    The columns show the spatial behavior of the surface $\SurfDomain$,
    $\metrTensCv{11}$, $\metrTensCv{22}$, and $\sqrt{\DET{\First}}$
    ($r=2$, $a=0$ first row; $r=2$, $a=0.5$, $k=5$ second row; $r=2$,
    $a=2.0$, $k=5$ third row).
    Note the completely different color scales between the case $a=0$
    and the cases $a>0$, where the metric tensor displays strong
    anisotropy ($\metrTensCv{11}\approx 20 \,(150)$ where
    $\metrTensCv{22}\approx 1$).
    The choice of an orthogonal reference frame ensures that the
    principal directions of anisotropy (the eigenvectors of $\First$) in the second order diffusion
    term do not vary in space.}
  \label{fig:surfaces}
\end{figure}

In this section we present numerical results on synthetic test cases
to support the statements of the previous sections by means of
experimental evidence.
Our test cases are grouped into four main categories.
The first two sets of experiments, Test Cases~1 and~2, are aimed at
showing the correctness of our implementation and the order of
convergence of the proposed VEM scheme up to fourth order of accuracy.
In the third set of experiments, Test Case~3, we explore the limits of
the VEM approach as the metric tensor $\First$ becomes more and more
anisotropic (large condition numbers) as a function of the regularity
of the surface.
Finally, Test Case~4 considers the use of stereographic projection to
build a two chart atlas for the sphere to show the
applicability of our approach in a multi-chart case.  

In the first three experiments we consider the surface provided by the
graph of the following height function, a simple trigonometric
perturbation of a portion of a sphere embedded in $\REAL^3$:
\begin{equation}\label{eq:surface-trig}
  \zcg=\height(\xcg,\ycg)=\sqrt{r - (\xcg)^2 - (\ycg)^2 +
    a \cos^2\left(k\frac{\pi}{2} ((\xcg)^2 + (\ycg)^2)\right)} \; ,
\end{equation}
where $r$ is the radius of the sphere, and $a$ and $k$ are the
amplitude and the frequency of the cosine trigonometric perturbation.
We use the Monge parametrization given by
$\MapU=\{\xcg=\xcl,\ycg=\ycl,\zcg=\height(\xcg,\ycg)\}$ and work on
the single chart represented by the domain $\RefDomain=\{(\xcl,\ycl):
\xcl,\ycl\ge 0 \mbox{ and } \sqrt{(\xcl)^2+(\ycl)^2}\le 1\}$.
For $r\rightarrow 1$ the metric tensor $\First$ tends to become
singular as one of the two tangent vectors increases indefinitely at
the boundary of the surface, leading to large spectral condition numbers $\Cond(\First)$.
Analogously, the condition number of $\First$ increases when the
frequency $k$ and the amplitude $a$ are increased.
Fig.~\ref{fig:surfaces} shows the three-dimensional plot of the
surface in the left panel, and the spatial distributions of
$\metrTensCv{11}$, $\metrTensCv{22}$, and $\sqrt{\DET{\First}}$,
respectively, in the next three columns.
The rows are relative to the case $r=2$, with the sphere ($a=0$) shown
in the top, while the trigonometric deformation of the sphere is shown
in the middle row for the case $a=0.5$ and $k=5$, and in the bottom
row for the case $a=2$ and $k=5$.
In this latter cases we note that $\First$ has a sinusoidal behavior
with $\metrTensCv{11}\approx 1$ in the regions where
$\metrTensCv{22}\approx 20(150)$, leading to $\Cond(\First)\approx
20(150)$.
Note that $\First$ does not enter the reaction term and has a small
effect in the advection term, as it amounts to a rotation and a
stretching of the advective field.
On the other hand, it has a large effect on the coercivity of the
diffusion bilinear form, and thus we concentrate on the latter.
The condition number of the VEM stiffness matrix $\mathcal{A}$ can be
bounded by~\cite{book:QuarteroniValli94,book:BrennerScott2008}:
\begin{equation*}
  \Cond(\mathcal{A})\le
  \max_{\sv\in\RefDomain}\left(\sqrt{\DET{\First(\sv)}}\,
    \Cond(\First^{-1}(\sv))
  \right)
  \Cond(\mathcal{A}_{\scriptscriptstyle{\textrm{LAP}}}),
\end{equation*}
where $\Cond(\mathcal{A}_{\scriptscriptstyle{\textrm{LAP}}})$ is the
condition number of the VEM stiffness matrix of Laplace equation.
In our case, we have that:
\begin{equation*}
  \sqrt{\DET{(\First)}}\;\Cond(\First^{-1})
  % =\frac{\max\{\metrTensCv{11},\metrTensCv{22}\}}{\min\{\metrTensCv{11},\metrTensCv{22}\}}
  % \sqrt{\max\{\metrTensCv{11},\metrTensCv{22}\}\min\{\metrTensCv{11},\metrTensCv{22}\}} 
  =\sqrt{\frac{\max\{\metrTensCv{11},\metrTensCv{22}\}^{3}}{\min\{\metrTensCv{11},\metrTensCv{22}\}}}.
\end{equation*}
Note that, $\sqrt{\DET{(\First)}}\;\Cond(\First^{-1})\ge 1$ is smooth
although possibly unbounded when $r\To 1$ (or $k$ and $a$ are large),
as mentioned above. 
We would like to remark that this is not a 
contraddiction of eq.~\eqref{eq:coercivity-continuity}, but rather a 
consequence of the fact that we use a single Monge parametrization.
The presence of the metric tensor in the equation always deteriorates,
possibly drastically, the condition number of the system matrix.

% stereographic projection
In Test Case~4, we consider Laplace equation (i.e., $\AdvVel=0$ and
$\MassCoef=0$) on $\SurfDomain=S^2$ and use two parametrizations, one
for the northern and one for the southern hemispheres, given by:
\begin{equation*}
  \MapU[N](\xcl,\ycl)
  =\left(\frac{2\xcl}{1+(\xcl)^2+(\ycl)^2},
    \frac{2\ycl}{1+(\xcl)^2+(\ycl)^2},
    \frac{1-(\xcl)^2-(\ycl)^2}{1+(\xcl)^2+(\ycl)^2}\right)
  =(\xcg,\ycg,\zcg)\,,
\end{equation*}
and:
\begin{equation*}
  \MapU[S](\xcl,\ycl)=\left(\frac{2\xcl}{1+(\xcl)^2+(\ycl)^2},
    \frac{2\ycl}{1+(\xcl)^2+(\ycl)^2},
    \frac{-1+(\xcl)^2+(\ycl)^2}{1+(\xcl)^2+(\ycl)^2}\right)
  =(\xcg,\ycg,\zcg)\,.
\end{equation*}
We proceed by discretizing the unit disk as reference domain once and
for all for both hemispheres, with a polygonal approximation of the
boundary. We then use the appropriate charts to express the VEM linear
and bilinear forms in the northern and southern cells.  We connect the
two domains together by means of a simple Jacobi domain decomposition
approach, and to avoid iterations we use the manufactured solution as
boundary condition at the domain interface.
Note that the use of curved edges would allow to solve the problem
without the need to decompose the computational domain, exploiting the
fact that the transition map for the two charts is readily available
and sufficiently regular.

In all the experiments, numerical errors are evaluated by defining a
manufactured solution $\Sol(\xcl,\ycl)\From\RefDomain\To\REAL$ and
calculating the resulting forcing function $\force(\xcl,\ycl)$ by
substitution into the original equation.
Using $\Sol(\xcl,\ycl)=\sin(2\pi\xcl)\sin(2\pi\ycl)$ and taking into account the contributions of the
metric $\First$, the general form
of $\force(\xcl,\ycl)$ can be evaluated as:
\begin{align*}
  \force(\xcl,\ycl)
  &=
  \sin(2\pi\ycl)\left(\frac{2\pi\ws_{(1)}\cos(2\pi\xcl)}{\sqrt{\metrTensCv{11}}}
  +\MassCoef\sin(2\pi\xcl)\right)+\\
  &+\frac{\pi\sin(2\pi\xcl)\cos(2\pi\ycl)\left(\metrTensCv{11}
    \DerParDS{\metrTensCv{22}}{\ycl}-\DerParDS{\metrTensCv{11}}{\ycl}\metrTensCv{22}\right)}{\metrTensCv{22}\DET{\First}}+\\
  &+\pi 
  \sin(2\pi\ycl)\left[\frac{
      \cos(2\pi\xcl)\left(\DerParDS{\metrTensCv{11}}{\xcl}
      \metrTensCv{22}-\metrTensCv{11}
      \DerParDS{\metrTensCv{22}}{\xcl}\right)}{\metrTensCv{11}\DET{\First}}+
    \frac{4\pi\sin(2\pi\xcl)(\metrTensCv{11}+\metrTensCv{22})}{\metrTensCv{22}\DET{\First}}\right]
  +\\
  &+\frac{2\pi\ws_{(2)}\sin(2\pi \xcl)\cos(2\pi\ycl)}{\sqrt{\metrTensCv{22}}},
\end{align*}
where $\hat{\AdvVel}=[w_{(1)},w_{(2)}]^{\T}$.

\subsection{Test Case 1}
\begin{figure}
  \centerline{
    \includegraphics[width=0.24\textwidth]{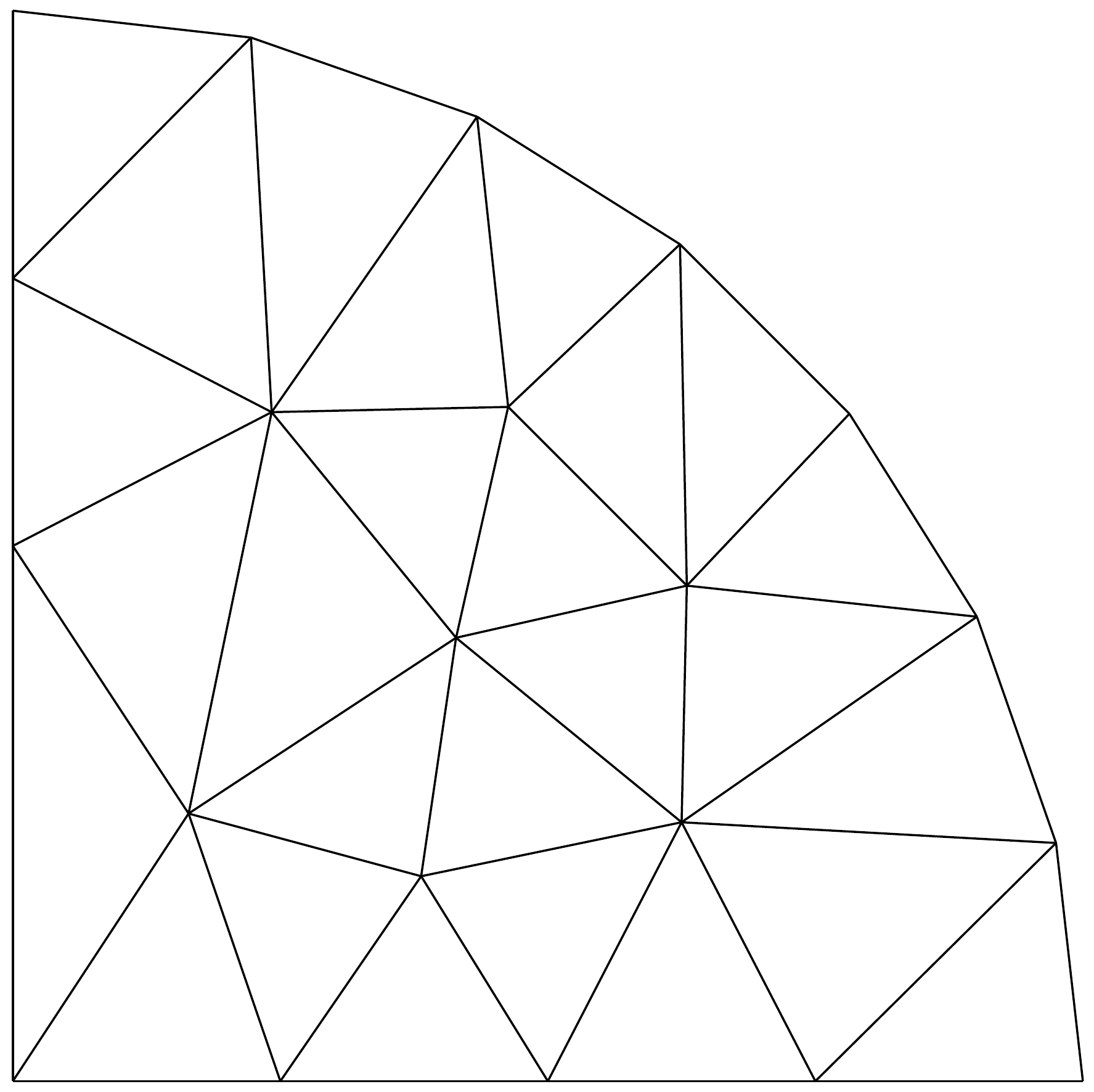}%%{triangle_mesh0}
    \includegraphics[width=0.24\textwidth]{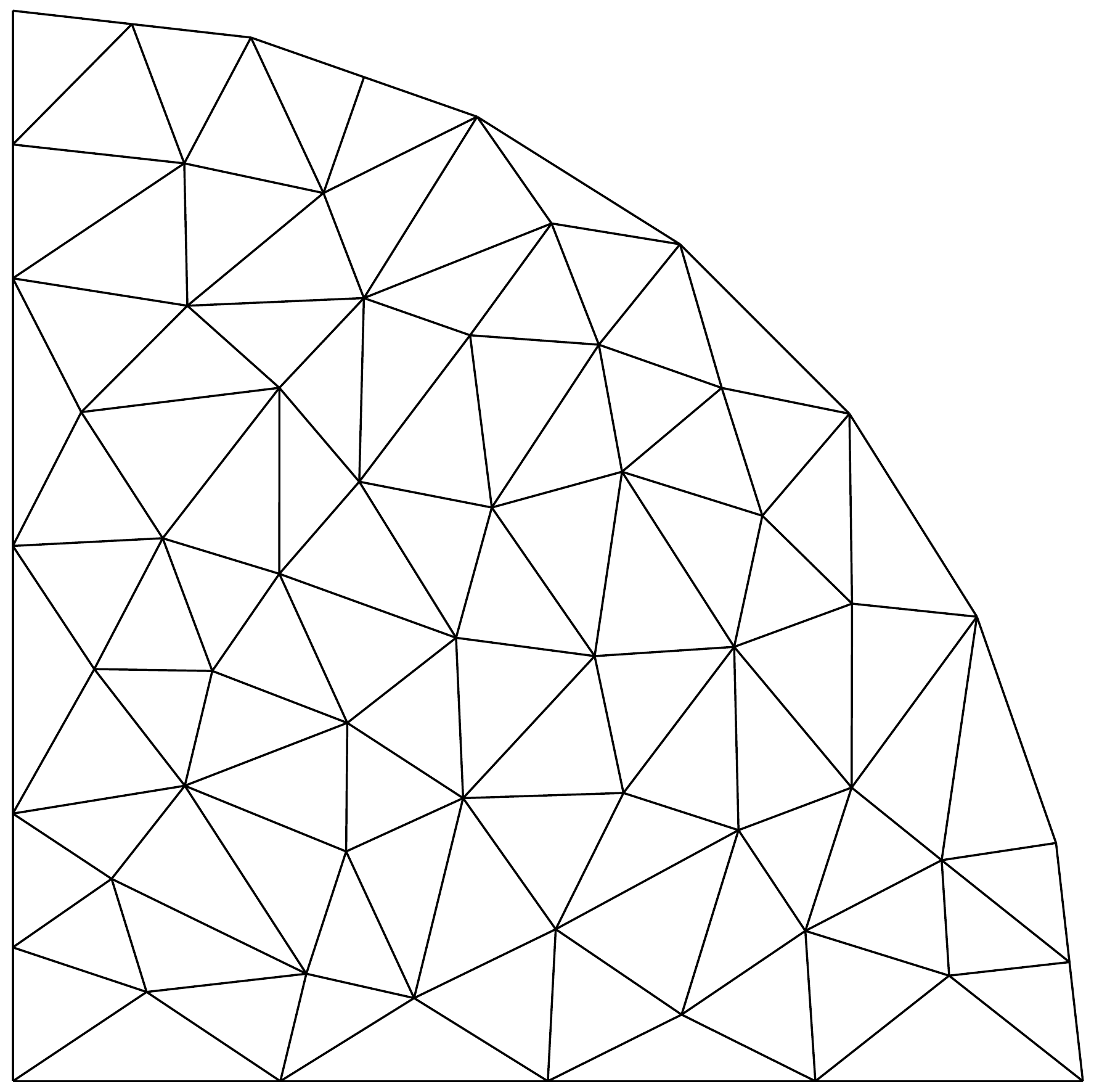}%%{triangle_mesh1}
    \hspace{2pt}
    \includegraphics[width=0.24\textwidth]{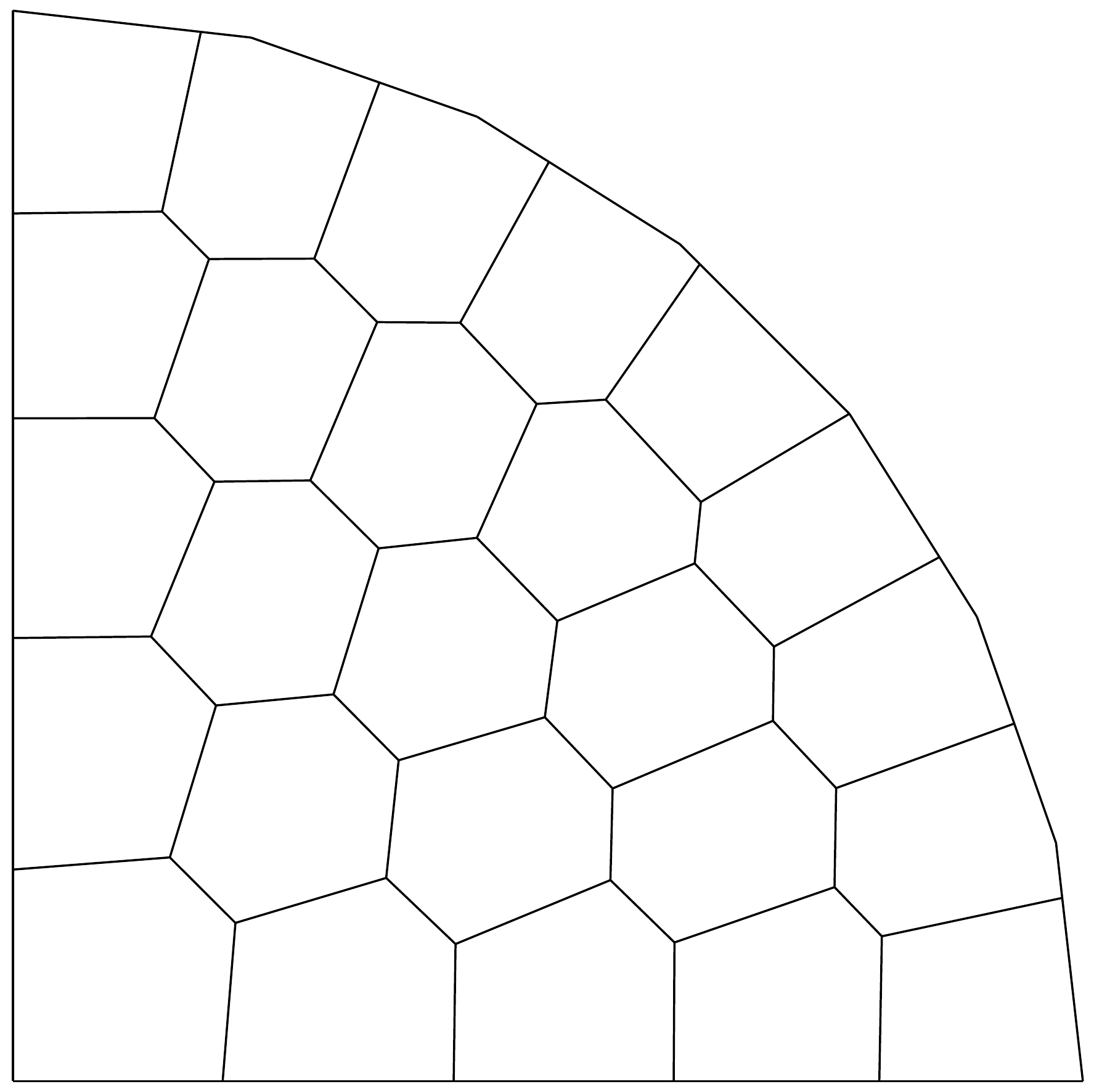}%%{polymesher_mesh0}
    \includegraphics[width=0.24\textwidth]{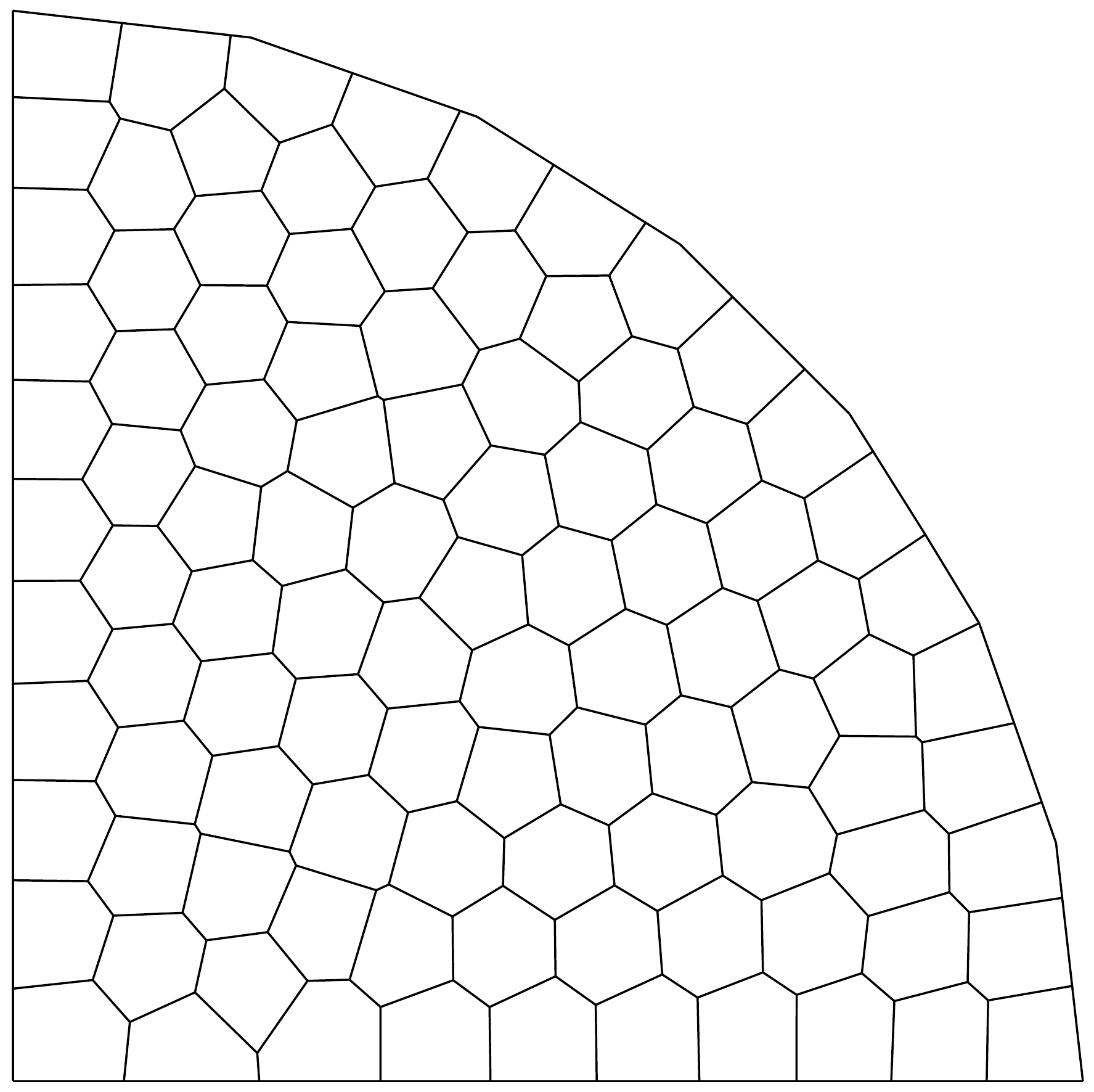}%%{polymesher_mesh1}
  }
  \caption{TC1: Level 0 and 1 triangulations (left panels) and polygonal
    meshes (right panels) of $\RefDomain$.}
  \label{fig:mesh0-test1}
\end{figure}

We use two families of polygonal meshes discretizing the domain (see
Fig.~\ref{fig:mesh0-test1}).
To avoid geometric error in the refinement process we uniformly
distribute 8 nodes on the curvilinear boundary of $\RefDomain$ and
approximate it with linear interpolation.
All the refined meshes are built on this geometry.
The meshes in the first family are constrained Delaunay triangulations
obtained using Triangle \cite{shewchuk96b,art:SHEWCHUK200221},
dividing by a factor 4 the area target of the elements at each
level.
The second family of meshes is obtained by means of PolyMesher
\cite{art:polymesher} by imposing approximately the same number of
elements of the triangulation at each corresponding level.
Note that the sides of the boundary elements may contain as extra node
one of the fixed vertices used to define the curvilinear boundary, and
are thus formed by more than one edge.

\begin{figure}
  \includegraphics[width=\textwidth]{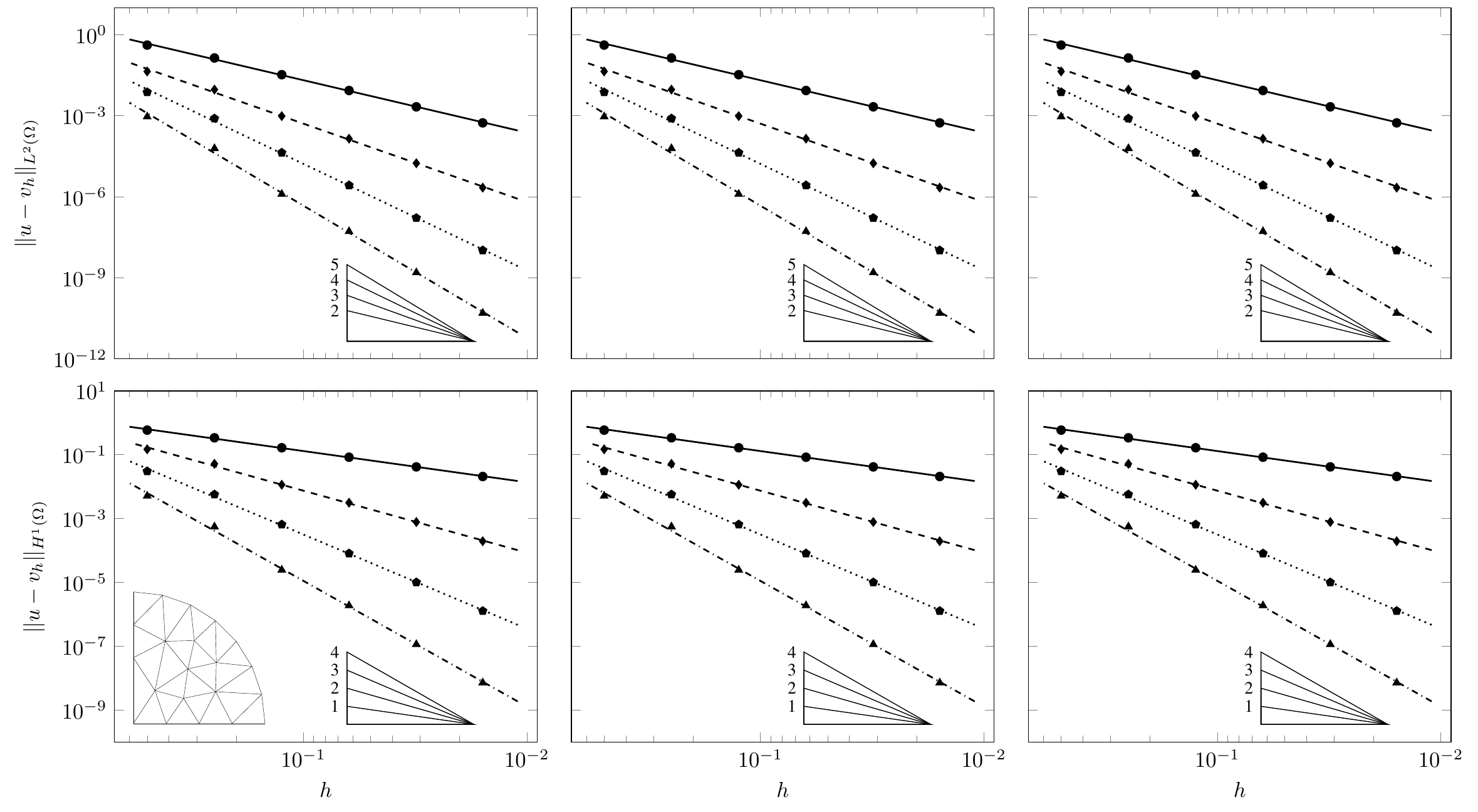}%%{tc1-conv_errL2-H1_tria_h}
  \caption{TC1: Convergence of $\LTWO$ (top) and $\HONE$ (bottom)
    errors vs $\hh$ on the triangulations.
    The convergence lines are obtained by approximating via
    least-squares all the point values.
    The different lines denote different polynomial orders from 1
    (solid line with circular data points) to 4 (dashed-dotted line with
    triangular data points).
    The optimal theoretical slope is represented by the lower right
    triangles.}
  \label{fig:tc1-tria}
\end{figure}
\begin{figure}
  \includegraphics[width=\textwidth]{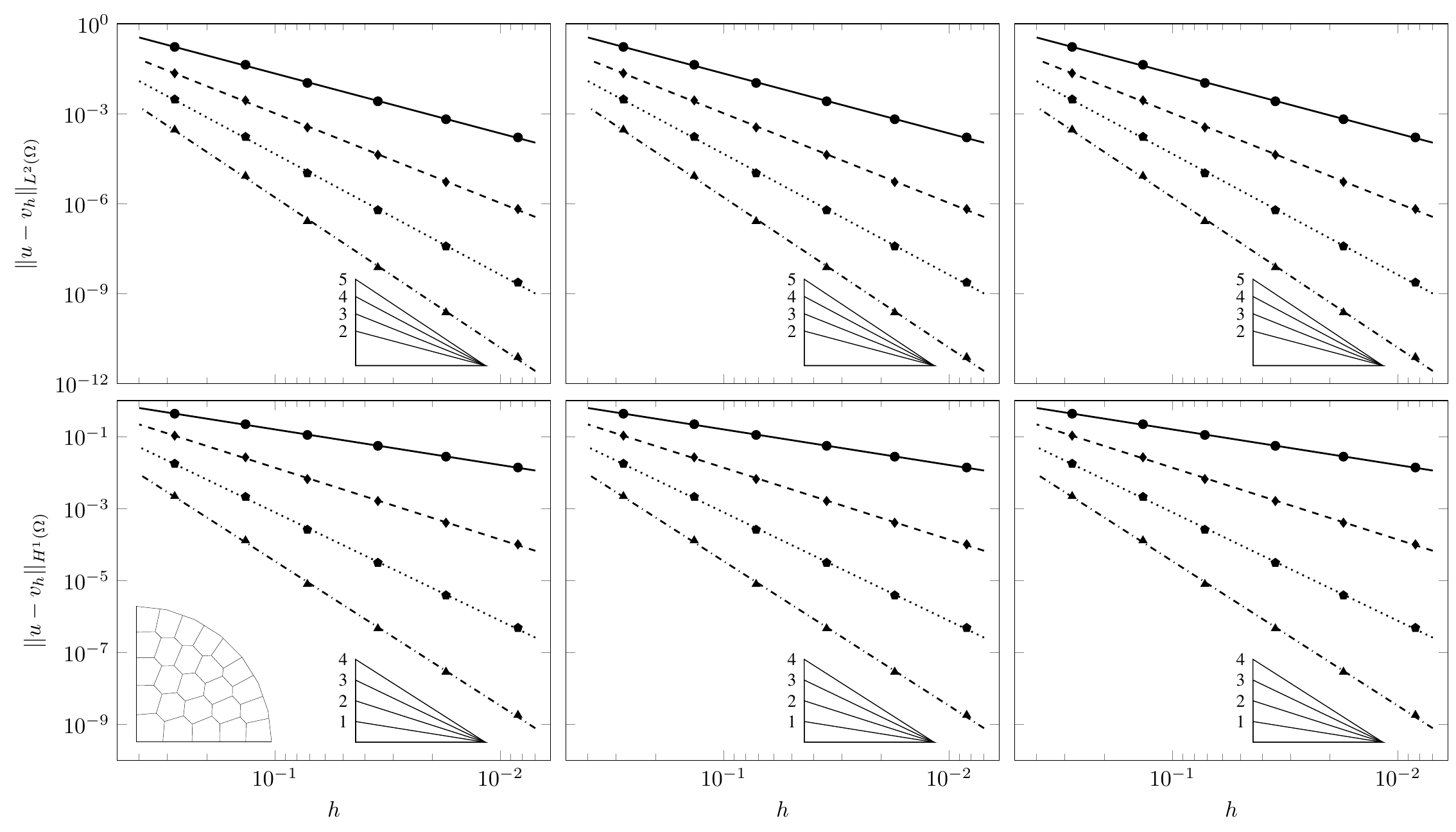}%%{tc1-conv_errL2-H1_poly_h}
  \caption{TC1: Convergence of $\LTWO$ (top) and $\HONE$ (bottom)
    errors vs $\hh$ on the polygonal meshes.
    The convergence lines are obtained by approximating via
    least-squares all the point values.
    The different lines denote different polynomial orders from 1
    (solid line with circular data points) to 4 (dashed-dotted line with
    triangular data points).
    The optimal theoretical slope is represented by the lower right
    triangles.}
  \label{fig:tc1-poly}
\end{figure}
Convergence is tested on four different grid refinement levels, in the
cases of $a=0$ and $r=1.1, 1.01, 1.001$.
Correspondingly, the results are reported in Fig.
\ref{fig:tc1-tria} and \ref{fig:tc1-poly} for the triangulations and
the polygonal meshes, respectively.
The experimental convergence rates are optimal for all the tested
polynomial orders, as can be seen from the figures and from the slopes
of the lines, which are obtained by approximating via least-squares
all the point values.
This confirms the theoretical expectations of the behavior of the VEM.

\subsection{Test Case 2}
\begin{figure} 
  \centerline{
  \includegraphics[width=0.24\textwidth]{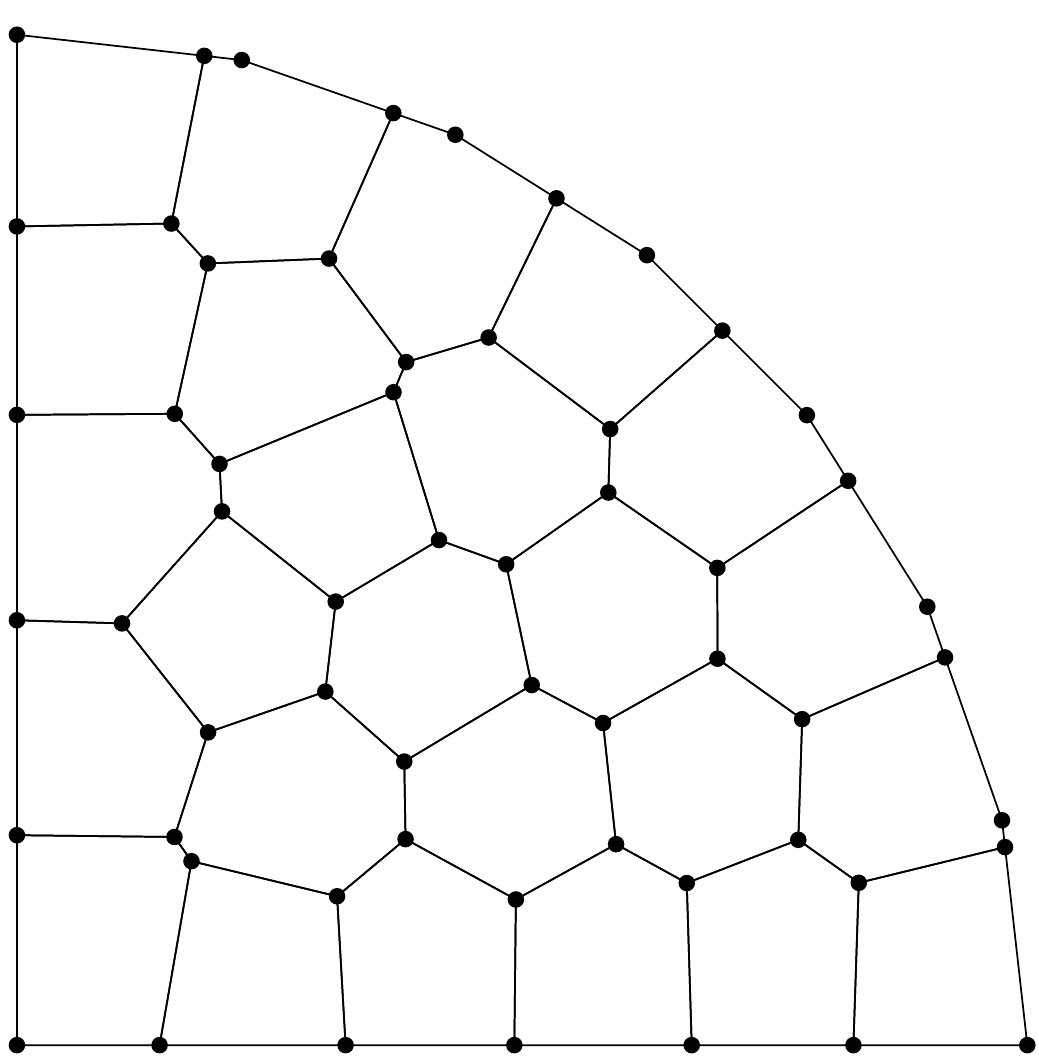} %% {polymesher_mesh0_tc2}
  \includegraphics[width=0.24\textwidth]{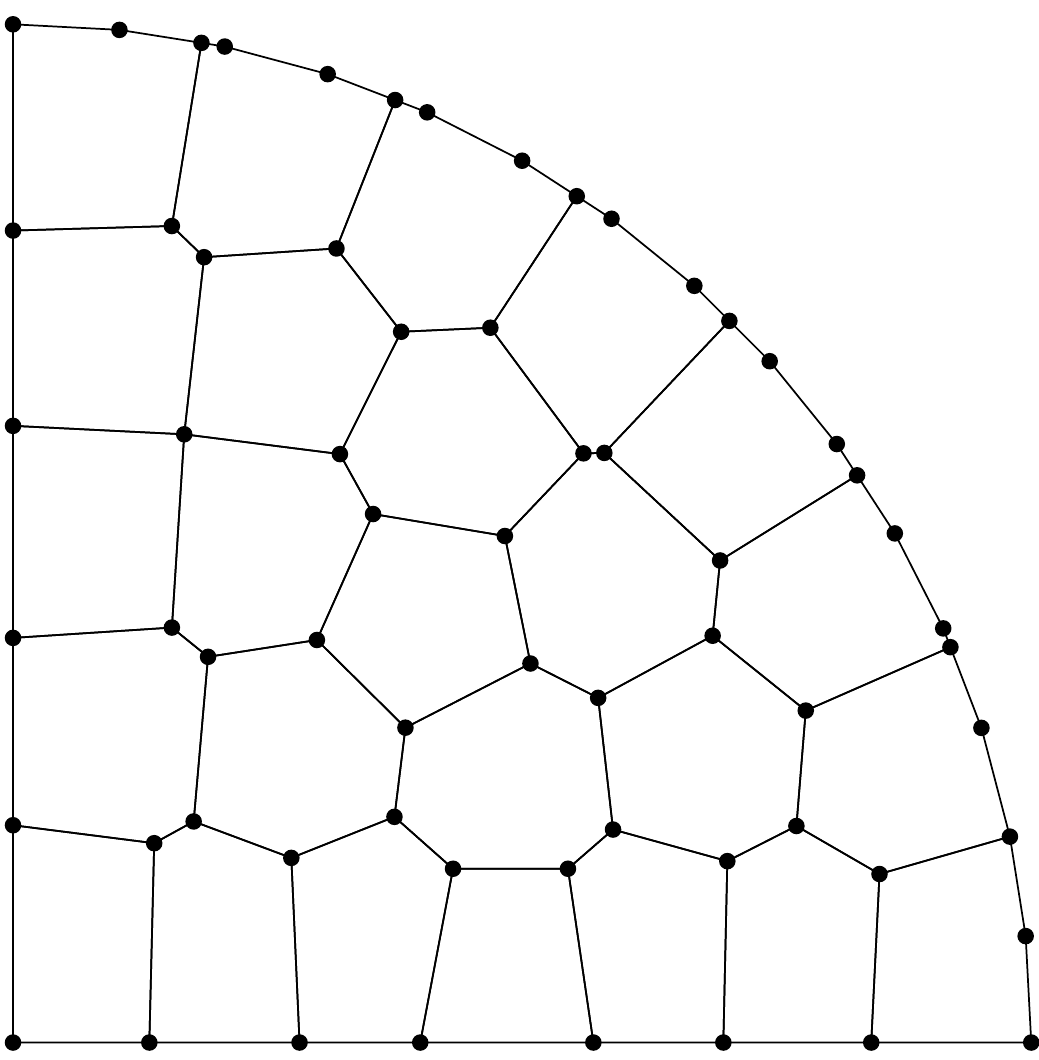} %% {polymesher_mesh1_tc2}
  \includegraphics[width=0.24\textwidth]{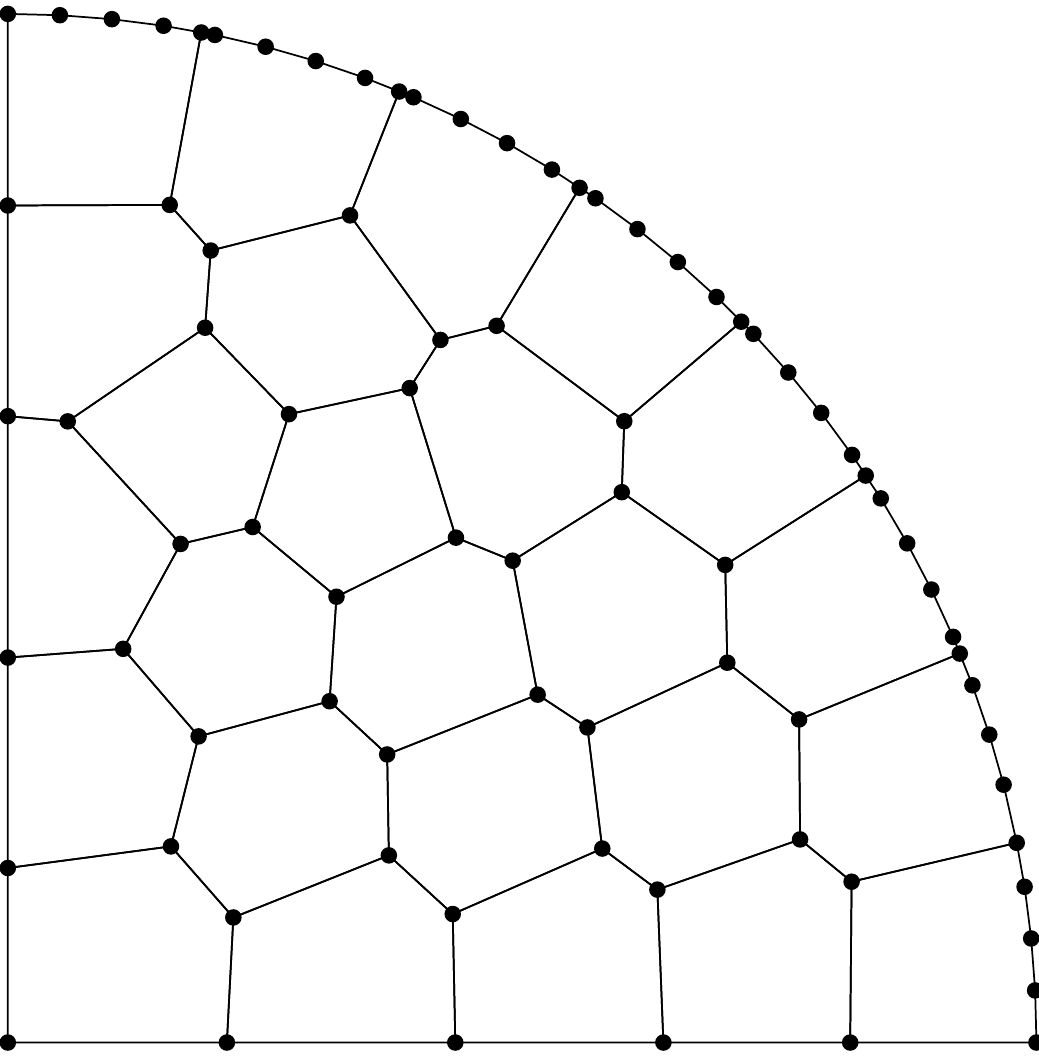} %% {polymesher_mesh2_tc2}
  }
  \caption{TC2: Level 0, 1, and 2 of polygonal meshes of $\RefDomain$
    with 25 cells and increased boundary resolution (8, 16, and 32
    nodes).}
  \label{fig:mesh012-test2}
\end{figure}

This test case is designed to verify the robustness of the scheme for
increasingly accurate approximations of the curvilinear boundary.
To this aim we look at the errors for a fixed mesh size and vary the
number of vertices used to discretize the curvilinear boundary, using
only the polygonal mesh.
Two sets of meshes are defined: one with approximately 25 elements and
the other with 100 elements.
In each set we consider 5 mesh levels characterized by different
approximations of the curvilinear boundary.
In the first level the boundary is discretized with 8 vertices, as in
the previous test case.
The subsequent levels are obtained by doubling each time the number of nodes located on the
curvilinear boundary to arrive at the final level with 128
vertices (see, e.g., Fig.~\ref{fig:mesh012-test2}).
Since the size of the cells remains approximately the same, the
boundary sides of the boundary elements are formed by an increasing
number of straight edges.
While the mesh levels approximate the curvilinear boundary with
increasing accuracy, the fact that the length of the edges of the
boundary elements becomes unbalanced may lead to increased errors in
the VEM solution.
However, we know from the literature~\cite{Brenner:2017:SEV}
that such unbalance may only affect the constant that appears in the 
error estimates.
Such constant is increased by a factor proportional to the square 
root of $\log(1+\ell)$, $\ell$ being the ratio between the maximum 
and the minimum edge length.
Hence, in our experiments we expect the errors to remain approximately
constant as we refine the boundary.

\begin{figure}
  \includegraphics[width=\textwidth]{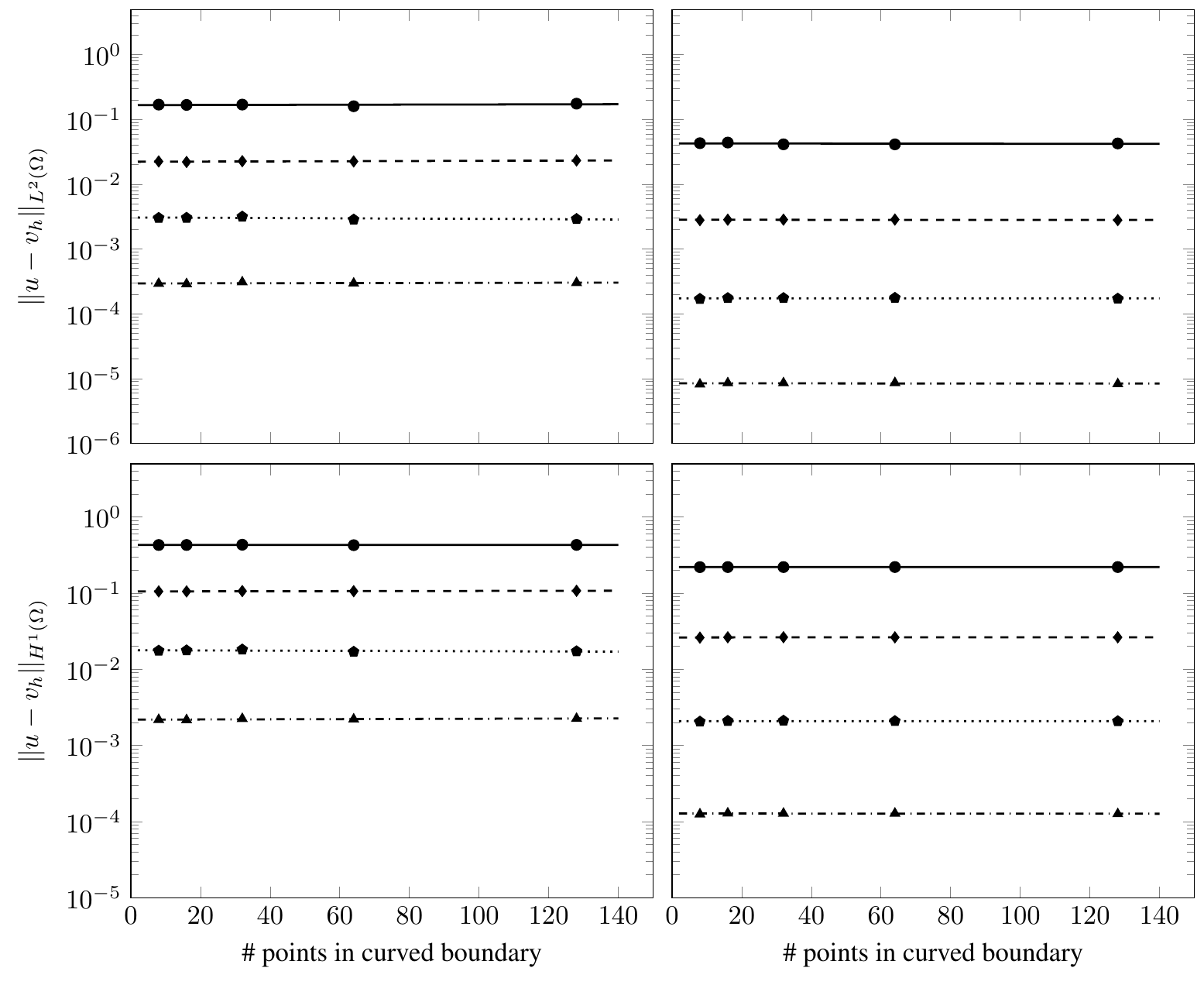} %% {tc2-errL2-H1_poly_h}
  \caption{TC2: $\LTWO$ (top) and $\HONE$ (bottom) errors vs number of
    points on the curvilinear boundary on the polygonal meshes
    (average of 25 cells, left panels; average of 100 cells, right
    panels).
    The different lines denote different polynomial orders from 1
    (solid line with circles) to 4 (dashed-dotted line with
    triangles). }
  \label{fig:tc2-poly}
\end{figure}

The results of the simulations are shown in Fig.~\ref{fig:tc2-poly},
where we report the $\LTWO$ and $\HONE$ errors (top and bottom,
respectively) with respect to the manufactured solutions as a function
of the number of points discretizing the curved boundary for the two
set of meshes (left and right).
It is evident that the proposed scheme is robust with respect the
increasing unbalance of the edge lengths as the errors for each
polynomial order remain approximately constant.

\subsection{Test Case 3}

This test case is aimed at studying the robustness of the scheme to
spatial variability and strong anisotropy of the diffusion bilinear
form with spatially variable anisotropy ratios.
We verify convergence of the proposed VEM scheme by solving our
equations on the same families of meshes of Test Case 1 with $r=2$,
$k=5$ and two different values for $a$, $a=0.5$ and $a=2$.
Fig.~\ref{fig:surfaces}, second and third rows, shows the surface
and the spatial distribution of $\metrTensCv{ii}$ for $k=5$ and
$a=0.5$ and $2$, respectively.
The spatial variability of the anisotropy ratio (ratio between
$\metrTensCv{11}$ and $\metrTensCv{22}$ or its inverse) varies between
1 and 150 for $a=2$ corresponding to large basis vectors for the
tangent plane.
This test case challenges the ability of the discretization scheme to
handle large and spatially varying anisotropy ratios.

\begin{figure}
  \includegraphics[width=\textwidth]{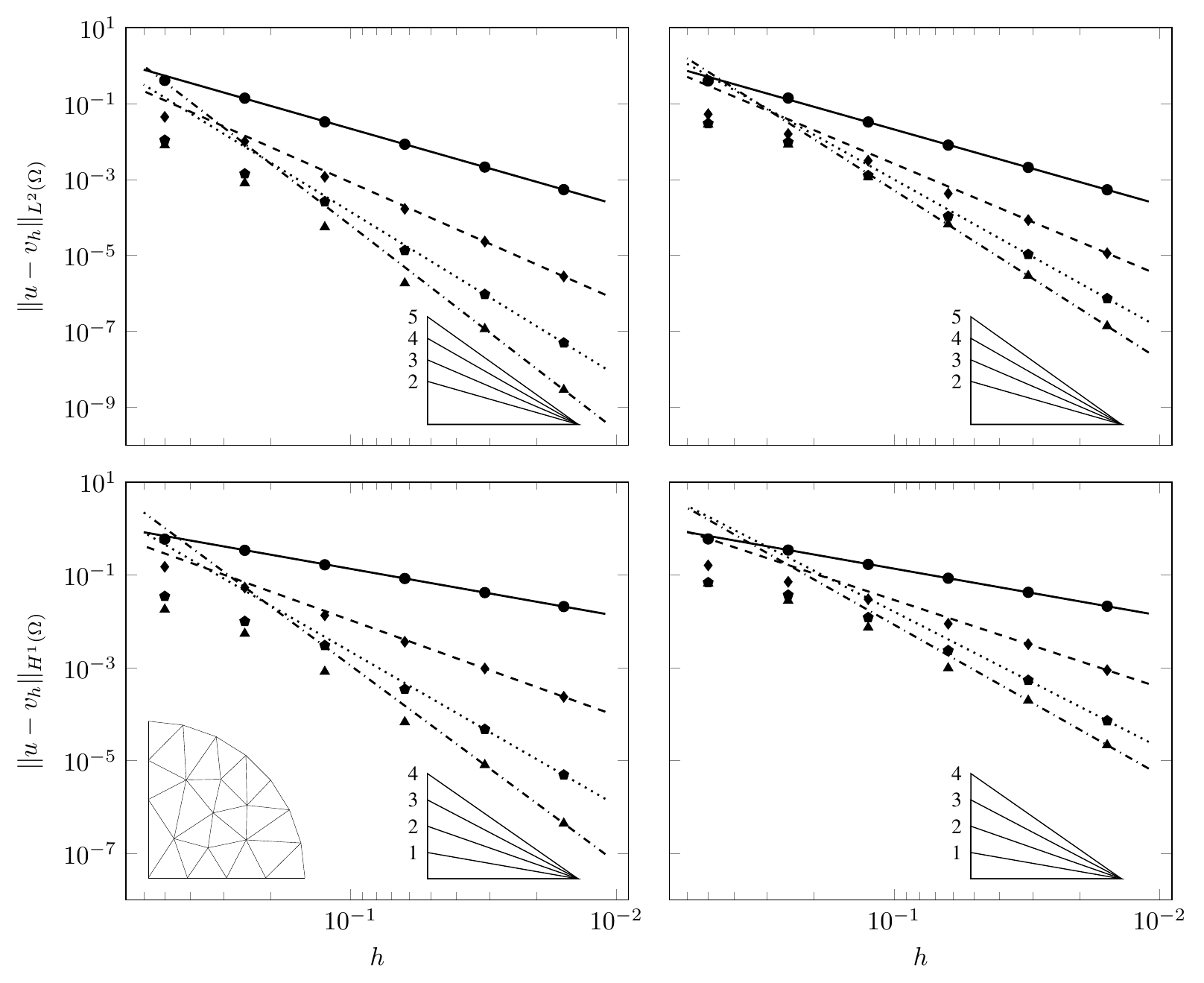} %% {tc3-conv_errL2-H1_tria_h}
  \caption{TC3: Convergence of $\LTWO$ (top) and $\HONE$ (bottom)
    errors vs $\hh$ on the triangulations ($a=0.5$, left panels;
    $a=2$, right panels).
    The convergence lines are obtained by approximating via
    least-squares only the last two point values.
    The different lines denote different polynomial orders from 1
    (solid line with circles) to 4 (dashed-dotted line with triangles).
    The optimal theoretical slope is represented by the lower right
    triangles. }
  \label{fig:tc3-tria}
\end{figure}

\begin{figure}
  \includegraphics[width=\textwidth]{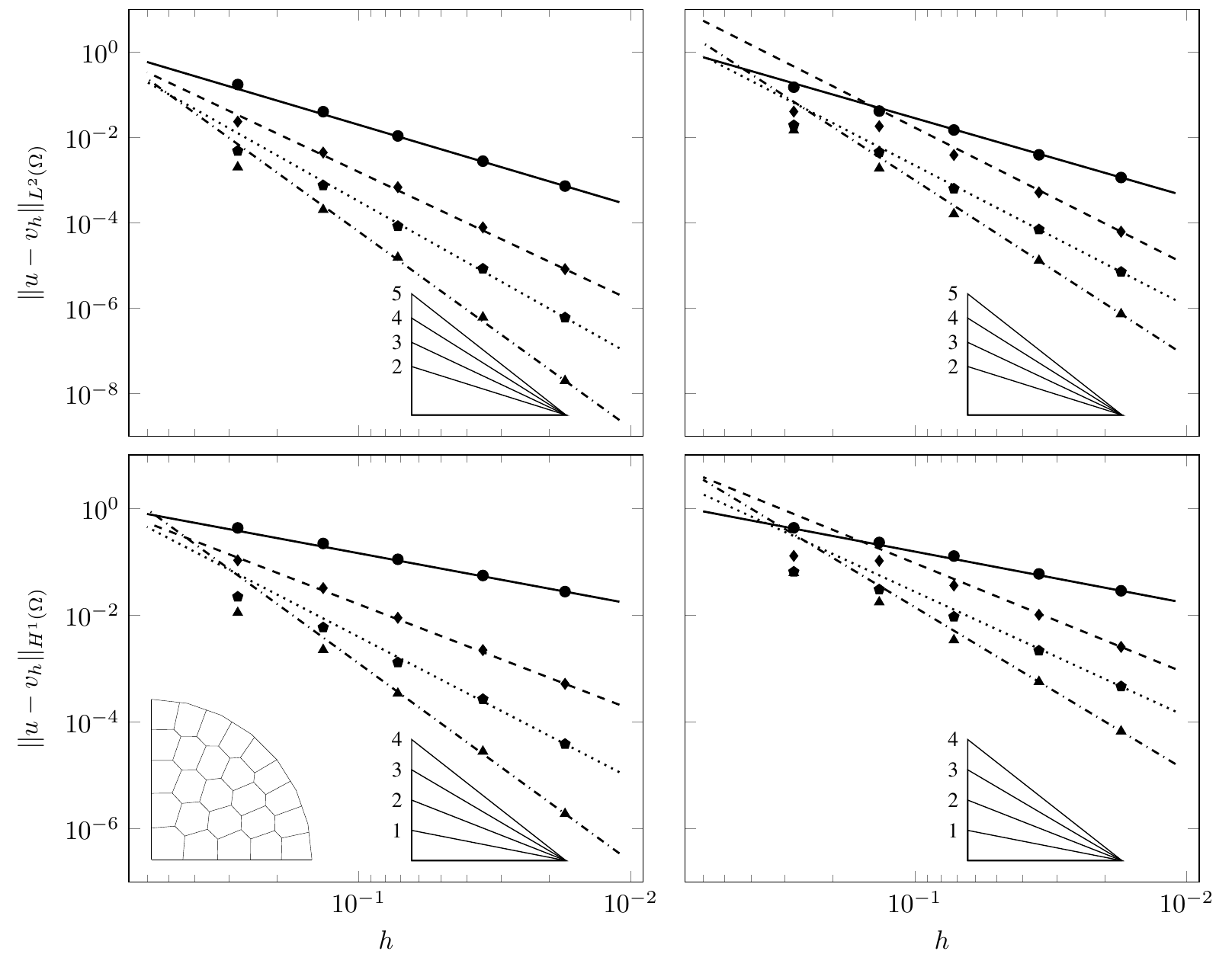} %% {tc3-conv_errL2-H1_poly_h}
  \caption{TC3: Convergence of $\LTWO$ (top) and $\HONE$ (bottom)
    errors vs $\hh$ on the polygonal meshes ($a=0.5$, left panels;
    $a=2$, right panels).
    The convergence lines are obtained by approximating via
    least-squares only the last two point values.
    The different lines denote different polynomial orders from 1
    (solid line with circles) to 4 (dashed-dotted line with
    triangles).
    The optimal theoretical slope is represented by the lower right
    triangles.}
  \label{fig:tc3-poly}
\end{figure}

Fig.~\ref{fig:tc3-tria} and~\ref{fig:tc3-poly} show the numerical
convergence of the $\LTWO$ and $\HONE$ norms of the error as a
function of $\hh$.
We note that asymptotic behavior of the error is reached as soon as
the mesh size is able to resolve the spatial scales of variation of
the metric tensor.
For this reason the convergence lines in the figures are obtained by
interpolation of the last two point values for each polynomial order.
Pre-asymptotic convergence is more evident for the higher order
polynomials.
We attribute this behavior to the smoothing effects of a lower order
interpolation.
Indeed, the first two or three point values for all polynomial orders
display at least the same convergence order of linear polynomials.

\begin{table}
  \begin{center}
    \begin{tabular}{|c|cc|cc|cc|cc|}
      \hline
      & \multicolumn{2}{c|}{$\PC{1}$}  & \multicolumn{2}{c|}{$\PC{2}$}  & \multicolumn{2}{c|}{$\PC{3}$}  & \multicolumn{2}{c|}{$\PC{4}$}  \\
      $\Lev$  & $\norm{\us-\vsh}{\LTWO}$ & $\Eoc$ & $\norm{\us-\vsh}{\LTWO}$ & $\Eoc$ & $\norm{\us-\vsh}{\LTWO}$ & $\Eoc$ & $\norm{\us-\vsh}{\LTWO}$ & $\Eoc$\\
      \hline
      \multicolumn{9}{|c|}{$a=0.5$}\\
      \hline
      $0$  & $4.09\,10^{-1}$ & $--$    & $4.47\,10^{-2}$ & $--$   & $1.09\,10^{-2}$ & $--$    & $8.03\,10^{-3}$ & $--$\\
      $1$  & $1.39\,10^{-1}$ & $1.556$ & $1.00\,10^{-2}$ & $2.154$& $1.42\,10^{-3}$ & $2.943$ & $8.08\,10^{-4}$ & $3.313$\\
      $2$  & $3.30\,10^{-2}$ & $2.076$ & $1.18\,10^{-3}$ & $3.086$& $2.62\,10^{-4}$ & $2.436$ & $5.60\,10^{-5}$ & $3.851$\\
      $3$  & $8.50\,10^{-3}$ & $1.956$ & $1.69\,10^{-4}$ & $2.809$& $1.35\,10^{-5}$ & $4.276$ & $1.86\,10^{-6}$ & $4.911$\\
      $4$  & $2.13\,10^{-3}$ & $1.995$ & $2.32\,10^{-5}$ & $2.861$& $9.47\,10^{-7}$ & $3.837$ & $1.16\,10^{-7}$ & $4.000$\\  
      $5$  & $5.43\,10^{-4}$ & $1.998$ & $2.80\,10^{-6}$ & $3.088$& $4.99\,10^{-8}$ & $4.303$ & $2.89\,10^{-9}$ & $5.399$\\
      \hline
      \multicolumn{9}{|c|}{$a=2.0$}\\
      \hline
      $0$ & $3.95\,10^{-1}$ & $--$    & $5.26\,10^{-2}$ & $--$    & $2.99\,10^{-2}$ & $--$   & $2.84\,10^{-2}$ & $--$\\ 
      $1$ & $1.40\,10^{-1}$ & $1.495$ & $1.59\,10^{-2}$ & $1.730$ & $9.61\,10^{-3}$ & $1.639$& $8.40\,10^{-3}$ & $1.755$\\
      $2$ & $3.28\,10^{-2}$ & $2.094$ & $3.14\,10^{-3}$ & $2.335$ & $1.28\,10^{-3}$ & $2.913$& $1.16\,10^{-3}$ & $2.852$\\
      $3$ & $8.07\,10^{-3}$ & $2.023$ & $4.27\,10^{-4}$ & $2.878$ & $1.07\,10^{-4}$ & $3.571$& $6.66\,10^{-5}$ & $4.125$\\
      $4$ & $2.09\,10^{-3}$ & $1.951$ & $8.55\,10^{-5}$ & $2.321$ & $1.07\,10^{-5}$ & $3.324$& $2.92\,10^{-6}$ & $4.512$\\
      $5$ & $5.38\,10^{-4}$ & $1.981$ & $1.14\,10^{-5}$ & $2.940$ & $7.38\,10^{-7}$ & $3.910$& $1.38\,10^{-7}$ & $4.462$\\
      \hline
    \end{tabular}
  \end{center}
  \caption{TC3: experimental errors and convergence rates for the
    triangular mesh set.}
  \label{tab:tc3-conv-tria}
\end{table}

\begin{table}
  \begin{center}
    \begin{tabular}{|c|cc|cc|cc|cc|}
      \hline
      & \multicolumn{2}{c|}{$\PC{1}$}  & \multicolumn{2}{c|}{$\PC{2}$}  & \multicolumn{2}{c|}{$\PC{3}$}  & \multicolumn{2}{c|}{$\PC{4}$}  \\
      $\Lev$  & $\norm{\us-\vsh}{\LTWO}$ & $\Eoc$ & $\norm{\us-\vsh}{\LTWO}$ & $\Eoc$ & $\norm{\us-\vsh}{\LTWO}$ & $\Eoc$ & $\norm{\us-\vsh}{\LTWO}$ & $\Eoc$\\
      \hline
      \multicolumn{9}{|c|}{$a=0.5$}\\
      \hline
      $0$ & $1.73\,10^{-1}$ & $--$    & $2.34\,10^{-2}$ & $--$    & $4.84\,10^{-3}$ & $--$    & $2.01\,10^{-3}$ & $--$\\   
      $1$ & $4.01\,10^{-2}$ & $2.023$ & $4.42\,10^{-3}$ & $2.305$ & $7.54\,10^{-4}$ & $2.568$ & $2.02\,10^{-4}$ & $3.177$\\
      $2$ & $1.08\,10^{-2}$ & $2.069$ & $6.83\,10^{-4}$ & $2.953$ & $8.36\,10^{-5}$ & $3.479$ & $1.54\,10^{-5}$ & $4.072$\\
      $3$ & $2.78\,10^{-3}$ & $1.888$ & $7.81\,10^{-5}$ & $3.012$ & $8.40\,10^{-6}$ & $3.191$ & $6.15\,10^{-7}$ & $4.471$\\
      $4$ & $7.26\,10^{-4}$ & $1.929$ & $8.23\,10^{-6}$ & $3.235$ & $6.05\,10^{-7}$ & $3.782$ & $2.00\,10^{-8}$ & $4.924$\\   
      $5$ & $1.81\,10^{-4}$ & $1.890$ & $9.07\,10^{-7}$ & $3.000$ & $4.33\,10^{-8}$ & $3.587$ & $6.78\,10^{-10}$ & $4.604$\\
      \hline
      \multicolumn{9}{|c|}{$a=2.0$}\\
      \hline
      $0$ & $1.50\,10^{-1}$ & $--$    & $4.03\,10^{-2}$ & $--$    & $1.92\,10^{-2}$ & $--$    & $1.46\,10^{-2}$ & $--$\\ 
      $1$ & $4.14\,10^{-2}$ & $1.780$ & $1.83\,10^{-2}$ & $1.092$ & $4.52\,10^{-3}$ & $2.002$ & $1.88\,10^{-3}$ & $2.836$\\
      $2$ & $1.49\,10^{-2}$ & $1.613$ & $3.87\,10^{-3}$ & $2.455$ & $6.31\,10^{-4}$ & $3.115$ & $1.61\,10^{-4}$ & $3.891$\\
      $3$ & $3.92\,10^{-3}$ & $1.859$ & $5.16\,10^{-4}$ & $2.799$ & $6.99\,10^{-5}$ & $3.056$ & $1.31\,10^{-5}$ & $3.486$\\
      $4$ & $1.15\,10^{-3}$ & $1.759$ & $6.20\,10^{-5}$ & $3.046$ & $7.07\,10^{-6}$ & $3.294$ & $7.32\,10^{-7}$ & $4.142$\\
      $5$ & $2.99\,10^{-4}$ & $1.832$ & $5.83\,10^{-6}$ & $3.214$ & $6.30\,10^{-7}$ & $3.290$ & $3.52\,10^{-8}$ & $4.130$\\
      \hline
    \end{tabular}
  \end{center}
  \caption{TC3: experimental errors and convergence rates for the
    polygonal mesh set.}
  \label{tab:tc3-conv-poly}
\end{table}

The high degree of anisotropy of this test case causes a loss of
convergence in the higher polynomial orders.
This behavior is more evident for the polygonal meshes.
To better quantify this convergence loss,
Tables~\ref{tab:tc3-conv-tria} and~\ref{tab:tc3-conv-poly} report the
convergence order sequence for the triangulations and polygonal
meshes, respectively.
We note that, for first and second order polynomials, almost optimal
order of convergence is reached for both mesh sets and both $a=0.5$
and $a=2$.
In the higher orders, a clear loss of at least half an order is
evident.
This can be attributed to difficulties in resolving the large
anisotropy ratios that are typical of this test case.
The fact that this occurs only for the higher orders is due to the
ill-conditioning of the resulting linear system that controls the
ratio between the norms of the residuals and the errors in the linear
system solver.

\subsection{Test Case 4}
\begin{figure}
  \centerline{
    \includegraphics[width=0.6\textwidth]{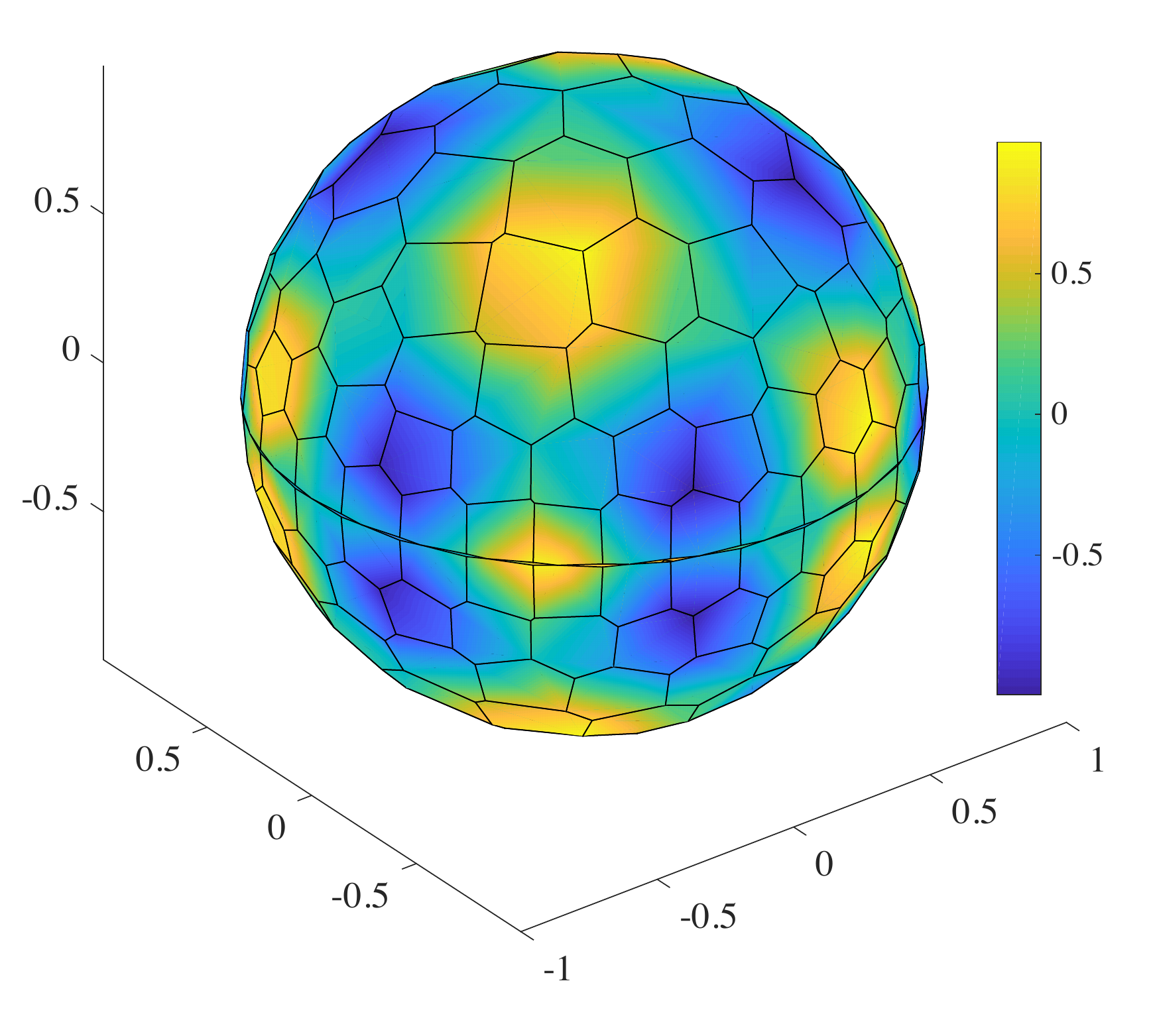} %%{sphere}
    }
  \caption{TC4: Numerical solution (nodal dofs only) on the coarsest
    mesh with linear interpolation from the nodal values.}
  \label{fig:sphere}
\end{figure}
 
This test case shows the convergence properties of the proposed
framework in a multiple charts setting.
We discretize the full unit disk with 5 polygonal mesh levels
following the strategy reported in Test Case~1, so that the level
$\Lev=0$ is characterized by 100 cells and the last one ($\Lev=4$) by
25600 cells.
The VEM solution, reconstructed on $\SurfDomain=S^2$ using nodal only
degrees of freedom, is shown in Fig.~\ref{fig:sphere} for the
coarsest mesh.

\begin{figure}
\includegraphics[width=\textwidth]{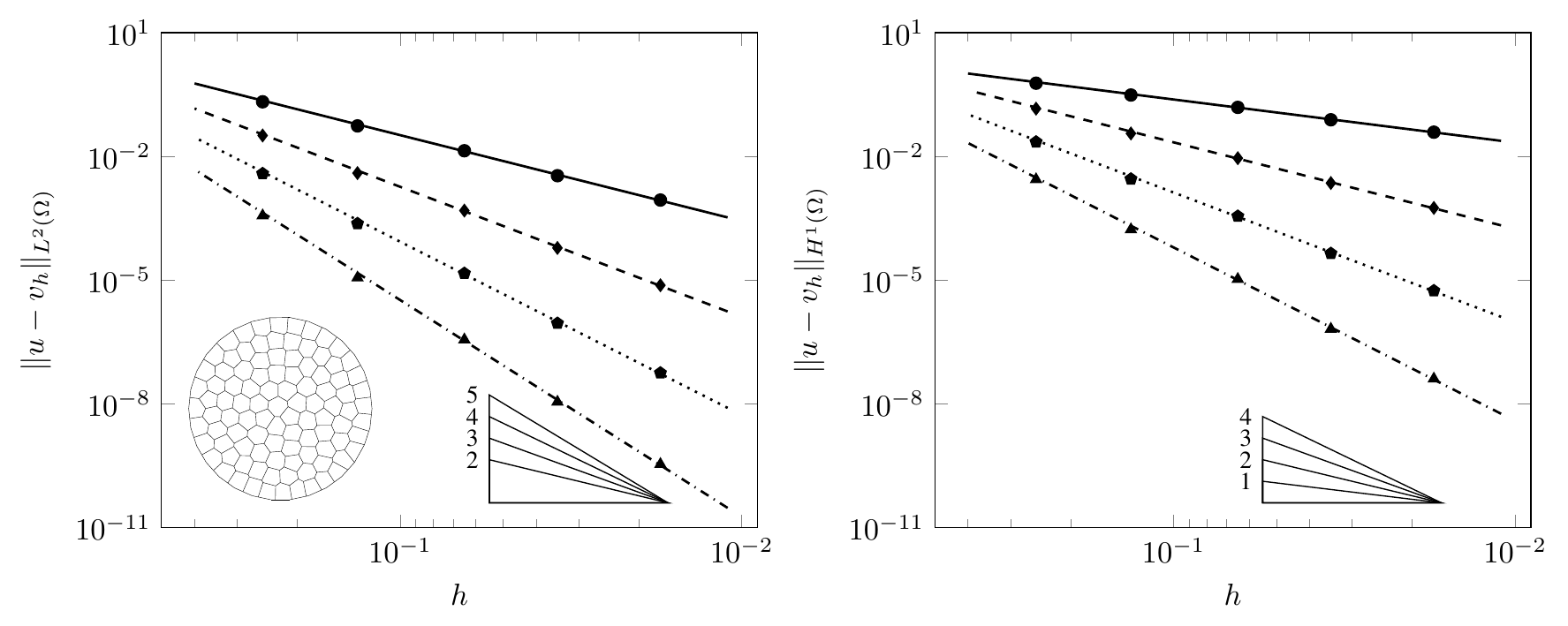} %% {tc4-conv_errL2-H1_poly_h}
  \caption{TC4: Convergence of $\LTWO$ (left) and $\HONE$ (right)
    errors vs $\hh$ on the polygonal meshes.
    The convergence lines are obtained by approximating via
    least-squares only the last three point values.
    The different lines denote different polynomial orders from 1
    (solid line with circles) to 4 (dashed-dotted line with
    triangles).
    The optimal theoretical slope is represented by the lower right
triangles. }
  \label{fig:tc4-poly}
\end{figure}

\begin{table}
  \begin{center}
    \begin{tabular}{|c|cc|cc|cc|cc|} \hline &
\multicolumn{2}{c|}{$\PC{1}$} & \multicolumn{2}{c|}{$\PC{2}$} &
\multicolumn{2}{c|}{$\PC{3}$} & \multicolumn{2}{c|}{$\PC{4}$} \\
$\Lev$ & $\norm{\us-\vsh}{\LTWO}$ & $\Eoc$ & $\norm{\us-\vsh}{\LTWO}$
& $\Eoc$ & $\norm{\us-\vsh}{\LTWO}$ & $\Eoc$ &
$\norm{\us-\vsh}{\LTWO}$ & $\Eoc$\\ \hline
      $0$ &  $2.12 \,10^{-1}$ & $--$ & $3.23 \,10^{-2}$  & $--$  & $3.85 \,10^{-3}$ & $--$ & $3.68 \,10^{-4}$ & $--$\\            
      $1$ &  $5.52 \,10^{-2}$ & $ 2.104$ & $ 3.96 \,10^{-3}$ & $ 3.284$ & $ 2.36 \,10^{-4}$ & $ 4.367$ & $ 1.1510^{-5} $ & $ 5.424$\\
      $2$ &  $ 1.38 \,10^{-2}$ & $ 1.924$ & $ 4.88 \,10^{-4}$ & $ 2.905$ & $ 1.4710^{-5}$ & $ 3.857$ & $ 3.61 \,10^{-7} $ & $ 4.801$\\
      $3$ &  $ 3.42 \,10^{-3}$ & $ 2.221$ & $ 6.0510^{-5}$ & $ 3.326$ & $ 9.02 \,10^{-7}$ & $ 4.443$ & $ 1.11 \,10^{-8} $ & $ 5.543$\\
      $4$ &  $ 8.78 \,10^{-4}$ &  $ 1.960$ & $ 7.53 \,10^{-6}$ & $ 3.002$ & $ 5.63 \,10^{-8}$ & $ 3.995$ & $ 3.44 \,10^{-10} $ & $ 5.009$\\
\hline      
    \end{tabular}
  \end{center}
  \caption{TC4: experimental errors and convergence rates for the
polygonal mesh set.}
  \label{tab:tc4-conv}
\end{table}

The experimental convergence on these mesh levels is reported for the
$\LTWO$ and $\HONE$ norms of the error as a function of $\hh$ in
Fig.~\ref{fig:tc4-poly} and in Table~\ref{tab:tc4-conv}.
The numerical results show that the proposed approach is functioning
as expected and that the use of two different charts doe not influence
the optimal convergence of the scheme.
Obviously, this is a very favorable case as the stereographical
projection produces charts and, if necessary, transition maps that are
sufficiently smooth to allow high order. In the future it will be
important to study how to derive charts and transition maps with
specified regularity for different surfaces, possibly starting from
the work of~\cite{art:Lindblom2016}.

%% SECTION 5
\section{Conclusions}
\label{sec:conclusions}

We have developed an arbitrary-order virtual element method for the
discretization of elliptic surface PDEs.
The approach employs a local parametrization of the surface to
properly re-define the PDE on the local chart. This allows the
straight-forward definition of a two-dimensional VEM discretization at
all polynomial orders, overcoming the difficult task of the consistent
approximation of the surface and of the distance function of its
tubular neighborhood.
The drawback of the approach is that the geometrically intrinsic form
of the PDE contains the metric information, which may be strongly
non-isotropic and highly variable in space, depending the regularity
of the surface.
The choice of the VEM scheme is motivated by the need to ensure
robustness and high order of convergence for these anisotropic and
spatially variable coefficients.

The developed scheme has been tested on several numerical examples
showing varying degrees of regularity.
In fact, optimal orders of convergence up to 5 has been reached for
surfaces with relatively small curvatures. Only when curvatures and
metric information become extremely large loss of convergence is
noticed.
This loss of convergence is related to the presence of strongly
anisotropic diffusion tensors and strongly aligned advective fields
due to the behavior of the metric tensor.  Handling strong anisotropy
is still a major challenge in the numerical solution of PDEs by the
virtual element method, and is left for future research.
This difficulty can be relaxed also by employing multiple charts that
decrease the anisotropic characteristic of the metric tensor. However,
proper regularity of the transition maps between the different charts
must be ensured to achieve full order convergence.
To verify the ability of our formulation to work with multiple charts
we tested the proposed scheme on the full sphere by employing two
charts arising from the stereographical projection.
Future work will be addressed to define for general surfaces
appropriate multiple charts with regular transition maps starting
from the work of~\cite{art:Lindblom2016}.

One of the major advantages of the developed VEM formulation is that
can be used efficiently to minimize geometric errors of curvilinear
boundaries.  We have tested our approach on an hemispherical surface
discretized by a fixed number of polygonal cells, where the boundary
edges formed by an increasing number of nodes. The resulting errors
were independent of the edge discretization, showing the robustness of
the VEM scheme in this situation.

%-------------------- acknowledgements --------------------------------

\begin{acknowledgements}
  EB and MP have been partially supported by project ``HYDROSEM''
  founded by Fondazione Cassa di Risparmio di Padova e Rovigo.
  GM has been partially supported by the ERC Project CHANGE, which has
  received funding from the European Research Council (ERC) under the
  European Unions Horizon 2020 research and innovation programme
  (grant agreement No 694515).
\end{acknowledgements}

%-------------------- bibliography -----------------------------------

%% \bibliographystyle{plain}
%% \bibliography{strings,new_add,vem}

\end{document}